\documentclass[draftclsnofoot,onecolumn,journal,romanappendices]{IEEEtran}

\usepackage{setspace}
\usepackage{cite}
\usepackage{textcomp}

\usepackage{amsmath,amssymb,graphicx,amsfonts}
\usepackage{latexsym,todonotes}
\usepackage{mathtools}
\usepackage{cite}
\usepackage{color}
\usepackage{arydshln}
\usepackage{subfigure}
\usepackage{graphbox}
\usepackage{bm}

\usepackage{enumitem}
\usepackage{cases}
\usepackage{dsfont}

\usepackage[vlined,boxed,linesnumbered]{algorithm2e}
\SetAlCapSkip{0.5em}
\usepackage{url}
\usepackage{hyperref}

\usepackage{mathrsfs,BOONDOX-cal}
\graphicspath{{./figs/}}

\newtheorem{theorem}{\bf Theorem}
\newtheorem{lemma}{\bf Lemma}
\newtheorem{proposition}{\bf Proposition}
\newtheorem{corollary}{\bf Corollary}

\newcounter{example}
\newenvironment{example}[1][]{\refstepcounter{example}\par\medskip
   \noindent 
{\bf Example~\theexample:}~#1\it}{\medskip}

\newcommand{\highlight}[1]{{#1}}

\newcommand{\plainproof}{\emph{Proof.}~}
\newcommand{\proofRemark}[1]{\emph{Proof of #1.}~}
\def\endproof{\hfill$\blacksquare$}

\newcounter{definition}
\newenvironment{definition}[1]{\refstepcounter{definition}
\par\medskip
\noindent \textbf{#1:~} \it}{\medskip \newline}

\newcommand*{\fullref}[2]{\hyperref[{#1}]{\ref*{#1} {\it(#2)}}}
\newcommand*{\partialref}[2]{\hyperref[{#1}]{#2}}

\makeatletter
\def\env@cases{%
  \let\@ifnextchar\new@ifnextchar
  \left.
  \def\arraystretch{1.2}%
  \array{@{}l@{\,\,}l@{}}%
}%
\makeatother

\begin{document}
\title{Coordinated Multi-Agent Patrolling with State-Dependent Cost Rates: Asymptotically Optimal Policies for Large-Scale Systems}
\author{Jing Fu, \IEEEmembership{Member, IEEE}, Zengfu Wang, and Jie Chen, \IEEEmembership{Fellow, IEEE}
\thanks{Jing Fu is with School of Engineering, STEM College, RMIT University, VIC3000, Australia (e-mail: jing.fu@rmit.edu.au). }
\thanks{Zengfu Wang is with the School of Automation, Northwestern Polytechnical University, Xi'an 710072, China, and with the Research \& Development Institute of Northwestern Polytechnical University in Shenzhen, Shenzhen 518057, China.(e-mail: wangzengfu@nwpu.edu.cn).}
\thanks{Jie Chen is with the Department of Electrical Engineering, City University of Hong Kong, Hong Kong, China. (e-mail: jichen@cityu.edu.hk). }
}

\maketitle

\begin{abstract}
We study a large-scale patrol problem with state-dependent costs and multi-agent coordination.
We consider heterogeneous agents, rather general reward functions, and the capabilities of tracking agents' trajectories.
Given the complexity and uncertainty of the practical situations for patrolling, we model the problem as a discrete-time Markov decision process (MDP) that consists of a large number of parallel stochastic processes.
We aim to minimize the cumulative patrolling cost over a finite time horizon.
The problem exhibits an excessively large size of state space, which increases exponentially in the number of agents and the size of geographical region for patrolling.
To reach practical solutions, we relax the dependencies between these parallel stochastic processes by randomizing all the state and action variables.
In this context, the entire problem can be decomposed into a number of sub-problems, each of which has a much smaller state space and can be solved independently. The solutions of these sub-problems can lead to efficient heuristics.
Unlike the past systems assuming relatively simple structure of the underlying stochastic process, here, tracking the patrol trajectories involves strong dependencies between the stochastic processes, leading to entirely different state and action spaces, transition kernels, and behaviours of processes,
rendering the existing methods inapplicable or impractical.
Further more, we prove that the performance deviation between the proposed policies and the possible optimal solution diminishes exponentially in the problem size, which also establishes the fact that the policies converge asymptotically at an exponential rate.
\end{abstract}


\section{Introduction}
\label{sec:intro}

\IEEEPARstart{M}{odern} technologies have enabled mission-based agents, such as robots, mobile sensors or unmanned-aerial-vehicles (UAVs), to explore unknown areas with appropriate patrol strategies. 
The patrolling problems have been studied in a wide range of practical scenarios driven by real-world applications, such as mobile sensor networks in smart cities~\cite{jamil2015smart}, green security~\cite{xu2021dual}, multi-UAV monitoring~\cite{mersheeva2015multi}, and emergency response in disasters~\cite{zhou2019online}. 
Relevant problems have been explored 
with situational awareness and objectives. 
For instance,  minimization of patrolling costs with persistent agents has been studied in \cite{rezazadeh2021submodular,wang2021modeling}, where the cost rate was defined as a function of the visited location and the idle time since the last visit to that location. 
Some other works aimed to maximize the opportunities to capture malicious intruders with moving agents~\cite{singh2003optimal,duan2021stochastic}.

In this paper, we study patrolling problems with coordinated heterogeneous agents and aim to minimize the cumulative cost, where the cost rates of the patrolled locations are dependent on the location's history of being patrolled. 
The patrolling history of a location can be specified as the idle time since the last visit, probabilities of capturing intruders, or some more general and practical history data related to the patrolling process, such as the error covariance  matrices of the (non-linear) Kalman filter when gathering information~\cite{jawaid2015submodularity,zhang2017sensor,chamon2021approximate}. 
In out present settings, this translates into state-dependent cost rates.

Patrolling problems with a single or identical agents have been widely studied for decades~\cite{singh2003optimal,wang2021modeling,duan2021stochastic,mallya2022priority}.
To mitigate the uncertainty of sensing responses, multi-agent patrolling was considered with assumed submodular reward functions~\cite{stranders2013near,lee2021mobile}, where the proposed algorithms were proved to achieve certain approximation ratios. Nonetheless, non-submodularity and/or large numbers of employed agents prevent the techniques from being applied to many practical cases with realistic applications, such as Kalman-filter-based uncertainty mitigation~\cite{jawaid2015submodularity,zhang2017sensor}.
In \cite{rezazadeh2021submodular}, a persistent patrolling process was measured by the rewards gained while visiting discretized locations of the patrolled area. The reward function was assumed to be concave and dependent only on the location and the idle time since the last visit to that location. 
In \cite{chamon2021approximate, xu2021dual}, agents (or sensors) were assumed to be able to detect any place at any time without tracking the moving trajectories of the agents and their largest detection coverage.
Here, we consider coordinated heterogeneous agents that patrol a large geographical area. 
We aim to minimize the cumulative costs with practical cost functions that are not necessarily concave or submodular and take consideration of the patrol history of the explored locations. 
This renders our method more flexible and  applicable.

In \cite{cassandras2012optimal, zhou2018optimal,pinto2022multiagent}, multi-agent persistent monitoring problems on a 1-D line were considered with fixed multiple targets and target states dependent on the targets' monitoring history, aiming at minimizing the uncertainties of the systems under consideration.
The authors proved that, in such a 1-D case, there exist piece-wise linear trajectories optimal for moving the agents.
In \cite{lin2014optimal,pinto2022multiagent}, similar techniques were discussed for 2-D planes, where the piece-wise linear trajectories are in general sub-optimal, with numerically demonstrated effectiveness of heuristic algorithms.
This paper considers similar scenarios except that we allow randomness of both target states and agent trajectories, even if a deterministic algorithm is employed.
This enables modeling of more complex detection information of moving targets with unknown or uncertain moving trajectories, other than the binary-style information - detected/monitored and not-detected/non-monitored - of fixed targets.

The coordination of heterogeneous agents, the practical cost rates and the size of the substrate patrolling area substantially complicate the problem and prevent conventional techniques from being applied.
To mitigate these difficulties, we employ a novel approach based on restless multi-armed bandit (RMAB) problem \cite{whittle1988restless},
in which the RMAB techniques can be utilized to achieve approximated optimality and to mitigate the curse of dimensionality.
In particular, RMAB problem is a Markov decision process (MDP) that consists of parallel \emph{bandit processes}, each of which is an MDP with binary actions. 
The cost rates that are dependent on the locations' patrolling history will be reflected in the state variables of the bandit processes, which we refer to as the state-dependent cost rates.
The RMAB technique has been widely used in cases where a large number of bandit processes compete for limited opportunities of being selected. The RMAB problem suffers the curse of dimensionality and has been proved to be PSPACE-hard in general~\cite{papadimitriou1999complexity}.
In~\cite{whittle1988restless}, for the continuous-time RMAB problem, it is conjectured that a simple policy, which was subsequently referred to as the \emph{Whittle index policy}, approached optimality when the size of the problem tends to infinity. This property is referred to as \emph{asymptotic optimality}. Later in \cite{weber1990index}, the conjecture was proved under a non-trivial condition related to the existence of a global attractor of the underlying stochastic process. 
The non-trivial condition, as well as the asymptotic optimality, has  been proved in several different cases~\cite{fu2016asymptotic,fu2020resource,fu2020energy}.
In \cite{brown2020index,balseiro2021dynamic}, discrete-time, finite-time horizon RMAB processes were analyzed through the \emph{Lagrangian dynamic programming (DP)} technique.
Although Lagrangian DP is usually not applicable to continuous-time or infinite time horizon cases, due to its computational complexity that increases polynomially in the number of the discrete time slots, 
it extends the technique in \cite{whittle1988restless} and  leads to proved asymptotic optimality for the discrete-time RMAB problems with finite time horizons.

Nevertheless, the multi-agent patrolling problem is far from a conventional RMAB process studied in the past \cite{whittle1988restless,brown2020index,fu2020resource}.
A key difference lies in the agent-trajectory-dependent constraints over the sets of eligible actions for the next decision epoch. 
Conventional RMAB processes describe problems of selecting $M$ over $N$ ($N > M$) parallel stochastic processes, each of which can or cannot be selected according to only its own state and is related to at most one constraint over its action variables.
While, in the multi-agent patrolling problem, the eligibility of moving an agent to another position is usually dependent on its own position (state), the positions of the other agents (states of other processes), and the overall topology of the patrolling region.
It follows with multiple constraints that are related to each parallel stochastic process and
a stronger dependency among their control variables.
The existing RMAB methods cannot be applied to the patrol problem nor will it ensure the existence of practical patrol policies with bounded performance degradation.

More specifically, we consider a problem consisting of multiple stochastic processes, each of which is a discrete-time MDP with finite state and action spaces. The state and action variables of these MDPs are coupled through multiple constraints that are linear in the action variables. 
We aim to minimize the expected cumulative total costs of all the MDPs in a finite time horizon and refer to such a problem as the Multi-Action Bandits with Multiple Linear Constraints (MAB-ML).
The term ``multi-action bandit" is used in the sense that the above-mentioned MDP acts as a bandit process with multiple actions.
A rigorous definition of MAB-ML is provided in Section~\ref{subsec:MAB-ML}.
Both the conventional (finite-time-horizon) RMAB and the patrol problem are special cases of such MAB-ML. 
MAB-ML is at least as difficult as RMAB, which is considered to be hard in general \cite{whittle1988restless,weber1990index,brown2020index,papadimitriou1999complexity}.


Our main contributions are summarized as follows.  
We formulate a discrete-time MAB-ML problem and
model the multi-agent patrolling problem as a special case of the MAB-ML problem.
In this vein, we incorporate the coordination of multiple agents and state-dependent cost rates and consider more general patrolling models, allowing heterogeneous agents, more general reward functions, and the capabilities of tracking agents' trajectories.
As alluded to above, the patrolling problem is complicated by the high-dimensionality of the state space, of which the size increases exponentially in the number of agents and the domain of the patrol region.
Since each parallel stochastic process is constrained by more than one action constraints, the available RMAB results are not applicable, nor can the existing technique ensure the existence of near-optimal patrol policies.

For the general MAB-ML, we adapt the \emph{Whittle relaxation technique} \cite{whittle1988restless} and the Lagrangian DP \cite{brown2020index} to decompose the patrolling problem into independent sub-problems with reduced complexity. 
We propose a patrolling policy that does not consume excessive computational or storage power when the problem size, measured by the number of agents and locations included in the patrol region, becomes large. 
This ensures the policy to be \emph{scalable}, a key attribute to large-scale patrolling problems. 
The proposed policy quantifies the marginal costs of moving an agent to regions by \emph{indices} that prioritize the movements with the least index values. 
We show that the complexity of implementing the proposed policy is at most linear-logarithmic in the number of agents and quadratic in the size of the patrolling region.

More importantly, by employing Freidlin's theorem \cite[Chapter 7]{freidlin2012random}, 
we prove that, under a mild condition related to the tie-case of the indices, such a policy approaches optimality when the {MAB-ML} problem size tends to infinity; that is, the policy is asymptotically optimal.
Furthermore, we show that the performance deviation between the proposed policy and the optimal policy diminishes exponentially in the size of the problem. 
This deviation bound is tighter than the $O(\sqrt{N})$ bound for the standard discrete-time finite-time-horizon RMAB problem achieved in \cite{brown2020index}.
These results indicate that the proposed policy quickly gets better when the problem becomes larger, and it is appropriate for large-scale problems.

The remainder of the paper is organized as follows.
In Section~\ref{sec:model}, a description of the patrol
model and the MAB-ML problem are introduced. 
In Section~\ref{sec:relaxation}, we relax MAB-ML through the Whittle relaxation and Lagrangian DP techniques, achieving intermediate results for further analysis.
In Section~\ref{sec:policies}, we propose index-based policies that are applicable to large-scale systems.
We provide in Section~\ref{sec:asym_opt} our main results for MAB-ML, as well as the patrol case.
In Section~\ref{sec:simulation}, numerical results are presented to demonstrate the theoretical results.
The paper concludes in Section~\ref{sec:conclusion}.

\section{Model}\label{sec:model}
For any positive integer $N$, let $[N]$ represent the set $\{1,2,\ldots,N\}$. Let $\mathbb{R}$, $\mathbb{R}_+$ and $\mathbb{R}_0$ be the sets of all real numbers, positive real numbers, and non-negative real numbers, respectively.
Similarly, $\mathbb{N}$, $\mathbb{N}_+$ and $\mathbb{N}_0$ are the sets of integers, positive integers and non-negative integers, respectively.
Let $\bm{1}^I_i$ represent a standard basis vector of $\mathbb{R}^I$ with all elements zero except the $i$th element equal to $1$.
\highlight{We use $\lVert \cdot \rVert$ as the Euclidean norm.}
Important notation of this paper is summarized in Appendix~\ref{app:notation}.

\subsection{Multi-Agent Patrol Problem}\label{subsec:patrolling}\color{black}
Consider a \emph{geographical region} consisting of $I\in\mathbb{N}_+$ different \emph{areas}
and $J\in\mathbb{N}_+$ different types of agents. 
Each type includes $M_j\in\mathbb{N}_+$, $j\in[J]$, agents that explore all the $I$ areas. Assume that the number of type-$j$ agents $M_j$ is no greater than the number of areas $I$.
Agents of different types can be considered as moving sensors with different functions or profiles, such as cameras and thermal sensors.
In each time slot $t \in [T]$ with a finite time horizon $T\in\mathbb{N}_+$, each agent, based on the employed patrol policy, moves to another area or stays in the original place and then scans (detects or explores) where it is before the next time slot. 
Consider a limited set of areas, $\mathscr{B}_{i,j}\in 2^{[I]}$, to which an agent of type $j$ can move from area $i$.
We refer to $\mathscr{B}_{i,j}$ as the set of \emph{neighbour areas} or, alternatively, the \color{black}\emph{neighbourhood} of area $i$ for agent type $j$. For all $i,i'\in[I]$ and $j\in[J]$, assume that $i\in \mathscr{B}_{i,j}$.
If $i\in\mathscr{B}_{i',j}$ then $i'\in\mathscr{B}_{i,j}$.
The neighbour areas are determined by intrinsic features of the patrol region, which can be instantiated as city areas connected by intricate \color{black}streets, aerial areas being detected by UAVs, wild areas formed by geographical borders, et cetera.
The reachability of a neighbourhood may vary across different agent types due to the maneuverability of different agents.
In each time slot, if an area has parked an agent, then other agents of the same type will not explore the same area.
That is, agent collision of the same type or repeated sensing within the same time unit is excluded.

For $i\in[I]$, $j\in[J]$ and $t\in[T]$, let a random variable $K_{i,j}(t)$ that takes values in a finite set $\mathscr{K}_{i,j}$ represent  the controller's knowledge of area $i$ by the end time slot $t-1$ that was collected or sensed by the agents of type $j$. The value of $K_{i,j}(t)$ can be dependent on the whole patrol and observation history associated with area $i$ and agent type $j$ by time $t$. 
\highlight{The process $\{K_{i,j}(t),t\in[T]\}$ evolves according to two different transition matrices that correspond to two cases: area $i$ is detected by an agent of type $j$ at time $t$ or not. 
The transition matrices for different $(i,j)\in[I]\times[J]$ are potentially different.
Such} $K_{i,j}(t)$ may be specified as the idle length by time $t$ since last time area $i$ was visited by an agent of type $j$, the error covariance matrix of (non-linear) Kalman filter by time $t$, the estimated probability of detecting an attacker through historical observations by time $t$, et cetera.
Define $I_{i,j}(t)$ as an indicator for the existence of a type-$j$ agent in area $i$ at time $t$. If there is an agent of type $j$ in area $i$ at the end of time slot $t-1$, then $I_{i,j}(t)=1$; otherwise, $I_{i,j}(t)=0$.

For $i\in[I]$ and $j\in[J]$, define $S_{i,j}(t) \coloneqq (K_{i,j}(t), I_{i,j}(t))$ as the state variable of a stochastic process associated with \emph{area-agent (AA) pair} $(i,j)$, and define $\bm{S}(t) \coloneqq (S_{i,j}(t): i\in[I],j\in[J])$.
There is a central controller of the patrol problem that makes decision at each time $t$ based on the value of $\bm{S}(t)$. 
More precisely, consider an action variable $a_{i,i',j}(\bm{S}(t),t)$, a function of the system state $\bm{S}(t)$ and $t$ and taking a value in $\{0,1\}$, for AA pair $(i,j)$ and an area $i'$ ($i,i'\in[I]$ and $j\in[J]$).
If $a_{i,i',j}(\bm{S}(t),t) =1$, then an agent of type $j$ will move to area $i$ from area $i'$ at the beginning of time $t$ and sense the area before time slot $t+1$; otherwise, no such agent will go from area $i'$ to area $i$ until the next time slot. Note that, if the agent is decided to stay in area $i$ at both time slots $t$ and $t+1$, then $a_{i,i,j}(\bm{S}(t),t) =1$. 
Since an agent can move only one step for each time slot, if $i'\in [I]\backslash\mathscr{B}_{i,j}$, then $a_{i,i',j}(\bm{S}(t),t) \equiv 0$. 
For an AA pair $(i,j)\in[I]\times[J]$ and time slot $t\in [T]$, define the action vector $\pmb{a}_{i,j}(\bm{S}(t),t) \coloneqq (a_{i,i',j}(\bm{S}(t),t): i'\in[I])$, and define $\pmb{a}(\bm{S}(t),t) \coloneqq (a_{i,i',j}(\bm{S}(t),t): i,i'\in[I],j\in[J])$ as the action vector for the entire system.
To track the trajectories of the agents, these action variables should satisfy 
\begin{equation}\label{eqn:constraint:exclusive}
\sum\nolimits_{i'\in\mathscr{B}_{i,j}}a_{i,i',j}(\bm{S}(t),t)\leq 1,~\forall i\in[I],j\in[J],t\in[T],
\end{equation}
and
\begin{equation}\label{eqn:constraint:neighbourhood}
\sum\nolimits_{i'\in\mathscr{B}_{i,j}}a_{i',i,j}(\bm{S}(t),t) = I_{i,j}(t),
\forall i\in[I],j\in[J],t\in[T].
\end{equation}
Inequality~\eqref{eqn:constraint:exclusive} indicates that there is at most one agent of the same type that can move to the same area during the same time slot. 
Based on the definition of the action vector and Constraint~\eqref{eqn:constraint:exclusive}, $\pmb{a}_{i,j}(\bm{S}(t),t)$ takes $|\mathscr{B}_{i,j}|+1$ possible values: $\bm{1}^I_{i'}$ for all $i'\in\mathscr{B}_{i,j}$ and $\bm{0}$ otherwise.
Equation~\eqref{eqn:constraint:neighbourhood} guarantees that, for each AA pair $(i,j)$, if there is an agent of type $j$ located in area $i$ at the end of time slot $t-1$, this agent cannot disappear. 
It must be moved to an area in the neighbourhood of $i$, including staying at $i$, at the beginning of the next time slot.

\begin{figure*}[t]
\centering
\subfigure[Initialization at the beginning of time slot $t = 1$.]{
\includegraphics[width = 0.25\textwidth]{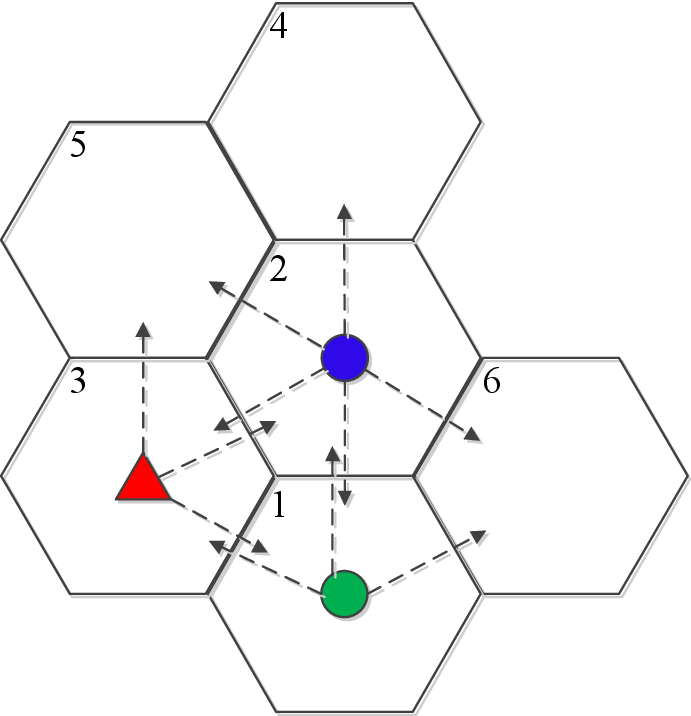}
\label{fig:examplesT0}}
\subfigure[During time slot $t = 1$.]{
\includegraphics[width = 0.25\textwidth]{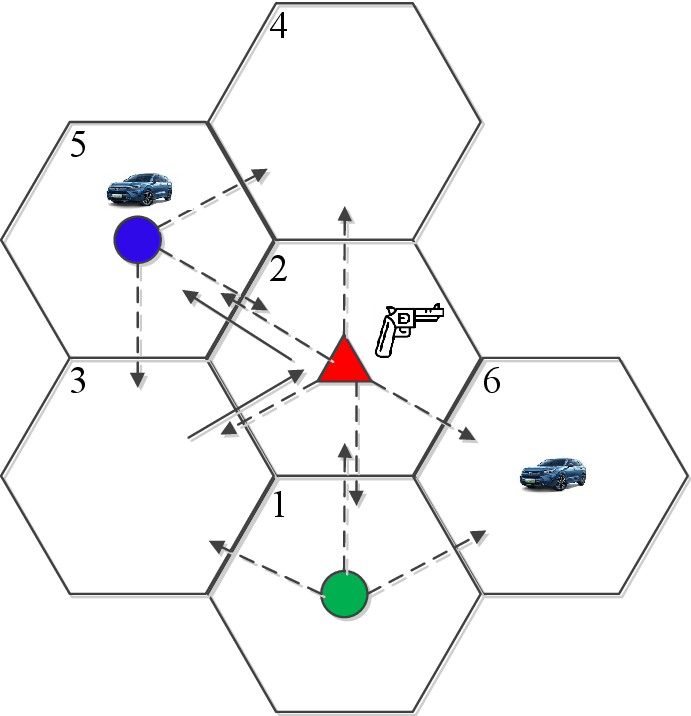}
\label{fig:examplesT1}}
\subfigure[During time slot $t = 2$.]{
\includegraphics[width = 0.25\textwidth]{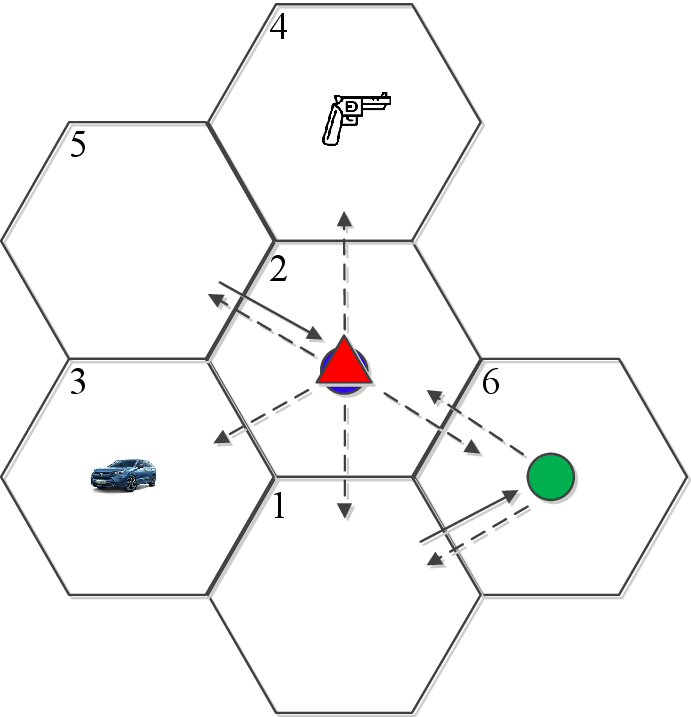}
\label{fig:examplesT2}}
\caption{An example with $I = 6$ areas and $J = 2$ types of agents: the round circles are agents of type $1$, and the triangle is an agent for type $2$.\label{fig:model_example}}
\end{figure*}

We illustrate the aforementioned concepts and terminologies using the following example.
\begin{example}
In Figure \ref{fig:model_example}, we consider a simple example for crime detection. 
There are $J=2$ types of agents: two for crime type 1 (possession of weapons) and one for type 2 (vehicle crime).
In the figure, Agents 1 and 2 are in type 1 denoted by a blue and a green circle, respectively, and Agent 3 is in type 2 and a red triangle.
The three agents are patrolling $I=6$ areas, aiming to minimize the expected total crime rates.
Here, each hexagon represents an area with crime rate $Q_{ij}(t) \in [0,1]$ ($(i,j)\in[I]]\times [J]$), which follows a beta distribution with parameters $(\alpha_{i,j}(t), \beta_{i,j}(t))$ at time $t\in[T]$.
At each time $t$, a random signal $D_{i,j}(t) \in \{0,1\}$ is drawn from a Bernoulli distribution with success probability $Q_{i,j}(t)$. If $D_{i,j}(t) = 1$, then a crime event of type $j$ is reported in area $i$; otherwise, no crime reported in area $i$. 
A reported crime event is either reported by a patrol agent (when the agent patrols in the corresponding area) or, when there is no patrol agent, reported by local residents or victims.
If a type-$j$ agent patrols in area $i$ at time $t$ and $D_{i,j}(t) = 1$ (the agent detects a crime and arrests the criminals), then $\alpha_{i,j}(t+1) = \alpha_{i,j}(t)$ and $\beta_{i,j}(t+1) = \beta_{i,j}(t)+\Delta_{\beta}$;
if a type-$j$ agent patrols in area $i$ at time $t$ and $D_{i,j}(t) =0$ (no criminal was arrested by the agent),
then $\alpha_{i,j}(t+1) = \alpha_{i,j}(t)+\Delta^1_{\alpha}$ and $\beta_{i,j}(t+1) = \beta_{i,j}(t)$.
Similarly,
if there is no type-$j$ agent patrols in area $i$ at time $t$ and $D_{i,j}(t)=1$, then $\alpha_{i,j}(t+1) = \alpha_{i,j}(t) + \Delta^2_{\alpha}$ and $\beta_{i,j}(t+1) = \beta_{i,j}(t)$; otherwise, $(\alpha_{i,j}(t+1),\beta_{i,j}(t+1)) = (\alpha_{i,j}(t),\beta_{i,j}(t))$.
Here, the parameters $\Delta_{\beta},\Delta^1_{\alpha}, \Delta^2_{\alpha}$ are positive constants.
In this context, $(\alpha_{i,j}(t), \beta_{i,j}(t))$ maintains sufficient statistics for measuring the crime rates at time $t$.
Let $K_{i,j}(t)=(\alpha_{i,j}(t), \beta_{i,j}(t))$, representing the controller's knowledge of area $i$ by time $t$.
In Figure \ref{fig:model_example}, the dashed arrows point to the areas the agents can possibly move to, and the solid arrows indicate the actual movements taken by the controller.

In the beginning of time slot $t = 1$ (in Figure~\ref{fig:examplesT0}), for AA pairs $(i,j)=(1,1)$, $(2,1)$ and $(3,2)$, we initialize their states $S_{i,j}(1)=((\alpha_{i,j}, \beta_{i,j}),1)$ with given $\alpha_{i,j},\beta_{i,j} $; and, for all the other AA pairs $(i,j)\in[I]\times[J]$, $S_{i,j}(1)=((\alpha_{i,j}, \beta_{i,j}),0)$. 
It means that, at time $t=1$, the three agents are initially located in areas  $1$, $2$, and $3$, and the controller's knowledge for all the six areas is $K_{i,j}(1)=(\alpha_{i,j}, \beta_{i,j})$.
Based on $\bm{S}(1)$, under some patrol policy, the controller decides that $a_{1,1,1}(\bm{S}(1),1) = 1$,
$a_{5,2,1}(\bm{S}(1),1) = 1$,
$a_{2,3,1}(\bm{S}(1),1) = 1$,
and $a_{i,i',j}(\bm{S}(1),1) = 0$ for the remaining $i,i'\in [I]$.
That is, during time slot $t=1$, the controller will keep Agent 2 in area $1$, move Agent 1 from area $2$ to $5$, 
and move Agent 3 from area $3$ to $2$,
as described in Figure~\ref{fig:examplesT1}.
During time slot $t=1$, the three agents will be moved to areas $1$, $5$, and $2$, and then they will patrol/explore these three areas before the next time slot. 
While patrolling, some crimes get reported with $D_{2,2}(1) = D_{5,1}(1) = D_{6,1}(1) = 1$, and $D_{i,j}(1) = 0$ for the other $(i,j) \in [I] \times [J]$. Based on the reported crimes (either reported by agents or local residents/victims), we obtain the values of $K_{i,j}(2)$ based on the transition rules described in the previous paragraph.

In the beginning of time slot $t=2$,
the states of AA pairs $(1,1)$, $(2,2)$, $(5,1)$ and $(6,1)$ become $S_{1,1}(2)=((\alpha_{1,1}(1)+\Delta^1_{\alpha}, \beta_{1,1}),1)$,
$S_{2,2}(2)=((\alpha_{2,2}(1), \beta_{2,2}+\Delta_{\beta}), 1)$,
$S_{5,1}(2)=((\alpha_{5,1}(1), \beta_{5,1}+\Delta_{\beta}), 1)$, and
$S_{6,1}(2)=((\alpha_{6,1}(1) + \Delta^2_{\alpha}, \beta_{6,1}), 1)$, respectively; 
and, for all the other AA pairs $(i,j)\in[i]\times[J]$, the states variables are $S_{i,j}(2)=((\alpha_{ij}(1), \beta_{ij}(1)),0)$. 
Based on the state vector $\bm{S}(2)$, the controller again decides the action variables $\pmb{a}(\bm{S}(2),2)$, and the agents will be moved accordingly.

Consider $a_{2,5,1}(\bm{S}(2),2) = 1$,
$a_{6,1,1}(\bm{S}(2),2) = 1$,
and $a_{i,i',1}(\bm{S}(2),2) = 0$ for all the other $i,i'\in [I]$ - moving Agents 1 and 2 from areas $5$ and $1$ to areas $2$ and $6$, respectively, during the time slot $t=2$.
For the type-$2$ agent, take action $a_{2,2,2}(\bm{S}(2),2) = 1$ with all the other action variables $a_{i,i',2}(\bm{S},2)=0$ - the agent remains in area 2.
Figure~\ref{fig:examplesT2} demonstrates the movements.
Similarly, after the movements in $t=2$, the agents will not move to other areas until $t=3$. 
While these agents patrol in time slot $t = 2$, the crime occurrence is represented by $D_{3,1}(2) = D_{4,2}(2) = 1$, 
and $D_{ij}(2) = 0$ for the other $(i,j) \in [I] \times [J]$.
Then, based on $\bm{S}(2)$, $\pmb{a}(\cdot,2)$, and the crime occurrence variables $D_{i,j}(2)$, the state vector $\bm{S}(3)$ can be obtained through the above-mentioned transition rules.
This process continues until reaching the horizon $t=T$.
\end{example}

For $i\in[I]$ and $j\in[J]$, define $\mathscr{S}_{i,j}$ as the set of all the possible values of $S_{i,j}(t)$ ($t\in[T]$), and define the state space of the stochastic process $\{\bm{S}(t), t\in[T]\}$ as $\mathscr{S}\coloneqq \prod_{i\in[I]}\prod_{j\in[J]}\mathscr{S}_{i,j}$, where $\prod$ is the Cartesian product. Recall that each value $s\in\mathscr{S}_{i,j}$ includes two parts, and we rewrite it as $s=(\mathcal{k}(s),\mathcal{i}(s))$. That is, if $S_{i,j}(t) = s$, $\mathcal{k}(s) = K_{i,j}(t)$ and $\mathcal{i}(s)=I_{i,j}(t)$.
The action vector $\pmb{a}(\cdot,\cdot)$ is a mapping from $\mathscr{S}\times [T]$ to $\{0,1\}^{I^2J}$, which determines a \emph{policy} of the patrol problem. To clarify the dependency between an employed policy, represented by $\phi$, the action variables and the state variables, for any $\pmb{s}\in\mathscr{S}$ and $t\in[T]$, we rewrite $\pmb{a}(\pmb{s},t)$ and $\bm{S}(t)$ as $\pmb{a}^{\phi}(\pmb{s},t)$ and $\bm{S}^{\phi}(t)$, respectively, with the added superscript $\phi$. Similarly, we rewrite $S_{i,j}(t)$, $K_{i,j}(t)$, $I_{i,j}(t)$, $\pmb{a}_{i,j}(\pmb{s},t)$ as $S^{\phi}_{i,j}(t)$, $K^{\phi}_{i,j}(t)$, $I^{\phi}_{i,j}(t)$ and $\pmb{a}^{\phi}_{i,j}(\pmb{s},t)$, respectively. 
Let $\Phi$ represent the set of all the policies $\phi$ determined by the associated action vector $\pmb{a}^{\phi}$.
Process $\{S_{i,j}(t),t\in[T]\}$ associated with each AA pair $(i,j)\in[I]\times [J]$ evolves with different transition probabilities when area $i$ is or is not patrolled by an agent of type $j$.

We aim to minimize the cumulative costs of the patrol process $\{\bm{S}^{\phi}(t), t\in[T]\}$ with the objective
\begin{equation}\label{eqn:obj}
\min_{\phi\in\Phi} \sum_{t\in[T]}\sum_{i\in[I]}\sum_{j\in[J]} \mathbb{E}\biggl[c_{i,j}\Bigl(S^{\phi}_{i,j}(t), \mathcal{e}^{\phi}_{i,j}\bigl(\bm{S}^{\phi}(t),t\bigr),t\Bigr)\biggr],
\end{equation}
subject to \eqref{eqn:constraint:exclusive} and \eqref{eqn:constraint:neighbourhood}, where $\mathcal{e}^{\phi}_{i,j}(\bm{S}^{\phi}(t),t)\coloneqq \mathds{1}\bigl\{\sum_{i'\in\mathscr{B}_{i,j}}a^{\phi}_{i,i',j}(\bm{S}^{\phi}(t),t) > 0\bigr\}$ indicating whether area $i$ is patrolled by an agent of type $j$ at time $t$, $c_{i,j}$ is a non-negative real-valued function representing an instantaneous cost rate, and $\mathbb{E}$ takes expectation over the random variables with a given policy $\phi$ and an initial distribution of $\bm{S}^{\phi}(1)$. 

\subsection{Multi-Action Bandits with Multiple Linear Action Constraints (MAB-ML)}\label{subsec:MAB-ML}

Neither Constraint~\eqref{eqn:constraint:exclusive} nor \eqref{eqn:constraint:neighbourhood} match those in conventional RMAB or multi-action bandit models but are essential for our patrol problem.
The conventional RMAB or multi-action bandits presumes a simple form of action constraints - the action vector $\pmb{a}^{\phi}_{i,j}(\bm{S}^{\phi}(t),t)$  of each bandit process $(i,j)\in[I]\times[J]$ is involved in at most one constraint.
That is, for the conventional RMAB case, each bandit process is coupled with others through at most one constraint (for each time slot). 
In other words, only weak dependencies exist between the bandit processes.
For that matter, a simplified analysis can be conducted and accordingly, algorithms, such as the classic Whittle index policy, were proposed.

For the multi-agent patrol problem, the action vector of each bandit process is subject to $|\mathscr{B}_{i,j}|\geq 2$ different constraints in \eqref{eqn:constraint:neighbourhood} led by different neighbourhoods and (time-varying) states of different AA pairs.
This imposes concrete dependencies intertwined across all the bandit processes and entirely different transition kernel(s) of the underlying MDP, which cannot be addressed by the conventional analysis.
In particular, the classic priority-style policies, including Whittle index policy, for RMAB problems are no longer applicable, as Constraints \eqref{eqn:constraint:exclusive} and \eqref{eqn:constraint:neighbourhood} fail to be applicable.

In what follows, we canonicalize the patrol problem. 
For each $(i,j)\in[I]\times[J]$, the discrete-time process $\{S^{\phi}_{i,j}(t),t\in[T]\}$ is an MDP with finite state space $\mathscr{S}_{i,j}$.
As indicated by Constraint~\eqref{eqn:constraint:exclusive}, its action variable $\pmb{a}^{\phi}_{i,j}(\bm{S}^{\phi}(t),t)$ takes $|\mathscr{B}_{i,j}|+1$ possible actions: $\bm{1}^I_{i'}$ for all $i'\in \mathscr{B}_{i,j}$ and \highlight{a zero vector}.
Process $\{S^{\phi}_{i,j}(t),t\in[T]\}$ evolves with different transition probabilities for $\pmb{a}^{\phi}_{i,j}(\bm{S}^{\phi}(t),t) = \bm{0}$ and $\pmb{a}^{\phi}_{i,j}(\bm{S}^{\phi}(t),t) \neq \bm{0}$.
Such a process $\{S^{\phi}_{i,j}(t),t\in[T]\}$ is a \emph{multi-action bandit process}, which reduces to a standard restless bandit process when $|\mathscr{B}_{i,j}| = 1$.

For $(i,j)\in[I]\times[J]$, we consider a subset $\bar{\mathscr{B}}_{i,j}\subset[I]$ and generalize \eqref{eqn:constraint:neighbourhood} to
\begin{equation}\label{eqn:constraint:canonical}
    \sum\nolimits_{i'\in\bar{\mathscr{B}}_{i,j}}\mathcal{w}_{i',i,j}a^{\phi}_{i',i,j}(\bm{S}^{\phi}(t),t) = g_{i,j}(S^{\phi}_{i,j}(t)),\\
\forall i\in[I],j\in[J],t\in[T],
\end{equation}
where $\mathcal{w}_{i,j}\in\mathbb{R}_0$ and $g_{i,j}(s)$ is a function $\mathscr{S}_{i,j}\rightarrow \mathbb{R}_0$.
For the above described $IJ$ multi-action bandit processes, we refer to the problem defined by \eqref{eqn:obj}, \eqref{eqn:constraint:exclusive} and \eqref{eqn:constraint:canonical} as the \emph{multi-action bandits with multiple linear constraints (MAB-ML)}. 
The patrol problem is a special case of MAB-ML by specifying  $|\mathscr{B}_{i,j}|+1$ actions (moving agents) for each multi-action bandit process, and setting $\bar{\mathscr{B}}_{i,j} = \mathscr{B}_{i,j}$, $\mathcal{w}_{i',i,j}=1$, and $g_{i,j}(s) = \mathcal{i}(s)$, where, for any $i,i'\in[I]$ and $j\in[J]$, if $i\in\mathscr{B}_{i',j}$ then $i'\in\mathscr{B}_{i,j}$.
MAB-ML is in general at least as difficult as the patrol problem. 
As explained in the first two paragraphs of this section, neither MAB-ML nor the patrol problem can be analysed through conventional RMAB techniques.
Note that MAB-ML reduces to RMAB by specifying $J=1$, $\mathscr{B}_{i,j}=\{1\}$ for all $(i,j)\in[I]\times[J]$, and $\bar{\mathscr{B}}_{i,j} = \emptyset$ and $g_{i,j}(\cdot) = 0$ for all $(i,j)\in[I]\times[J]$ except $\bar{\mathscr{B}}_{1,j} = [I]$, $g_{1,j}(\cdot) = M$, and $\mathcal{w}_{i',1,j}=1$ for all $i'\in[I]$. 
Here, $M\in[I]$ is a positive integer for the standard RMAB. 
Hence, solving MAB-ML is at least as difficult as RMAB, which is in general hard with intractable solutions~\cite{whittle1988restless,weber1990index,brown2020index,ninomora2020verification,fu2024restless} and is proved to be PSPACE-hard for the infinite time horizon case~\cite{papadimitriou1999complexity}.

Recall that the patrol problem, a special case of MAB-ML, already lies outside the scope of RMAB or multi-action RMAB.
Nor can the past RMAB techniques or achievements be applicable to MAB-ML.
In this paper, Theorems~\ref{theorem:asym_opt}, \ref{theorem:IND_asym_opt}, and \ref{theorem:exp_convergence}, Lemma~\ref{lemma:indexability}, and Propositions~\ref{prop:strong_duality} and \ref{prop:decomposition} are applicable in general to MAB-ML; and Lemma~\ref{lemma:MAI_asym_opt} and Corollary~\ref{coro:MAI_asym_opt} are specified results for the patrol case.


\section{Lagrangian Dynamic Programming and Relaxation}\label{sec:relaxation}

MAB-ML, including the patrol problem as a special case, exhibits a large state space that exponentially increases in the number of agents and the exploring areas. Conventional MDP techniques such as the value iteration cannot achieve optimality without consuming excessively large amount of computational and storage power. We resort to simple heuristic policies that are applicable to large-scale systems and achieve an established near-optimality.

To proceed,  we randomize the action variables and relax \eqref{eqn:constraint:canonical} to
\begin{equation}\label{eqn:constraint:neighborhood:relax}
\sum\nolimits_{i'\in \bar{\mathscr{B}}_{i,j}\mathcal{w}_{i',i,j}}\mathbb{E}\Big[a^{\phi}_{i',i,j}(\bm{S}^{\phi}(t),t)\Bigr] = \mathbb{E}\Big[g_{i,j}(S^{\phi}_{i,j}(t))\Big],\\~\forall i\in[I], j\in[J], t\in[T],
\end{equation}
where the action variables $a^{\phi}_{i,i',j}(\bm{S}^{\phi}(t),t)$ take values in $\{0,1\}$ with policy-dependent probabilities. 
Let $\bar{a}^{\phi}_{i,i',j}(\pmb{s},t)\in[0,1]$ ($\pmb{s}\in\mathscr{S}$, $i,i'\in[I]$, $j\in[J]$, $t\in[T]$) represent the probability of taking $a^{\phi}_{i,i',j}(\pmb{s},t)=1$ under policy $\phi$, and let $\bar{\pmb{a}}^{\phi}_{i,j}(\pmb{s},t)\coloneqq (\bar{a}^{\phi}_{i,i',j}(\pmb{s},t):i'\in[I])$.
In this context, for $i\in[I]$, $j\in[J]$, $t\in[T]$, $\pmb{s}\in\mathscr{S}$, $\phi\in\tilde{\Phi}$ and any $i'\in[I]\backslash\mathscr{B}_{i,j}$, $\bar{a}^{\phi}_{i,i',j}(\pmb{s},t) \equiv 0$.
Define a set $\bar{\Phi}$ of all the policies $\phi$, each of which is determined by $\bar{a}^{\phi}_{i,i',j}(\pmb{s},t)$ for all $\pmb{s}\in\mathscr{S}$, $i,i'\in[I]$, $j\in[J]$ and $t\in[T]$, and define $\tilde{\Phi}$ as the set of all the policies $\phi\in\bar{\Phi}$ such that, for $(i,j)\in[I]\times[J]$, $\pmb{s}\in\mathscr{S}$, and $t\in[T]$, $\bar{\pmb{a}}^{\phi}_{i,j}(\pmb{s},t)$ takes values in 
\begin{equation}\label{eqn:define_A}
\mathscr{A}_{i,j}\coloneqq\biggl\{\pmb{a}\in[0,1]^I \biggl| \sum_{i'\in\mathscr{B}_{i,j}}a_{i'} \leq 1 \\\text{ and } a_{i'}=0 \text{ for all }i'\in[I]\backslash\mathscr{B}_{i,j}\biggr\}.
\end{equation}
Consider the problem
\begin{equation}\label{eqn:obj:relax}
    \min_{\phi\in\tilde{\Phi}} \sum_{t\in[T]}\sum_{i\in[I]}\sum_{j\in[J]} \mathbb{E}\biggl[c_{i,j}\Bigl(S^{\phi}_{i,j}(t),\mathcal{e}^{\phi}_{i,j}\bigl(\bm{S}^{\phi}(t),t\bigr),t\Bigr)\biggr],
\end{equation}
subject to \eqref{eqn:constraint:neighborhood:relax}.
We refer to the problem described in \eqref{eqn:obj:relax} and \eqref{eqn:constraint:neighborhood:relax} as the \emph{relaxed} version of the original MAB-ML described in \eqref{eqn:obj}, \eqref{eqn:constraint:exclusive} and \eqref{eqn:constraint:canonical}.
In this context, if a policy $\phi\in\Phi$ is applicable to the original problem (satisfying \eqref{eqn:constraint:exclusive} and \eqref{eqn:constraint:canonical}), then $\phi\in\tilde{\Phi}$ and $\phi$ is also applicable to the relaxed problem (satisfying \eqref{eqn:constraint:neighborhood:relax}).
Note that, for the relaxed problem, Constraint~\eqref{eqn:constraint:exclusive} remains unchanged and is included in the action space $\mathscr{A}_{i,j}$ defined in \eqref{eqn:define_A}.
The minimum of the relaxed problem constitutes a lower bound for that of the original problem (described in \eqref{eqn:obj}, \eqref{eqn:constraint:exclusive} and \eqref{eqn:constraint:canonical}).
The dual function of the relaxed problem is
\begin{multline}\label{eqn:dual_func}
L(\pmb{\gamma}) = \min_{\phi\in\tilde{\Phi}} \sum_{t\in[T]} \sum_{j\in[J]}
\sum_{i\in[I]}\mathbb{E}\biggl[c_{i,j}\Bigl(S^{\phi}_{i,j}(t),\mathcal{e}^{\phi}_{i,j}\bigl(\bm{S}^{\phi}(t),t\bigr),t\Bigr)\\
+\gamma_{i,j,t}{g_{i,j}(S^{\phi}_{i,j}(t))}-\sum_{i'\in\bar{\mathscr{B}}^{-1}_{i,j}}\gamma_{i',j,t}\mathcal{w}_{i,i',j}a^{\phi}_{i,i',j}(\bm{S}^{\phi}(t),t) \biggr] ,
\end{multline}
where $\pmb{\gamma}\in\mathbb{R}^{IJT}$ are the Lagrange multipliers~\cite{bertsekas2014constrained} for Constraint~\eqref{eqn:constraint:neighborhood:relax}, and $\bar{\mathscr{B}}^{-1}_{i,j}\coloneqq\{i'\in[I] | i\in\bar{\mathscr{B}}_{i',j}\}$.
For the patrol problem where $\mathscr{B}_{i,j}$ is specified as the neighbourhood of AA pair $(i,j)$ and $\bar{\mathscr{B}}_{i,j} = \mathscr{B}_{i,j}$, we have $\bar{\mathscr{B}}^{-1}_{i,j}= \mathscr{B}_{i,j}=\bar{\mathscr{B}}_{i,j}$.
The right hand side of \eqref{eqn:dual_func} is the minimized cumulative expected costs of process $\{\bm{S}^{\phi}(t), t\in[T]\}$ with the expected cost in state $\pmb{s}\in\mathscr{S}$ for action $\pmb{a}=(a_{i,i',j}:i,i'\in[I],j\in[J])\in \prod_{i\in[I],j\in[J]}\mathscr{A}_{i,j}$ given by
\begin{equation}    \sum_{j\in[J]}\biggl[\sum_{i\in[I]}\Bigl(c_{i,j}\Bigl(s_{i,j},\mathcal{e}\bigl(\pmb{a}_{i,j}\bigr),t\Bigr)+\gamma_{i,j,t}{g_{i,j}(s_{i,j})}\\-\sum\nolimits_{i'\in\bar{\mathscr{B}}^{-1}_{i,j}}\gamma_{i',j,t}\mathcal{w}_{i,i',j}a_{i,i',j} \Bigr)\biggr],
\end{equation}
where we recall that the action variable $a_{i,i',j}\in[0,1]$ represents the probability of  taking action $(i,i',j)$, $\pmb{a}_{i,j}\coloneqq (a_{i,i',j}:i'\in[I])$, $\mathcal{e}(\pmb{a}_{i,j})\coloneqq 1-\prod_{i'\in [I]}(1-a_{i,i',j})$ representing the probability of taking a non-zero action for process $\{S^{\phi}_{i,j}(t),t\in[T]\}$, and, 
for $e\in(0,1)$, $(i,j)\in[I]\times[J]$, $s\in\mathscr{S}_{i,j}$ and $t\in[T]$, $c_{i,j}(s,e,t)\coloneqq ec_{i,j}(s,1,t)+(1-e)c_{i,j}(s,0,t)$.
For the patrol problem, as a special case of MAB-ML, $\mathcal{e}(\pmb{a}_{i,j})$ represents the probability of moving a type-$j$ agent to area $i$.


\subsection{Decomposition of the Relaxed MAB-ML}\label{subsec:decomposition}

For $i\in[I]$, $j\in[J]$ and $\pmb{\gamma}_{i,j}\coloneqq (\gamma_{i',j,t}\in\mathbb{R}: i'\in \{i\}\cup\bar{\mathscr{B}}^{-1}_{i,j},t\in[T])$, define
\begin{equation}\label{eqn:define_func_i}
L^{\phi}_{i,j}(\pmb{\gamma}_{i,j}) \coloneqq \sum_{t\in[T]}
\sum_{\pmb{s}\in\mathscr{S}}
C^{\pmb{\gamma}_{i,j}}_{i,j}(s_{i,j},\bar{\pmb{a}}^{\phi}_{i,j}(\pmb{s},t),t)  \mathbb{P}\Bigl\{\bm{S}^{\phi}(t) = \pmb{s}\Bigr\},
\end{equation}
where, for $s\in\mathscr{S}_{i,j}$ and $\pmb{a}\in\mathscr{A}_{i,j}$,
\begin{equation}\label{eqn:define_c_gamma}
    C^{\pmb{\gamma}_{i,j}}_{i,j}(s,\pmb{a},t)\\\coloneqq c_{i,j}(s,\mathcal{e}(\pmb{a}),t)+ \gamma_{i,j,t}g_{i,j}(s) - \sum_{i'\in\bar{\mathscr{B}}^{-1}_{i,j}}\gamma_{i',j,t}\mathcal{w}_{i,i',j}a_{i'}.
\end{equation}
Equation~\eqref{eqn:dual_func} can be rewritten as
\begin{equation}\label{eqn:dual_func:2}
L(\pmb{\gamma}) = \min\nolimits_{\phi\in\tilde{\Phi}} \sum_{i\in[I]} \sum_{j\in[J]} L^{\phi}_{i,j}(\pmb{\gamma}_{i,j}),
\end{equation}
where $\pmb{\gamma}=(\gamma_{i,j,t}:i\in[I],j\in[J],t\in[T])$ with $\gamma_{i,j,t}\in\mathbb{R}$ being the Lagrange multiplier for \eqref{eqn:constraint:neighborhood:relax}, and $\pmb{\gamma}_{i,j}= (\gamma_{i',j,t}: i'\in \{i\}\cup\bar{\mathscr{B}}^{-1}_{i,j},t\in[T])$.

Following the idea of Whittle relaxation in \cite{whittle1988restless}, there is no constraint that restricts the value of each $L^{\phi}_{i,j}( \pmb{\gamma}_{i,j})$ on the right hand side of \eqref{eqn:dual_func:2} once the other $L^{\phi}_{i',j'}(\pmb{\gamma}_{i',j'})$ for $(i',j')\neq (i,j)$ are known. 
More precisely, consider the following proposition.
\begin{proposition}\label{prop:decomposition}
For given $\pmb{\gamma}\in\mathbb{R}^{IJT}$,
\begin{equation}\label{eqn:dual_func:3}
L(\pmb{\gamma}) =\min_{\phi\in\tilde{\Phi}} \sum_{i\in[I]} \sum_{j\in[J]} L^{\phi}_{i,j}(\pmb{\gamma}_{i,j})= \sum_{i\in[I]}\sum_{j\in[J]}\min_{\phi_{i,j}\in\tilde{\Phi}}   L^{\phi_{i,j}}_{i,j}(\pmb{\gamma}_{i,j}),
\end{equation}
where $\pmb{\gamma}_{i,j}= (\gamma_{i',j,t}: i'\in \{i\}\cup\bar{\mathscr{B}}^{-1}_{i,j},t\in[T])$.
\end{proposition}
The proof of Proposition~\ref{prop:decomposition} is provided in Appendix~\ref{app:prop:decomposition}.

Observe that the minimization over $\phi\in\tilde{\Phi}$ for each $L^{\phi}_{i,j}(\pmb{\gamma}_{i,j})$ in \eqref{eqn:dual_func:3} is equivalent to the minimization over variables  
\begin{equation}\label{eqn:define_alpha}
    \alpha^{\phi}_{i,i',j}(s,t)\coloneqq \mathbb{E}\Bigl[a^{\phi}_{i,i',j}(\bm{S}^{\phi}_{i,j}(t),t)~|~S^{\phi}_{i,j}(t)=s\Bigr]
\end{equation}
in $[0,1]$ for all $i\in[I]$, $s\in\mathscr{S}_{i,j}$ and $t\in[T]$. The value of $\alpha^{\phi}_{i,i',j}(s,t)$ is interpreted as the probability of taking the action $i'\in\mathscr{B}_{i,j}$
at time $t$ under policy $\phi$, given $S^{\phi}_{i,j}(t) = s$. 
The minimization over $\phi\in\tilde{\Phi}$ in the dual function can be decomposed into $IJ$ independent sub-problems.
We refer to $\min_{\phi\in\tilde{\Phi}}L^{\phi}_{i,j}(\pmb{\gamma}_{i,j})$ as the \emph{sub-problem}~$(i,j)$.
Compared to the original MAB-ML problem, the state space of each sub-problem has been significantly reduced. 

For $(i,j)\in[I]\times[J]$ and given $\pmb{\gamma}_{i,j}\in\mathbb{R}^{\lvert\{i\}\cup\bar{\mathscr{B}}^{-1}_{i,j}\rvert T}$, define $V^{\pmb{\gamma}_{i,j}}_{i,j}(s,t)$ ($s\in\mathscr{S}_{i,j},t\in[T]$)
as the solution of the Bellman equations
\begin{equation}\label{eqn:subproblem:bellman}
V^{\pmb{\gamma}_{i,j}}_{i,j}(s,t) = \min_{\pmb{a}\in\mathscr{A}_{i,j}}\Bigl\{C^{\pmb{\gamma}_{i,j}}_{i,j}(s,\pmb{a},t) + \sum_{s'\in\mathscr{S}_{i,j}}P_{i,j,t}(s,\mathcal{e}(\pmb{a}),s')V^{\pmb{\gamma}_{i,j}}_{i,j}(s',t+1)\Bigr\},~\forall s\in\mathscr{S}_{i,j},t\in[T],
\end{equation}
where, for $\pmb{a}=(a_{i'}: i'\in[I])\in\mathscr{A}_{i,j}$,  $\mathcal{e}(\pmb{a})= 1-\prod_{i'\in\mathscr{B}_{i,j}}(1-a_{i'})$ represents the probability of taking a non-zero action,
 $V^{\pmb{\gamma}_{i,j}}_{i,j}(s,T+1)\equiv 0$ for all $s\in\mathscr{S}_{i,j}$, and  $P_{i,j,t}(s,\mathcal{e}(\pmb{a}),s')$ are the transition probabilities from state $s$ to $s'$.
For $\mathcal{e}(\pmb{a}) =0$ or $1$, 
the transition probabilities $P_{i,j,t}(s,\mathcal{e}(\pmb{a}),s')$ are given values.
While $P_{i,j,t}(s,\mathcal{e}(\pmb{a}),s')$ for $\mathcal{e}(\pmb{a}) \in (0,1)$ is a linear combination of $P_{i,j,t}(s,0,s') $ and $P_{i,j,t}(s,1,s')$; that is, for $\mathcal{e}(\pmb{a})=e\in(0,1)$, $P_{i,j,t}(s,e,s') = e P_{i,j,t}(s,1,s') + (1-e)P_{i,j,t}(s,0,s')$.
We can solve \eqref{eqn:subproblem:bellman} for all $(i,j)\in[I]\times[J]$ independently.
From \cite{ross1992applied}, there exists an optimal solution $\phi^*$ for sub-problem $(i,j)$ meeting the condition that, for any $s\in\mathscr{S}_{i,j}$ and $t\in[T]$, $\alpha^{\phi^*}_{i,i',j}(s,t) = a_{i'}$ where $\pmb{a}\in\mathscr{A}_{i,j}$ minimizes the right hand side of \eqref{eqn:subproblem:bellman}. We refer to $V^{\pmb{\gamma}_{i,j}}_{i,j}(s,t)$ ($s\in\mathscr{S}_{i,j}$ and $t\in[T+1]$) as the \emph{value function} of sub-problem $(i,j)$ for given $\pmb{\gamma}_{i,j}$.

The dual version of the relaxed problem described in \eqref{eqn:obj:relax} and \eqref{eqn:constraint:neighborhood:relax} is
\begin{equation}\label{eqn:dual_problem:1}
\max_{\pmb{\gamma}\in\mathbb{R}^{IJT}} L(\pmb{\gamma}) 
= \max_{\pmb{\gamma}\in\mathbb{R}^{IJT}} \sum_{i\in[I]}\sum_{j\in[J]}\min_{\phi\in\tilde{\Phi}}L^{\phi}_{i,j}(\pmb{\gamma}_{i,j}),
\end{equation}
where the equality is achieved based on \eqref{eqn:dual_func:3}. For $(i,j)\in[I]\times[J]$ and given probability $\pi^0_{i,j}(s)$ of the initial state $S^{\phi}_{i,j}(1)=s$, $\min_{\phi\in\tilde{\Phi}}L^{\phi}_{i,j}(\pmb{\gamma}_{i,j}) = \sum_{s\in\mathscr{S}_{i,j}}\pi^0_{i,j}(s) V^{\pmb{\gamma}_{i,j}}_{i,j}(s,1)$.
Since $V^{\pmb{\gamma}_{i,j}}_{i,j}(s,t)$ ($s\in\mathscr{S}_{i,j}$ and $t\in[T]$) are the solutions of the equations in \eqref{eqn:subproblem:bellman}, the problem in \eqref{eqn:dual_problem:1} achieves the same maximum as 
\begin{equation}\label{eqn:dual_problem:2:1}
\max\nolimits_{\begin{subarray}
~\pmb{\gamma}\in \mathbb{R}^{IJT},\\
\bm{V}\in\mathbb{R}^{N}\end{subarray}}
\sum_{i\in[I]} \sum_{j\in[J]}\sum_{s\in\mathscr{S}_{i,j}}\pi^0_{i,j}(s)V_{i,j}(s,1) 
\end{equation}
subject to
\begin{multline}\label{eqn:dual_problem:2:2}
V_{i,j}(s,t) \leq C^{\pmb{\gamma}_{i,j}}_{i,j}(s,\pmb{\alpha},t) + \sum_{s'\in\mathscr{S}_{i,j}}P_{i,j,t}(s,\mathcal{e}(\pmb{\alpha}),s')V_{i,j}(s',t+1),\\
~\forall i\in[I],j\in[J], t\in[T], s\in\mathscr{S}_{i,j}, 
\pmb{\alpha}\in \tilde{\mathscr{A}}_{i,j},
\end{multline}
and
\begin{equation}\label{eqn:dual_problem:2:3}
V_{i,j}(s,T+1) = 0,~\forall i\in[I],j\in[J],s\in\mathscr{S}_{i,j}, 
\end{equation}
where $\tilde{\mathscr{A}}_{i,j}\coloneqq \{0,1\}^I\cap \mathscr{A}_{i,j}$, $\pi^0_{i,j}(s)$ is the given probability of the initial state $s\in\mathscr{S}_{i,j}$, $V_{i,j}(s,t)$ ($i\in[I],j\in[J],t\in[T],s\in\mathscr{S}_{i,j}$) and $\pmb{\gamma}$ are the unknown variables, and $\bm{V}\coloneqq (V_{i,j}(s,t):i\in[I],j\in[J],s\in\mathscr{S}_{i,j},t\in[T+1])$ with cardinality $P= (T+1)\sum_{(i,j)\in[I]\times[J]}|\mathscr{S}_{i,j}|$.  
Let $\pmb{\gamma}^*\in\mathbb{R}^{IJT}$ and $\bm{V}^*\in\mathbb{R}^{P}$ represent an optimal solution of \eqref{eqn:dual_problem:2:1}-\eqref{eqn:dual_problem:2:3}. For any $(i,j)\in[I]\times[J]$, $t\in[T]$ and $s\in\mathscr{S}_{i,j}$, $V^*_{i,j}(s,t) = V^{\pmb{\gamma}^*_{i,j}}_{i,j}(s,t)$, where $\pmb{\gamma}^*_{i,j} = (\gamma^*_{i',j,t}:i'\in \{i\}\cup \bar{\mathscr{B}}^{-1}_{i,j},t\in[T])$ with $\gamma^*_{i',j,t}$ the corresponding elements of $\pmb{\gamma}^*$.

The linear optimization described in \eqref{eqn:dual_problem:2:1}-\eqref{eqn:dual_problem:2:3}, with $(T+1)\sum_{(i,j)\in[I]\times[J]} |\mathscr{S}_{i,j}| +IJT$ unknown variables, can be solved by conventional linear programming methods, such as the interior point methods \cite{cohen2021solving,illes2002pivot}.
Note that this linear optimization problem differs from the original MAB-ML problem in \eqref{eqn:obj}, \eqref{eqn:constraint:exclusive} and \eqref{eqn:constraint:canonical}, whose optimal solution is in general not applicable to the original MAB-ML problem.
Instead, solving the linear problem is an intermediate step for proposing a near-optimal scheduling policy for the original MAB-ML.
We provide in Section~\ref{subsec:policies} the detailed steps of proposing such scheduling policies and, in Section~\ref{subsec:complexity}, analyze their computational complexities.

\subsection{Strong Duality}\label{subsec:strong_duality}
We prove in this subsection that the complementary slackness gap between optimality of the primary and dual problems becomes zero under a threshold policy.
That is, this threshold policy achieves optimality of the relaxed problem described in \eqref{eqn:obj:relax} and \eqref{eqn:constraint:neighborhood:relax}.
The threshold policy is in general not applicable to the original patrol problem described in \eqref{eqn:obj}, \eqref{eqn:constraint:exclusive} and \eqref{eqn:constraint:canonical}, but it quantifies marginal costs of moving agents to different directions. 
We start with solving sub-problem $(i,j)\in[I]\times[J]$. For $(i,j)\in[I]\times[J]$, $t\in[T]$,  $i'\in\mathscr{B}_{i,j}$, and  $\gamma\in\mathbb{R}$, define
\begin{equation}\label{eqn:define_theta}
\theta_{i,i',j,t}(\gamma)  \coloneqq 
\left\{\begin{cases}
    \gamma
    \mathcal{w}_{i,i',j}, & \text{if }i'\in\mathscr{B}_{i,j}\cap \bar{\mathscr{B}}^{-1}_{i,j},\\
    0, & \text{otherwise}.
\end{cases}\right.
\end{equation}
\begin{lemma}\label{lemma:indexability}
For $(i,j)\in[I]\times[J]$ and $\pmb{\gamma}_{i,j}\in\mathbb{R}^{\lvert\{i\}\cup\bar{\mathscr{B}}^{-1}_{i,j}\rvert T}$, a policy $\phi_{i,j}(\pmb{\gamma}_{i,j})\in\tilde{\Phi}$ is optimal to $\min\nolimits_{\phi_{i,j}\in \tilde{\Phi}} L^{\phi}_{i,j}(\pmb{\gamma}_{i,j})$ (sub-problem $(i,j)$) if, for $t\in[T]$, $i'\in\mathscr{B}_{i,j}$ and $s\in\mathscr{S}_{i,j}$,
\begin{equation}\label{eqn:lemma:indexability}
\alpha^{\phi_{i,j}(\pmb{\gamma}_{i,j})}_{i,i',j}(s,t)
\left\{\begin{array}{ll}
=1, & 
\text{if } \theta_{i,i',j,t}(\gamma_{i',j,t}) 
> \vartheta^{\pmb{\gamma}_{i,j}}_{i,j}(s,t) \text{ and } 
i' = \min\arg\max_{i''\in\mathscr{B}_{i,j}}\theta_{i,i'',j,t}(\gamma_{i'',j,t}),\\
\in[0,1], & \text{if } \theta_{i,i',j,t}(\gamma_{i',j,t}) = \vartheta^{\pmb{\gamma}_{i,j}}_{i,j}(s,t) \text{ and } 
i' \in \arg\max_{i''\in\mathscr{B}_{i,j}}\theta_{i,i'',j,t}(\gamma_{i'',j,t}), \\
=0, & \text{otherwise},
\end{array}\right.
\end{equation}
where 
\begin{equation}\label{eqn:vartheta}
     \vartheta^{\pmb{\gamma}_{i,j}}_{i,j}(s,t) \coloneqq c_{i,j}(s,1,t)-c_{i,j}(s,0,t)\\+\sum_{s'\in\mathscr{S}_{i,j}}\Bigl(P_{i,j,t}(s,1,s')-P_{i,j,t}(s,0,s')\Bigr)V^{\pmb{\gamma}_{i,j}}_{i,j}(s',t+1).
\end{equation}
\end{lemma}
The proof of Lemma~\ref{lemma:indexability} is provided in Appendix~\ref{app:lemma:indexability}.
In the special case for multi-agent patrolling, for $(i,j)\in[I]\times[J]$ and $i'\in\mathscr{B}_{i,j}$, since $\mathscr{B}_{i,j}=\bar{\mathscr{B}}^{-1}_{i,j}$ and $\mathcal{w}_{i,i',j}=1$, $\theta_{i,i',j,t}(\gamma_{i',j,t}) = \gamma_{i',j,t}$.

Lemma~\ref{lemma:indexability} implies that the optimal solution for each sub-problem exists in a threshold form \eqref{eqn:lemma:indexability}, similar to the \emph{Whittle indexability} for an RMAB problem.
For a continuous-time RMAB problem, the existence of Whittle indexability remains an open question in general, which has been proved under non-trivial conditions in \cite{nino2001restless,nino2007dynamic,ninomora2020verification}.
For the discrete-time case, threshold-form solutions exist in general for sub-problems of RMABs~ \cite{brown2020index}.
As explained in Section~\ref{subsec:MAB-ML},
if $J=1$, $\mathscr{B}_{i,j}=\{1\}$ for all $(i,j)\in[I]\times[J]$, and $\bar{\mathscr{B}}_{i,j} = \emptyset$ and $g_{i,j}(\cdot) = 0$ for all $(i,j)\in[I]\times[J]$ except $\bar{\mathscr{B}}_{1,j} = [I]$, $g_{1,j}(\cdot) = M$, and $\mathcal{w}_{i',1,j}=1$ for all $i'\in[I]$, then MAB-ML reduces to RMAB, and process $\{S^{\phi}_{i,j}(t), t\in[T]\}$ becomes a restless bandit process.
In this special case, the conclusion of Lemma~\ref{lemma:indexability} reduces to the \emph{modified indexability} discussed in~\cite{brown2020index}.
Nonetheless, as explained in Section~\ref{subsec:MAB-ML}, MAB-ML in general is not RMAB, nor can the definition of Whittle indexability in \cite{whittle1988restless} or the related results, such as \cite{nino2001restless,nino2007dynamic,ninomora2020verification,brown2020index}, be applicable here.

For MAB-ML, we have the following definition about having a threshold-form optimal solution,  referred to as \emph{indexability}.

\begin{definition}{Indexability}\label{definition:indexability}
For $(i,j)\in[I]\times[J]$, we say that process $\{S^{\phi}_{i,j}(t),t\in[T]\}$ is indexable if, for given $\pmb{\gamma}_{i,j}\in\mathbb{R}^{\lvert\{i\}\cup\bar{\mathscr{B}}^{-1}_{i,j}\rvert T}$, there exists an optimal solution $\phi_{i,j}(\pmb{\gamma}_{i,j})$ for the sub-problem $\min_{\phi\in\tilde{\Phi}}L^{\phi}_{i,j}(\pmb{\gamma}_{i,j})$ such that \eqref{eqn:lemma:indexability} and \eqref{eqn:vartheta} are satisfied for all $s\in\mathscr{S}_{i,j}$, $i'\in\mathscr{B}_{i,j}$ and $t\in[T]$.

In particular, for $t\in[T]$, given $\pmb{\gamma}\in\mathbb{R}^{IJT}$ and $(i,j)\in[I]\times[J]$, we say that process $\{S^{\phi}_{i,j}(t),t\in[T]\}$ is \emph{strictly indexable} when, for any $i_1,i_2\in\mathscr{B}_{i,j}$, $s_1\in\mathscr{S}_{i_1,j}$ and $s_2\in\mathscr{S}_{i_2,j}$, $\vartheta^{\pmb{\gamma}_{i_1,j}}_{i_1,j}(s_1,t)=\vartheta^{\pmb{\gamma}_{i_2,j}}_{i_2,j}(s_2,t)$ if  $i_1=i_2$ and $s_1=s_2$ where $\pmb{\gamma}_{i,j}=(\gamma_{i',j,t}:i'\in\{i\}\cup\bar{\mathscr{B}}^{-1}_{i,j},t\in[T])$ with $\gamma_{i',j,t}$ the corresponding elements of $\pmb{\gamma}$.

If process $\{S^{\phi}_{i,j}(t),t\in[T]\}$ is (strictly) indexable for all $(i,j)\in[I]\times[J]$, we say the MAB-ML process is (strictly) indexable.
\end{definition}

Unlike conventional RMAB problems that involve only one Lagrange multiplier, the threshold-form policy described in \eqref{eqn:lemma:indexability} decides its action variables for process $\{S^{\phi}_{i,j}(t),t\in[T]\}$ by comparing $\vartheta^{\pmb{\gamma}_{i,j}}_{i,j}(s,t)$ with multiple multipliers associated with all the elements in $\mathscr{B}_{i,j}$. 
The nested relationship between the many Lagrange multipliers and the action variables complying \eqref{eqn:lemma:indexability} negatively affects the clarity about how the threshold-form policy leads to an optimal solution of the relaxed problem or, more importantly, a near-optimal policy to the original MAB-ML.

In Proposition~\ref{prop:strong_duality}, we prove that a linear combination of policies that satisfy~\eqref{eqn:lemma:indexability} and \eqref{eqn:vartheta} achieves optimality of the relaxed problem. Later in Section~\ref{sec:policies}, we explain how the linear combination, acting as a policy to the relaxed problem (that is, satisfy Constraint~\eqref{eqn:constraint:exclusive} and \eqref{eqn:constraint:neighborhood:relax}), leads to a near-optimal policy to MAB-ML with proved asymptotic optimality.

\begin{proposition}\label{prop:strong_duality}
There exist $M\in\mathbb{N}_+$, policies $\phi^*_1,\phi^*_2,\ldots,\phi^*_M\in\tilde{\Phi}$ that satisfy \eqref{eqn:lemma:indexability} and \eqref{eqn:vartheta} for $\pmb{\gamma}=\pmb{\gamma}^*$, 
and a probability vector $\bm{\pi}^*\in[0,1]^M$ such that,
for $i\in[I]$, $j\in[J]$ and $t\in[T]$,
\begin{equation}\label{eqn:prop:strong_duality:6}
\sum_{m\in[M]}\pi^*_m\biggl(\sum_{i'\in\bar{\mathscr{B}}_{i,j}}\mathcal{w}_{i',i,j}\mathbb{E}\Bigl[\alpha^{\phi^*_m}_{i',i,j}\bigl(S^{\phi^*_m}_{i',j}(t),t\bigr)\Bigr] \\- \mathbb{E}\Bigl[g_{i,j}\bigl(S^{\phi^*_m}_{i,j}(t)\bigr)\Bigr]\biggr)
=0,
\end{equation}
where $\pmb{\gamma}^*\in\mathbb{R}^{IJT}$ are the optimal dual variables of the relaxed problem (an optimal solution of \eqref{eqn:dual_problem:2:1}-\eqref{eqn:dual_problem:2:3}).
\end{proposition}
The proof of Proposition~\ref{prop:strong_duality} is provided in Appendix~\ref{app:prop:strong_duality}.

Based on Proposition~\ref{prop:strong_duality}, optimality of the relaxed problem is achieved by a linear combination of policies $\phi^*_1,\phi^*_2,\ldots,\phi^*_M\in\tilde{\Phi}$, which are in the form of \eqref{eqn:lemma:indexability} and optimal to all the sub-problems. 
More precisely, consider a policy $\phi^*$ that tags each process $\bigl\{S^{\phi}_{i,j}(t),t\in[T]\bigr\}$ for $(i,j)\in[I]\times[J]$ with a number $m$ randomly selected from $[M]$, for which the probability of selecting $m$ is $\pi_m$. If a process $\bigl\{S^{\phi}_{i,j}(t),t\in[T]\bigr\}$ is tagged with $m\in[M]$, then apply policy $\phi^*_m$ to it. 
From Proposition~\ref{prop:strong_duality}, this policy $\phi^*$ achieves the minimum of the relaxed problem described in \eqref{eqn:obj:relax} and \eqref{eqn:constraint:neighborhood:relax}.
Note that $\phi^*$ is usually not applicable to the original MAB-ML, because it does not necessarily comply with Constraint~\eqref{eqn:constraint:canonical}.

\section{Scheduling Policies and Computational Complexity Analysis}\label{sec:policies}


\subsection{Scheduling Policies}\label{subsec:policies}
For the original MAB-ML, when $t\in[T]$ and $S^{\phi}_{i,j}(t) = s\in\mathscr{S}_{i,j}$, consider a policy $\phi\in\Phi$ that prioritizes the action  $\pmb{a}^{\phi}_{i,j}(\bm{S}^{\phi}(t),t)) = \pmb{1}^I_{i'}$ with $(i,j)\in[I]\times[J]$ and $i'\in\mathscr{B}_{i,j}$ according to an ascending order of
\begin{equation}\label{eqn:index:1}
\eta_{i,i',j}(s,t) \coloneqq 
\left\{\begin{cases}
    \vartheta^{\pmb{\gamma}^*_{i,j}}_{i,j}(s,t)- \theta_{i,i',j,t}(\gamma^*_{i',j,t}),&
    \text{if }i'\in\mathscr{B}_{i,j},\\
    +\infty,&\text{otherwise}
\end{cases}\right.
\end{equation}
where $\vartheta^{\pmb{\gamma}}_{i,j}(s,t)$ is defined in \eqref{eqn:vartheta}, and $\gamma^*_{i,j,t}$ and $V^*_{i,j}(s,t)$ ($(i,j)\in[I]\times[J], t\in[T],s\in\mathscr{S}_{i,j}$) are the optimal solutions of the dual problem described in \eqref{eqn:dual_problem:2:1}-\eqref{eqn:dual_problem:2:3}. 
The value $\eta_{i,i',j}(s,t)$ represents the marginal cost of taking the corresponding action, and, along with the tradition of the RMAB community, we refer to it as the \emph{index} assigned to the action. All the parameters on the right hand side of \eqref{eqn:index:1} are known \emph{a priori} or computable. Let $\pmb{\eta}\coloneqq (\eta_{i,i',j}(s,t): (i,i',j)\in[I]^2\times[J],s\in\mathscr{S}_{i,j},t\in[T])$.

The values of the indices $\pmb{\eta}$ can be calculated in an offline manner or updated online. In fine-tuned situations, we can also approximate the index values through learning techniques, such as Q-learning \cite{fu2019towards} and the upper-confidence-bound (UCB) algorithm \cite{wang2023optimistic}.

\subsubsection{Index Policy}\label{subsubsec:index_policy}
We propose an \emph{index policy} based on the action priorities imposed by the indices.
Let $a^{\rm IND}_{i,i',j}(\bm{S}^{\rm IND}(t),t)$ ($(i,i',j)\in[I]^2\times [J]$, $t\in[T]$) represent the action variables of the index policy. 
Define an \emph{action} as a tuple $(i,i',j)$, which, in the patrol case, represents a movement of a type-$j$ agent from area $i'$ to $i$.
We refer to such a tuple as the \emph{action} or \emph{movement} $(i,i',j)$.
For $j\in[J]$, $\phi\in\Phi$ and $t\in[T]$, maintain a set of actions (movements) $\mathscr{M}^{\phi}_j(t)\coloneqq \bigl\{(i,i',j)\in[I]^2| ~g_{i',j}(S^{\phi}_{i',j}(t))>0, i'\in\mathscr{B}_{i,j}\bigr\}$, representing the set of possible actions (movements) under policy $\phi$.
We will not take any of the other actions $(i,i',j)\notin \mathscr{M}^{\phi}_j(t)$ at time $t$ (that is, $a^{\phi}_{i,i',j}(\bm{S}^{\phi}_{i,i',j}(t),t) = 0$), because of \eqref{eqn:constraint:canonical}.
\begin{definition}{Action (Movement)
 Ranking}\label{define:movement_ranking}
For $\phi\in\Phi$, $t\in[T]$, and each $j\in[J]$, rank all the  actions (movements) $(i,i',j)\in\mathscr{M}^{\phi}_j(t)$ in the ascending order of their indices $\eta_{i,i',j}\bigl(S^{\phi}_{i,j}(t),t\bigr)$ at time $t$.
For the same $(i,j)\in[I]\times[J]$, if $i_1,i_2\in[I]$ with $i_1<i_2$ and $\eta_{i,i_1,j}\bigl(S^{\phi}_{i,j}(t),t\bigr)= \eta_{i,i_2,j}\bigl(S^{\phi}_{i,j}(t),t\bigr)$, then the action (movement) $(i,i_1,j)$ proceeds $(i,i_2,j)$.
Other tie cases can be considered in an arbitrary manner.
\end{definition}
Although different tie-breaking rules may affect the performance of the index policy in certain cases, all the theoretical results presented in this paper apply to arbitrary tie-breaking scenarios.

According to the \partialref{define:movement_ranking}{action ranking} for $\phi\in\Phi$, $j\in[J]$ and $t\in[T]$, we denote the rank of the action $(i,i',j)$ as $\mathcal{r}^{\phi}_j(i,i',t)$.
In this context, the index policy is such that, for $i,i'\in[I]$ and $j\in[J]$,
\begin{equation}\label{eqn:index_policy}
a^{\rm IND}_{i,i',j}\bigl(\bm{S}^{\rm IND}(t),t\bigr) \\=
\left\{\begin{cases}
    1,& \text{if } \sum\nolimits_{\begin{subarray}~i''\in \mathscr{B}_{i,j}: \\ \mathcal{r}^{\text{IND}}_j(i,i'',t) < \mathcal{r}^{\text{IND}}_j(i,i',t) \end{subarray}}\!\!\!\!\!\!\!\!\!\!\!\!\!\!\!\!\!a^{\rm IND}_{i,i'',j}\bigl(\bm{S}^{\rm IND}(t),t\bigr) = 0, \text{ and }
    \sum\nolimits_{\begin{subarray}~i''\in \bar{\mathscr{B}}_{i',j}: \\ \mathcal{r}^{\text{IND}}_j(i'',i',t) < \mathcal{r}^{\text{IND}}_j(i,i',t) \end{subarray}}< g_{i',j}(S^{\rm IND}_{i',j}(t)),
    \\
    0, & \text{otherwise.}
\end{cases}\right.
\end{equation}
For action $(i,i',j)$, the equality condition in \eqref{eqn:index_policy} requires 
that no action prior to $(i,i',j)$ is taken for process $\{S^{\text{IND}}_{i,j}(t),t\in[T]\}$ (constrained by \eqref{eqn:constraint:exclusive}). 
The inequality condition in \eqref{eqn:index_policy} 
is led by \eqref{eqn:constraint:canonical}.
If both conditions are satisfied, then we take action $(i,i',j)$.

We present in Algorithm~\ref{algo:index_policy} the pseudo-code of implementing the index policy, where, according to the above mentioned \partialref{define:moveement_ranking}{action ranking} for $j\in[J]$, $\phi$=IND and $t\in[T]$, we refer to the $\mathcal{r}$th movement as action $\mathcal{r}$ or action $(i_{\mathcal{r},j},i'_{\mathcal{r},j},j)$.

\IncMargin{1em}
\begin{algorithm}\small
\linespread{0.5}\selectfont

\SetKwProg{Fn}{Function}{}{End}
\SetKwInOut{Input}{Input}\SetKwInOut{Output}{Output}
\SetAlgoLined
\DontPrintSemicolon

\Input{Ranked actions for each $j\in[J]$ and $\bm{S}^{\rm IND}(t)$.}
\Output{$a^{\rm IND}_{i,i',j}\bigl(\bm{S}^{\rm IND}(t),t\bigr)$ for all $ (i,i',j)\in \bigcup_{j\in[J]}\mathscr{M}^{\text{IND}}_j(t)$.}

\Fn{IndexPolicy}{
	$a^{\rm IND}_{i,i',j}\bigl(\bm{S}^{\rm IND}(t),t\bigr)\gets 0$ for all $(i,i',j)\in\mathscr{M}^{\text{IND}}_j(t)$ and $j\in[J]$
 \tcc*{Initialization}\;
	\For{$j=1,2,\ldots,J$}{
            $q(i)\gets 0$ for all $i\in\bigl\{i| (i,i',j)\in\mathscr{M}^{\text{IND}}_j(t)\bigr\}$\;
            $p(i')\gets g_{i',j}(S^{\text{IND}}_{i',j}(t))$ for all $i'\in\bigl\{i'| (i,i',j)\in\mathscr{M}^{\text{IND}}_j(t)\bigr\}$\;
            \For{$\mathcal{r}=1,2,\ldots, |\mathscr{M}^{\text{IND}}_j(t)|$}{
	 	    \If{($i_{\mathcal{r}}\notin\bar{\mathscr{B}}_{i'_{\mathcal{r}},j}$ or $p(i'_{\mathcal{r}})  \geq \mathcal{w}_{i_{\mathcal{r}},i'_{\mathcal{r}},j}$)   and $q(i_{\mathcal{r},j})=0$\label{algo:availability}}{
                        $a^{\rm IND}_{i_{\mathcal{r},j},i'_{\mathcal{r},j},j}\bigl(\bm{S}^{\rm IND}(t),t\bigr) \gets 1$ \;
                        $q(i_{\mathcal{r},j}) \gets 1$\;
                        $p(i'_{\mathcal{r},j}) \gets p(i'_{\mathcal{r},j}) -\mathcal{w}_{i_{\mathcal{r}},i'_{\mathcal{r}},j}\mathds{1}\{i_{\mathcal{r}}\in\bar{\mathscr{B}}_{i'_{\mathcal{r}},j}\}$\;
                        
                }

            }
        }
        \Return $a^{\rm IND}_{i,i',j}\bigl(\bm{S}^{\rm IND}(t),t\bigr)$ for all $(i,i',j)\in\bigcup_{j\in[J]}\mathscr{M}^{\text{IND}}_j(t)$\;
    }
\caption{Pseudo-code for implementing the index policy.}\label{algo:index_policy}
\end{algorithm}
\DecMargin{1em}

Recall that $\pmb{\gamma}^*$ and $\bm{V}^*$ used for computing the indices in \eqref{eqn:index:1} are obtained by solving the linear program \eqref{eqn:dual_problem:2:1}-\eqref{eqn:dual_problem:2:3}, for which the number of unknown variables is linear in $T$, $I$, $J$ and $\sum_{(i,j)\in[I]\times[J]}|\mathscr{S}_{i,j}|$.
The complexity of ranking all the areas is $I\log I$, and, as described in Algorithm~\ref{algo:index_policy}, the complexity of implementing the index policy is linear in $I$ and $J$. 

We will prove in Section~\ref{subsec:asym_opt} that the index policy approaches optimality of the relaxed problem as the problem size tends to infinity and the sub-optimality diminishes exponentially in the problem size.

If MAB-ML reduces to RMAB (as explained in Seciton~\ref{subsec:MAB-ML}), then, in this very special case, the index policy satisfies \eqref{eqn:constraint:canonical} and is applicable to the original MAB-ML (RMAB). 
Nonetheless, in general, the index policy does not necessarily satisfy Constraint~\eqref{eqn:constraint:canonical} (or \eqref{eqn:constraint:neighbourhood}) - it may not be applicable to the original MAB-ML (or the original patrol problem).
For example, for some agents located in areas $i'$ with low $\gamma^*_{i,i',j}$ for all the $i\in\mathscr{B}_{i',j}$, their neighbourhoods may be fully occupied by others without leaving vacant areas for them.

In Section~\ref{sec:asym_opt}, we will provide a detailed discussion on asymptotically optimal policies, which are applicable to the general MAB-ML.
In Section~\ref{subsubsec:MAI} (the following subsection), we will focus on the patrol problem and propose a policy, derived from the index policy with adapted movements and applicable to the original patrol problem. 
In the scope of the patrol problem, such a policy is asymptotically optimal under a mild condition. The detail discussion for asymptotic optimality will also be provided in Section~\ref{sec:asym_opt}.

\subsubsection{Movement-Adapted Index Policy}\label{subsubsec:MAI}
In this subsection, we focus on the special case of MAB-ML - the scope of the patrol problem described in \eqref{eqn:obj}, \eqref{eqn:constraint:exclusive} and \eqref{eqn:constraint:neighbourhood}.
In this special case,  recall that, for all $(i,j)\in[I]\times[J]$, $\mathscr{B}_{i,j}=\bar{\mathscr{B}}_{i,j} = \bar{\mathscr{B}}^{-1}_{i,j}$, $\mathcal{w}_{i',i,j}=1$ for all $i'\in\bar{\mathscr{B}}_{i,j}$, $g_{i,j}(S^{\phi}_{i,j}(t))=\mathcal{i}(S^{\phi}_{i,j}(t)) = I^{\phi}_{i,j}(t)$ for any $t\in[T]$, and an action $(i,i',j)$ ($i\in\mathscr{B}_{i,j}$) represents a movement that moves a type-$j$ agent from area $i'$ to $i$.

To comply with Constraint~\eqref{eqn:constraint:neighbourhood}, we can adapt part of the movements determined by the index policy. 
In particular, 
for $t\in[T]$ and each $j\in[J]$, let $\mathscr{I}^{\phi}_j(t)$ represent a subset of $[I]$ for which $I^{\phi}_{i,j}(t)=1$ but \eqref{eqn:constraint:neighbourhood} is not satisfied under a policy $\phi$.
In other words, for any $i\in\mathscr{I}^{\rm IND}_j(t)$, all areas $i'\in\mathscr{B}_{i,j}$ have been determined to locate agents from other areas $i''\neq i$ at time $t$. 
For $i,i'\in[I]$, $j\in[J]$ and $t\in[T]$, given an action vector $\pmb{a}\in \prod_{(i,j)\in[I]\times[J]}\tilde{\mathscr{A}}_{i,j}$, 
if movement $(i,i',j)$ is taken with $a_{i,i',j}=1$, then define a variable $\varpi_{i,j}(\pmb{a}) \coloneqq i'$, representing the origin of the agent to move to area $i$ under action $\pmb{a}$; for all the other $(i,j)\in[I]\times[J]$ with $\sum_{i'\in\mathscr{B}_{i,j}}a_{i,i',j} = 0$, define $\varpi_{i,j}(\pmb{a})\coloneqq 0$. 
Let 
$\pmb{\varpi}(\pmb{a}) \coloneqq (\varpi_{i,j}(\pmb{a}): i\in[I], j\in[J])$.

For any $i\in\mathscr{I}^{\rm IND}_j(t)$, given the action vector $\pmb{a}^{\rm IND}(\bm{S}^{\rm IND}(t),t) = \pmb{a}$, we can select an area $i'\in\mathscr{B}_{i,j}$ but no longer take movement $(i',\varpi_{i',j}(\pmb{a}),j)$ determined by the index policy but replace it with movement $(i',i,j)$.  That is, the agent in area $i$ moves to area $i'$. 
In this way, the agent originally located in $\varpi_{i',j}(\pmb{a})$ needs to find another area in its neighbourhood to move to.
If area $\varpi_{i',j}(\pmb{a})$ coincidently has a vacant area in its neighbour, move the agent in $\varpi_{i',j}(\pmb{a})$ to the vacant place. 
Otherwise, we repeat the process of selecting an area in the neighbourhood of $\varpi_{i',j}(\pmb{a})$ and replacing the original movement by a new one.
We can keep replacing the original movements iteratively until reaching a vacant area. 
In this context, from the movements determined by the index policy, we can reach a new policy by iteratively replacing movements until Constraint~\eqref{eqn:constraint:neighbourhood} is satisfied.
We refer to this process of iteratively replacing the movements determined by the index policy until Constraint~\eqref{eqn:constraint:neighbourhood} are satisfied as the \emph{movement-adaption} process.

\begin{figure}[t]
\centering
\includegraphics[width=0.4\linewidth]{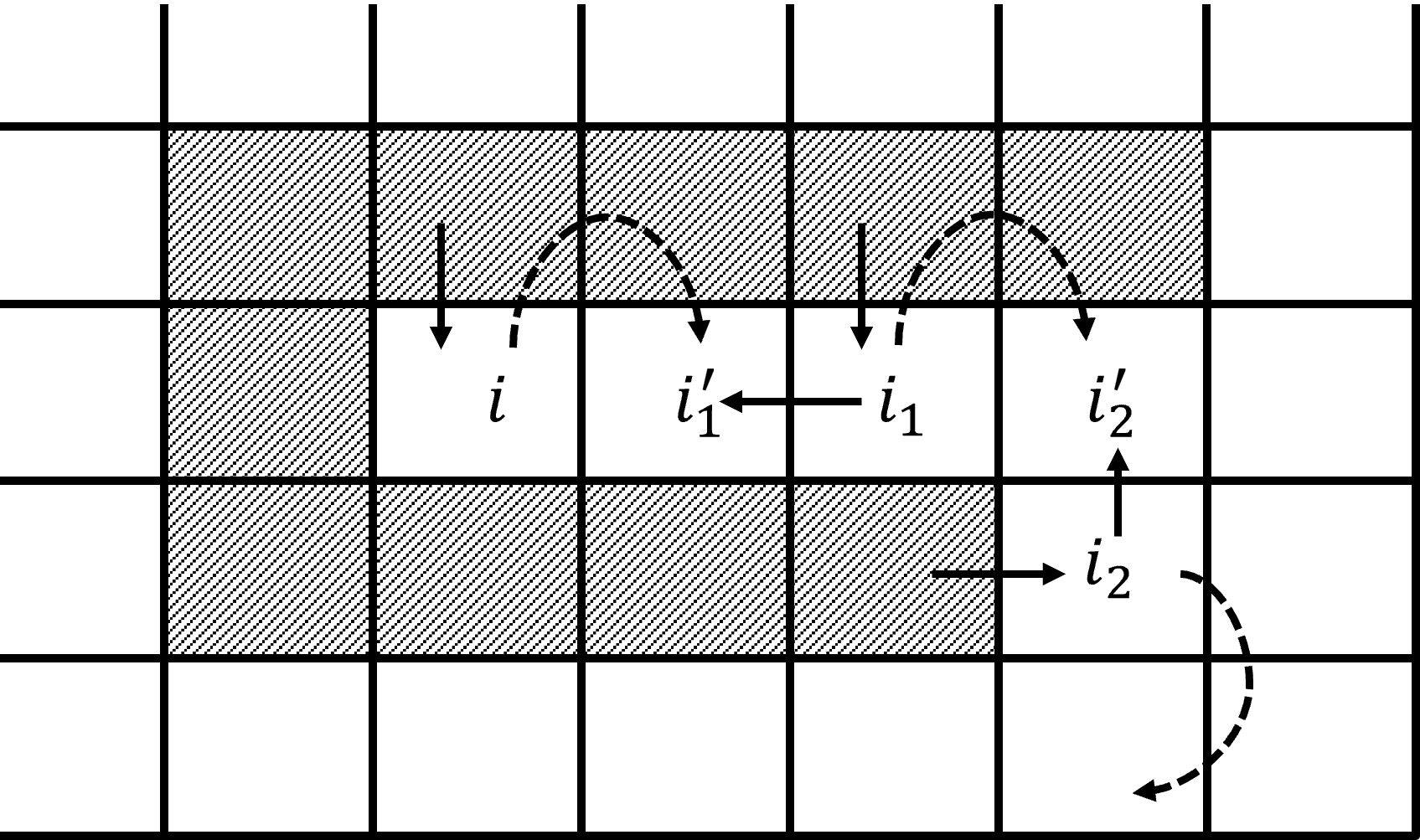}
\caption{An example of movement adaption. }\label{fig:movement_adaption}
\end{figure}

In Figure~\ref{fig:movement_adaption}, we provide a simple illustration for the movement adaption.
Areas $i,i_1$ and $i_2$ are explored by type-$j$ agents in time slot $t-1$ and, for time slot $t$, the index policy decides to move the agents in areas $i_1$ and $i_2$ to areas $i'_1$ and $i'_2$, respectively. In the figure, the solid arrows indicate the movements determined by the index policy, and the shadow parts are areas occupied by some other type-$j$ agents. Based on the index policy, all the neighbourhood of area $i$ have been occupied by some agents without a vacant place for the agent in area $i$ to stay in time $t$. We can adapt the movements by iteratively replacing $(i'_1,i_1,j)$ with $(i'_1,i,j)$ and $(i'_2,i_2,j)$ with $(i'_2,i_1,j)$, and then move the agent in area $i_2$ to any of its vacant neighbours. 
The dashed arrows represent the movements after the movement adaption.

We propose in Algorithm~\ref{algo:movement_adaption} the pseudo-code for a movement-adaption method and provide in the following the explanations of the proposed steps.
We refer to the policy described in Algorithm~\ref{algo:movement_adaption} as the \emph{movement-adapted index} (MAI) policy.
Let $a^{\rm MAI}_{i,i',j}(\pmb{s},t)$ ($i,i'\in[I],j\in[J],t\in[T],\pmb{s}\in\mathscr{S}$) represent the action variables of the MAI policy. 
For $i\in[I]$ and $j\in[J]$, let $d_{i,j}(\pmb{s},\pmb{\varpi})$ represent the \emph{distance} from area $i$ to a vacant area given the state $\bm{S}^{\rm MAI}(t)=\pmb{s}$ and 
$\pmb{\varpi}\Bigl(\pmb{a}^{\rm MAI}(\bm{S}^{\rm MAI}(t),t)\Bigr)=\pmb{\varpi}$ with respect to agent-type $j$; more precisely, define
\begin{equation}\label{eqn:distance}
    d_{i,j}(\pmb{s},\pmb{\varpi}) = \left\{\begin{cases}
        0,& \text{if }\mathcal{i}(s_{i,j}) = 1, \text{ and } \exists i'\in\mathscr{B}_{i,j} \text{ such that } \varpi_{i',j} = 0,\\       \min\nolimits_{\begin{subarray}~i'\in\mathscr{B}_{i,j}:\\\varpi_{i',j}\neq i\end{subarray}}d_{\varpi_{i',j}}+ 1, & \text{if }\mathcal{i}(s_{i,j}) = 1,  \text{ and }  \forall i'\in\mathscr{B}_{i,j},~\varpi_{i',j} > 0, \\
        \infty, &\text{otherwise}.
    \end{cases}\right.
\end{equation}
Recall that, for the patrol problem, we assume $|\mathscr{B}_{i,j}| \geq 1 $ for all $(i,j)\in[I]\times[J]$. 
In the second case in \eqref{eqn:distance}, for each $(i,j)\in[I]\times[J]$, the set $\{i'\in\mathscr{B}_{i,j}| \varpi_{i',j} \neq i\} \neq \emptyset$.

Given $\pmb{s}\in\mathscr{S}$ and $\pmb{\varpi}\in [I]^{IJ}$, we can initialize $d_{i,j}(\pmb{s},\pmb{\varpi}) = \infty$ for all $(i,j)\in[I]\times[J]$.
We then find all the $(i,j)\in[I]\times[J]$ for which $\mathcal{i}(s_{i,j})=1$ and there exists $i'\in\mathscr{B}_{i,j}$ such that $\varpi_{i',j} = 0$. 
For such $(i,j)$, update $d_{i,j}(\pmb{s},\pmb{\varpi})=0$. 
Then iteratively obtain $d_{i',j}(\pmb{s},\pmb{\varpi})$ for the direct neighbours $i'\in\mathscr{B}_{i,j}$ of area $i$ with previously determined $d_{i,j}(\pmb{s},\pmb{\varpi})<\infty$.

In Algorithm~\ref{algo:movement_adaption}, we  initialize $\pmb{a}^{\rm MAI}\bigl(\bm{S}^{\rm MAI}(t),t\bigr)$ with the action vector $\pmb{a}^{\rm IND}\bigl(\bm{S}^{\rm MAI}(t),t\bigr)$ by calling Algorithm~\ref{algo:index_policy}, where $\bm{S}^{\rm MAI}(t)$ is substituted for $\bm{S}^{\rm IND}(t)$. 
Based on the initialized $\pmb{a}^{\rm MAI}\bigl(\bm{S}^{\rm MAI}(t),t\bigr)=\pmb{a}$ and the given $\bm{S}^{\rm MAI}(t) = \pmb{s}$, we can calculate $\varpi_{i,j}(\pmb{a})$ and $d_{i,j}(\pmb{s},\pmb{\varpi}(\pmb{a}))$ for all $(i,j)\in[I]\times [J]$.

For each $j\in[J]$ and $i\in\mathscr{I}^{\rm MAI}_j(t)$, 
in Lines~\ref{line:while}-\ref{line:endWhile} of Algorithm~\ref{algo:movement_adaption}, we find a \emph{path} from area $i$ to a vacant area through the movement adaption.
In particular, we start with area $i$ and replace a movement $(i^*,\varpi,j)$ with $(i^*,i,j)$, for which $\varpi = \varpi_{i^*,j}\bigl(\pmb{a}^{\rm MAI}(\pmb{s},t)\bigr)$, $i^*$ is in $\mathscr{B}_{i,j}$, and area $\varpi$ has the smallest distance to a vacant area (see Line~\ref{line:select_neighbour}).
After this movement replacement, we find a place for the agent from area $i$ to stay in time slot $t$ while the agent in area $\varpi$ may have no place to go.
Then we focus on $\varpi$ in the next while loop with appropriately updated variables.
We will iteratively replace the selected movements until reaching a vacant area (that is, the stop condition in Line~\ref{line:reach_edge} is satisfied).
The output action variables $\pmb{a}^{\rm MAI}(\pmb{s},t)$ form the MAI policy and satisfy Constraint~\eqref{eqn:constraint:neighbourhood}. The MAI policy is applicable to the original patrol problem described in \eqref{eqn:obj}, \eqref{eqn:constraint:exclusive} and \eqref{eqn:constraint:neighbourhood}.

In Section~\ref{subsec:asym_opt}, we prove that, under a mild condition related to the fine-tuned values of the cost rates, when the problem size grows to infinity, the index policy, the MAI policy and the optimality of the relaxed problem coincide with each other. Since MAI is applicable to the original patrol problem, MAI also approaches optimality of the original patrol problem. 

\IncMargin{1em}
\begin{algorithm*}\small
\linespread{1.5}\selectfont

\SetKwProg{Fn}{Function}{}{End}
\SetKwInOut{Input}{Input}\SetKwInOut{Output}{Output}
\SetAlgoLined
\DontPrintSemicolon

\Input{The state vector $\bm{S}^{\rm MAI}(t) = \pmb{s}$ at time $t$.}
\Output{The action variables $a_{i,i',j}^{\rm MAI}\bigl(\bm{S}^{\text{MAI}}(t),t\bigr)$ for $(i,i',j)\in\bigcup_{j\in[J]}\mathscr{M}^{\text{MAI}}_j(t)$.}

\Fn{MAIPolicy}{
	$a_{i,i',j}^{\rm MAI}(\pmb{s},t)\gets a_{i,i',j}^{\rm IND}(\pmb{s},t)$  for $(i,i',j)\in\bigcup_{j\in[J]}\mathscr{M}^{\text{MAI}}_j(t)$\tcc*{Initialize the action variables by calling  Algorithm~\ref{algo:index_policy}.}
    Based on the action variables $a_{i,i',j}^{\rm MAI}(\pmb{s},t)$, initialize $\mathscr{I}^{\rm MAI}_j(t)$ for all $j\in[J]$, $\varpi_{i,j}\bigl(\pmb{a}^{\rm MAI}(\pmb{s},t)\bigr)$ for all $(i,j)\in[I]\times[J]$, and $d_{i,j}\Bigl(\pmb{s},\pmb{\varpi}\bigl(\pmb{a}^{\rm MAI}(\pmb{s},t)\bigr)\Bigr)$ described in \eqref{eqn:distance} for all $(i,j)\in[I]\times [J]$.\;
    $\varpi_{i,j}\gets \varpi_{i,j}\bigl(\pmb{a}^{\rm MAI}(\pmb{s},t)\bigr)$ for all $(i,j)\in[I] \times [J]$\;
    $d_{i,j}\gets d_{i,j}\Bigl(\pmb{s},\pmb{\varpi}\bigl(\pmb{a}^{\rm MAI}(\pmb{s},t)\bigr)\Bigr)$ for all $(i,j)\in[I] \times [J]$\;
	\For{$j=1,2,\ldots,J$}{
        \For{$\forall i\in\mathscr{I}^{\rm MAI}_j(t)$}{
            $\bar{i}\gets i$\;
            \While{$\bar{i} >0$}{\label{line:while}
                $i^*\gets \min\arg\min_{i'\in\mathscr{B}_{\bar{i},j}}d_{\varpi_{i',j},j}$\label{line:min}\;
                $\varpi \gets \varpi_{i^*,j}$ \label{line:select_neighbour}\; 
                $a^{\rm MAI}_{i^*,\bar{i},j}(\pmb{s},t)\gets 1$\;
                $a^{\rm MAI}_{i^*,\varpi,j}(\pmb{s},t)\gets 0$
                \tcc*{Replace movement $(i^*,\varpi,j)$ with $(i^*,\bar{i},j)$}
                \If {$d_{\varpi,j}=0$}{ \label{line:reach_edge}
                    Find an $i'\in\mathscr{B}_{\varpi,j}$ such that $\varpi_{i',j}= \infty$\;       
                    $a^{\rm MAI}_{i',\varpi,j}(\pmb{s},t)\gets 1$
                    \tcc*{Take movement $(i',\varpi,j)$}
                    $\varpi_{i',j}\gets \varpi$\;
                    $\varpi_{i^*,j}\gets \bar{i}$\;
                    $\bar{i}\gets 0$\;
                }\Else{
                    Exchange the values of variables $\varpi_{i^*,j}$ and $\bar{i}$\;
                }
            \label{line:endWhile}}
            Based on the updated $\pmb{\varpi}\coloneqq (\varpi_{i,j}: i\in[I], j\in[J])$, update the values of $d_{i',j}(\pmb{s},\pmb{\varpi})$ for all $i'\in[I]$ that satisfy \eqref{eqn:distance}.\label{line:update_distance}\;
            $d_{i',j}\gets d_{i'j}(\pmb{s},\pmb{\varpi})$ for all $i'\in[I]$\;
        }        
    }
    \Return $a_{i,i',j}^{\rm MAI}\bigl(\bm{S}^{\text{MAI}}(t),t\bigr)$ for $(i,i',j)\in\bigcup_{j\in[J]}\mathscr{M}^{\text{MAI}}_j(t)$\;
}
\caption{Pseudo-code for the movement adaption.}\label{algo:movement_adaption}
\end{algorithm*}
\DecMargin{1em}
\subsection{Computational Complexity of the Policies}\label{subsec:complexity}

\subsubsection{Computational complexity of the index policy}\label{subsubsec:complexity:index_policy}
The values of the indices $\eta_{i,i',j}(s,t)$ for all $i,i'\in[I]$, $j\in[J]$, $s\in\mathscr{S}_{i,j}$, and $t\in[T]$ are computed a priori in an offline manner. 
From the definition in \eqref{eqn:index:1}, the indices are in closed forms when the optimal solution $\pmb{\gamma}^*$ and $\bm{V}^*$ of the dual problem described in \eqref{eqn:dual_problem:2:1}-\eqref{eqn:dual_problem:2:3} are computed. 
Since the dual problem in \eqref{eqn:dual_problem:2:1}-\eqref{eqn:dual_problem:2:3} is a linear programming problem with $\bar{P}=(T+1)\sum\nolimits_{(i,j)\in[I]\times[J]}|\mathscr{S}_{i,j}| + IJT$ variables and $2T \sum\nolimits_{(i,j)\in[I]\times[J]}|\mathscr{S}_{i,j}| = \Omega(\bar{P})$ constraints, 
it can be solved by conventional linear programming methods, such as the interior point methods \cite{cohen2021solving,illes2002pivot} with complexity no worse than $O(\bar{P}^{2.5})$. 

As described in~\partialref{define:movement_ranking}{Action Ranking}, we rank the actions for each $j\in[J]$ according to the pre-computed $\pmb{\eta}$ and the states $S^{\text{MAI}}_{i,j}(t)$ ($(i,j)\in[I]\times[J]$) at time $t$. 
For $\phi\in\Phi$ and $t\in[T]$, define $\Psi^{\phi}(t) \coloneqq \sum_{j\in[J]}|\mathscr{M}^{\phi}_j(t)|$.
At time $t\in[T]$, for $j\in[J]$, the computational complexity to maintain the set $\mathscr{M}^{\text{MAI}}_j(t)$ and rank the $|\mathscr{M}^{\text{MAI}}_j(t)|$ actions is no worse than $O(\Psi^{\text{MAI}}(t))$ and $O(\Psi^{\text{MAI}}(t)\log \Psi^{\text{MAI}}(t))$, respectively.

The index policy can be constructed through the steps in Algorithm~\ref{algo:index_policy}, for which the computational complexity is $O(\Psi^{\text{MAI}}(t))$.
Taking $B \coloneqq\max_{i,j}|\mathscr{B}_{i,j}|$, we obtain that $\Psi^{\text{MAI}}(t) \leq \sum_{j\in[J]}M_jB$ for any $t\in[T]$.
Hence, the total computational complexity for ranking the actions and construct the index policy is $O(\sum_{j\in[J]}M_jB\log(\sum_{j\in[J]}M_jB))$.

\subsubsection{Computational complexity of MAI}
Recall that, as described in Section~\ref{subsubsec:MAI}, MAI is proposed for the patrol problem, a special case of MAB-ML, and is based on the indices defined in \eqref{eqn:index:1}.

With given indices $\pmb{\eta}$, the computational complexity for MAI consists of three cascading parts: ranking the movements (the actions for the general MAB-ML) based on the pre-computed indices $\pmb{\eta}$ (described in \eqref{eqn:index:1}), determining the action variables for the index policy with the ranked movements (pseudo-code in Algorithm~\ref{algo:index_policy}), and the movement adaption process (pseudo-code in Algorithm~\ref{algo:movement_adaption}).
All the three parts are processed online with observed system state $\bm{S}^{\text{MAI}}(t)$ at time $t$.
As alluded to above, the total complexity for the first two parts is $O(\sum_{j\in[J]}M_jB\log(\sum_{j\in[J]}M_jB))$.

The third part, movement adaption, is processed through Algorithm~\ref{algo:movement_adaption}, which includes initialization steps and  three levels of loops.
The computational complexity of the initialization steps is no worse than $O(IJB+\sum\nolimits_{j\in[J]}M_jB)$.
The external loop with variable $j\in[J]$ iterates $J$ times.
The second-level loop with variable $i\in\mathscr{I}^{\text{MAI}}_j(t)$ repeats $O(I)$ times. 
The internal loop from Line~\ref{line:while} to \ref{line:endWhile} iterates at most $d_{\bar{i},j}=O(I)$ times, each of which has complexity $O(B)$ mainly incurred by the minimum operation in Line~\ref{line:min}.
In Line~\ref{line:update_distance}, updating the area distances for each $j\in[J]$ costs $O(IB+M_jB)$. 
Since $M_j \leq I$, the worst case complexity is $O(IB)$.
Hence, the overall computational complexity of Algorithm~\ref{algo:movement_adaption} is at most $O(I^2JB)$.

For each time slot $t$, the computational complexity for implementing MAI (the above-mentioned three cascading parts) is 
$O\bigl(MB \log (MB)+I^2JB\bigr)$,
where $M=\sum_{j\in[J]}M_j$ representing the total number of agents, and $I$ and $J$ are the numbers of different areas and agent types, respectively.

\section{Asymptotic Optimality}\label{sec:asym_opt}
\subsection{Asymptotic Regime}\label{subsec:asym_regime}

For the patrol problem, consider a large number $h\in\mathbb{N}_+$ of highly sophisticated sensors that are used to precisely explore tiny \emph{sub-areas} located in a geographical region. In particular, we divide each area $i\in[I]$ into $h$ sub-areas, referred to as the sub-areas $(i,k)$ ($i\in[I]$, $k\in[h]$), each of which is related to $J$ stochastic processes $\{S^{\phi,k}_{i,j}(t), t\in[T]\}$ ($j\in[J]$) with state space $\mathscr{S}_{i,j}$. The transition probabilities of process $\{S^{\phi,k}_{i,j}(t): t\in[T]\}$ are the same as $\{S^{\phi}_{i,j}(t): t\in[T]\}$ defined in Section~\ref{sec:model}. 
The number of agents of type $j\in[J]$ is considered to be $M_j = h M_j^0$ for some $M_j^0\in\mathbb{N}_+$. 
Within a unit time, an agent in the sub-area $(i,k)$ ($i\in[I]$, $k\in[h]$) can move to any sub-area in the neighbourhood of area $i$ (that is, the agent can move to any sub-area $(i',k')$ with $i'\in\mathscr{B}_{i,j}$ and $k'\in[h]$). 
We refer to $h$ as the \emph{scaling parameter} of the patrol problem and refer to the problem with the scaling parameter as the \emph{patrol problem scaled by $h$}.
As $h$ increases, the patrol region is divided into increasingly many sub-areas explored by compatible numbers of different agents.
More sub-areas and agents indicate larger geographical region or the same region with higher precision of collected information and more careful exploration. The patrol system defined in Section~\ref{subsec:patrolling} is a special case of the scaled patrol problem with $h=1$.

For the general case, MAB-ML, the problem can also be scaled by $h\in\mathbb{N}_+$.
Similarly, for $(i,j)\in[I]\times[J]$, consider $h$ processes $\{S^{\phi,k}_{i,j}(t): t\in[T]\}$ ($k\in[h]$) that have the same state and action spaces, transition matrices under different actions, and cost functions as $\{S^{\phi}_{i,j}(t),t\in[T]\}$.
Note that the initial states of the $h$ processes can be different and follow a given probability distribution $\pmb{\pi}^0_{i,j} = (\pi^0_{i,j}(s):s\in\mathscr{S}_{i,j}) \in [0,1]^{|\mathscr{S}_{i,j}|}$.
We refer to these processes $\{S^{\phi,k}_{i,j}(t): t\in[T]\}$ ($k\in[h]$) as the sub-processes associated with $(i,j)$.
Let $\bm{S}^{\phi,h}(t)\coloneqq (S^{\phi,k}_{i,j}(t):i\in[I],j\in[j],k\in[h])$, which takes values in 
\begin{equation}
    \mathscr{S}^h\coloneqq \prod\nolimits_{(i,j)\in[I]\times[J]}\Bigl(\mathscr{S}_{i,j}\Bigr)^h.
\end{equation}
More precisely, for $i,i'\in[I]$, $j\in[J]$, $k\in[h]$, $\pmb{s}\in\mathscr{S}^h$, and $t\in[T]$, define action variable $a^{\phi,k}_{i,i',j}(\pmb{s},t)$ for process $\{S^{\phi,k}_{i,j}(t),t\in[T]\}$. 
It takes binary values: if $a^{\phi,k}_{i,i',j}(\pmb{s},t) = 1$ then, when $\bm{S}^{\phi,k}(t) = \pmb{s}$, action $(i,i',j)$ is taken; otherwise, this action is not taken. 
If $i'\in[I]\backslash \mathscr{B}_{i,j}$, then $a^{\phi,k}_{i,i',j}(\pmb{s},t)\equiv 0$.
For the patrolling case, action $(i,i',j)$ is interpreted as moving an agent of type $j$ from a sub-area of area $i'$ to sub-area $(i,k)$ at the beginning of time slot $t$ and explores there before the next time slot.
For $(i,j)\in[I]\times[J]$, $k\in[h]$, $\pmb{s}\in\mathscr{S}^h$, $\phi\in\Phi^h$ and $t\in[T]$, define $\pmb{a}^{\phi,k}_{i,j}(\pmb{s},t)\coloneqq (a^{\phi,k}_{i,i',j}(\pmb{s},t): i'\in[I])$.


In such a scaled general MAB-ML system, let $\Phi^h$ represent the set of all the policies determined by such action variables.
We generalize the MAB-ML problem described in \eqref{eqn:obj}, \eqref{eqn:constraint:exclusive} and \eqref{eqn:constraint:canonical} to its scaled version,
\begin{equation}\label{eqn:obj:h}
    \min_{\phi\in\Phi^h} \frac{1}{h}\sum_{t\in[T]}\sum_{i\in[I]}\sum_{j\in[J]}\sum_{k\in[h]}\mathbb{E}\Bigl[c_{i,j}\Bigl(S^{\phi,k}_{i,j}(t),\mathcal{e}^{\phi,k}_{i,j}\bigl(\bm{S}^{\phi,h}(t),t\bigr),t\Bigr)\Bigr],
\end{equation}
where $\mathcal{e}^{\phi,k}_{i,j}(\bm{S}^{\phi,h}(t),t)\coloneqq \mathds{1}\bigl\{\sum_{i'\in\mathscr{B}_{i,j}}a^{\phi,k}_{i,i',j}\bigl(\bm{S}^{\phi}(t),t\bigr) > 0\bigr\}$,
subject to
\begin{equation}\label{eqn:constraint:exclusive:h}
    \sum_{i'\in\mathscr{B}_{i,j}}a^{\phi,k}_{i,i',j}\Bigl(\bm{S}^{\phi,h}(t),t\Bigr) \leq 1,~\forall i\in[I],j\in[J],t\in[T],k\in[h],
\end{equation}
and
\begin{equation}\label{eqn:constraint:neighbourhood:h}
    \frac{1}{h}\sum_{i'\in\bar{\mathscr{B}}_{i,j}}\sum_{k\in[h]}\mathcal{w}_{i',i,j}a^{\phi,k}_{i',i,j}\bigl(\bm{S}^{\phi,h}(t),t\bigr) = \frac{1}{h}\sum_{k\in[h]} g_{i,j}\bigl(S^{\phi,k}_{i,j}(t)\bigr),\\~\forall i\in[I],j\in[J],t\in[T].
\end{equation}
The $1/h$ in \eqref{eqn:obj:h}, \eqref{eqn:constraint:exclusive:h} and \eqref{eqn:constraint:neighbourhood:h} is used to make all the expressions finite when $h\rightarrow +\infty$. 

In the patrol case, 
Constraint~\eqref{eqn:constraint:exclusive:h} guarantee that each sub-area is explored by at most one agent of the the same type at the same time slot.
With Constraint~\eqref{eqn:constraint:neighbourhood:h}, each agent must move to its neighbourhood, including staying in its current location, within a unit time. 

The MAB-ML problem described in \eqref{eqn:obj}, \eqref{eqn:constraint:exclusive} and \eqref{eqn:constraint:canonical} is a special case of the problem~\eqref{eqn:obj:h}, \eqref{eqn:constraint:exclusive:h} and \eqref{eqn:constraint:neighbourhood:h} by setting $h=1$.
We refer to the problem described in \eqref{eqn:obj:h}, \eqref{eqn:constraint:exclusive:h} and \eqref{eqn:constraint:neighbourhood:h} as the \emph{scaled} MAB-ML problem.
We refer to the case with $h\rightarrow +\infty$ as the \emph{asymptotic regime}.

Similar to the special case in Section~\ref{sec:relaxation}, we randomize the action variables $a^{\phi,k}_{i,i',j}(\pmb{s},t)$ and relax \eqref{eqn:constraint:neighbourhood:h} to
\begin{equation}\label{eqn:constraint:neighbourhood:h:relax}
    \frac{1}{h}\sum_{i'\in\bar{\mathscr{B}}_{i,j}}\sum_{k\in[h]}\mathcal{w}_{i',i,j}\alpha^{\phi,k}_{i',i,j}\bigl(S^{\phi,k}_{i',j}(t),t\bigr) \\= \frac{1}{h}\sum_{k\in[h]} \mathbb{E}\Bigl[g_{i,j}\bigl(S^{\phi,k}_{i,j}(t)\bigr)\Bigr],~\forall i\in[I],j\in[J],t\in[T],
\end{equation}
where $\alpha^{\phi,k}_{i,i',j}(s,t)\coloneqq \mathbb{E}\bigl[a^{\phi,k}_{i,i',j}\bigl(\bm{S}^{\phi,h}(t),t\bigr)|S^{\phi,k}_{i,j}(t) = s\bigr]$, representing the probability of taking $a^{\phi,k}_{i,i',j}\bigl(\bm{S}^{\phi,h}(t),t\bigr)=1$ when $S^{\phi,k}_{i,j}(t) = s$. Define $\tilde{\Phi}^h$ as the set of all the policies $\phi$ determined by such action variables $\alpha^{\phi,k}_{i,i',j}(s,t)$ for all $i,i'\in[I]$, $j\in[J]$, $k\in[h]$, $t\in[T]$ and $s\in\mathscr{S}_{i,j}$ with satisfied \eqref{eqn:constraint:exclusive:h}.
We refer to the problem
\begin{equation}\label{eqn:obj:h:relax}
\min_{\phi\in\tilde{\Phi}^h} \frac{1}{h}\sum_{t\in[T]}\sum_{i\in[I]}\sum_{j\in[J]}\sum_{k\in[h]} \mathbb{E}\Bigl[c_{i,j}\Bigl(S^{\phi,k}_{i,j}(t),\mathcal{e}^{\phi,k}_{i,j}\bigl(\bm{S}^{\phi,k}_{i,j},t\bigr),t),t\Bigr) \Bigr],
\end{equation}
subject to  \eqref{eqn:constraint:neighbourhood:h:relax} as the relaxed version of MAB-ML in the scaled system. 

The relaxed problem in the scaled system can be decomposed along the same lines as described in Section~\ref{sec:relaxation}, resulting in $IJh$ sub-problems that are associated with the $IJh$ processes $\Bigl\{S^{\phi,k}_{i,j}(t),t\in[T]\Bigr\}$ ($i\in [I]$, $j\in[J]$ and $k\in[h]$).
We still use $\pmb{\gamma}\in\mathbb{R}^{IJT}$ to represent the Lagrange multipliers for \eqref{eqn:constraint:neighbourhood:h:relax}.
More precisely, for $i\in[I]$, $j\in[J]$, $k\in[h]$, $\phi\in\tilde{\Phi}^h$, and $\pmb{\gamma}_{i,j}\in\mathbb{R}^{\lvert\{i\}\cup\bar{\mathscr{B}}^{-1}_{i,j}\rvert T}$, 
define 
\begin{equation}\label{eqn:define_func_i:h}
L^{\phi,k}_{i,j}( \pmb{\gamma}_{i,j}) \\\coloneqq \sum_{t\in[T]}
\sum_{s\in\mathscr{S}_{i,j}}\sum_{\pmb{a}\in\mathscr{A}_{i,j}}
C^{\pmb{\gamma}_{i,j}}_{i,j}(s,\pmb{\alpha}^{\phi,k}_{i,j}\bigl(s,t\bigr),t)  \mathbb{P}\Bigl[S^{\phi,k}_{i,j}(t) = s\Bigr],
\end{equation}
and the dual function of the relaxed problem in \eqref{eqn:obj:h:relax} and \eqref{eqn:constraint:neighbourhood:h:relax} is
\begin{equation}\label{eqn:dual_func:h}
    L^h(\pmb{\gamma}) = \frac{1}{h}\min_{\phi\in\tilde{\Phi}^h} \sum_{i\in[I]} \sum_{j\in[J]}\sum_{k\in[h]} L^{\phi,k}_{i,j}(\pmb{\gamma}_{i,j}) \\= \sum_{i\in[I]} \sum_{j\in[J]}\sum_{k\in[h]}\min_{\phi_{i,j}\in\tilde{\Phi}^h} L^{\phi_{i,j},k}_{i,j}(\pmb{\gamma}_{i,j}),
\end{equation}
where $\pmb{\gamma}=(\gamma_{i,j,t}:i\in[I],j\in[J],t\in[T])\in\mathbb{R}^{IJT}$ are the Lagrange multipliers, $\pmb{\gamma}_{i,j}=(\gamma_{i',j,t}:i'\in\{i\}\cup\bar{\mathscr{B}}^{-1}_{i,j},t\in[T])$ with $\gamma_{i',j,t}$ the corresponding elements of $\pmb{\gamma}$, and the second equality is a direct result of Proposition~\ref{prop:decomposition}.
Similar to the special case in Section~\ref{sec:relaxation}, we refer to $\min_{\phi_{i,j}\in\tilde{\Phi}^h} L^{\phi_{i,j},k}_{i,j}( \pmb{\gamma}_{i,j})$ as the sub-problem associated with $(i,j,k)\in[I]\times[J]\times[h]$. 
Given $(i,j)\in[I]\times[J]$, the sub-problems are all the same for $k\in[h]$.

It follows that the linear problem in \eqref{eqn:dual_problem:2:1}-\eqref{eqn:dual_problem:2:3} achieves the same minimum as the dual problem of the relaxed problem in the scaled system upon the same solution $\pmb{\gamma}^*$ and $\bm{V}^*$.
For $(i,j)\in[I]\times[J]$ and any $k\in[h]$, Lemma~\ref{lemma:indexability} and Proposition~\ref{prop:strong_duality} can be directly applied to the sub-problems in the scaled system with added superscript $k$ to the variables $\pmb{\alpha}^{\phi}_{i,j}(\cdot,\cdot)$, $S^{\phi}_{i,j}(t)$ and $L^{\phi}_{i,j}(\pmb{\gamma}_{i,j})$.
The index policy can be applied along similar lines as those described in Section~\ref{subsubsec:index_policy}.
For the scaled patrol problem, MAI defined in Section~\ref{subsubsec:MAI} is also directly applicable.
In Appendix~\ref{app:algorithms:h}, we provide the pseudo-codes for the index policy and MAI in the scaled system for the general MAB-ML and the patrol case (special case of MAB-ML), respectively.
The computational complexity of both policies can be analysed through the same steps, but is scaled by the parameter $h$ compared to the special case discussed in Section~\ref{subsec:complexity}. 
The complexity for implementing the index policy (including ranking the actions) and MAI (for the patrol case) are thus $O(M^0B\log (M^0B)+MB)$ and $O(MB\log(MB) + I^2JBh\bar{M})$, respectively, where $M=\sum\nolimits_{j\in[J]}M^0_j h$ is the total number of agents in the scaled system, $M^0=\sum\nolimits_{j\in[J]}M^0_j=M/h$, and $\bar{M}\in[h]$
is the maximal number of agents in a single area that have no place to be located at a time slot under the index policy. 
Based on the discussion in Appendix~\ref{app:lemma:MAI_asym_opt}, 
the difference between MAI and the index policy becomes negligible as $h$ approaches infinity, and $\bar{M} = o(h)$.

\subsection{Asymptotic Optimality}\label{subsec:asym_opt}


Consider an integer $M\in\mathbb{N}_+$, policies $\phi^*_1,\phi^*_2,\ldots,\phi^*_M\in\tilde{\Phi}^h$ such that satisfy \eqref{eqn:lemma:indexability} and \eqref{eqn:vartheta}  for all $(i,j,k)\in[I]\times[J]\times[h]$ with $\pmb{\gamma}$ and $\phi_{i,j}(\pmb{\gamma}_{i,j})$ replaced by $\pmb{\gamma}^*$ and $\phi^*_m$ ($m\in[M]$), respectively,
and a probability vector $\bm{\pi}^*\in[0,1]^M$ such that \eqref{eqn:prop:strong_duality:6} is satisfied. 
The existence of such $M$, $\phi^*_m$, and $\pi^*_m$ is guaranteed by Proposition~\ref{prop:strong_duality}.
Let $\pmb{\phi}^*\coloneqq (\phi^*_m: m\in[M])$.
Define a \emph{mixed policy} $\psi(M,\pmb{\phi}^*,\bm{\pi}^*)$ by mixing these policies $\phi^*_1,\phi^*_2,\ldots,\phi^*_M$. 
We say a mixed policy $\psi(M,\bm{\phi}^*,\bm{\pi}^*)$ is employed if, for each $(i,j)\in[I]\times[J]$ and $m\in[M-1]$, 
$\lfloor \pi^*_m h\rfloor$ of the sub-processes associated with $(i,j)$ are under the policy $\phi^*_m$, and the remaining $h - \sum_{m\in[M-1]}\lfloor \pi_m h \rfloor$ sub-processes follow the policy $\phi^*_M$. 
Let $k^-(m) \coloneqq \sum_{m'\in[m]} \lfloor \pi^*_{m'} h\rfloor$ for $m\in[M-1]$, $k^-(0)=0$ and $k^-(M) = h$.
Without loss of generality, we assume that processes $\{S^{\psi(M,\pmb{\phi}^*,\bm{\pi}^*),k}_{i,j}(t), t\in[T]\}$ with $k=k^-(m-1)+1, k^-(m-1)+2,\ldots, k^-(m)$ are under the policy $\phi^*_m$.
Such a mixed policy $\psi(M,\bm{\phi}^*,\bm{\pi}^*)$, abbreviated as $\psi^*$, satisfies \eqref{eqn:constraint:exclusive:h} but not necessarily complies with \eqref{eqn:constraint:neighbourhood:h:relax} - it is not necessarily applicable to the relaxed MAB-ML problem.
It serves as an intermediate step and will be used to establish a relationship between a policy applicable to the original problem (that is, a policy $\phi\in\Phi^h$ satisfying \eqref{eqn:constraint:exclusive:h}) and an optimal solution to the relaxed problem. 

For $\phi\in\tilde{\Phi}^h$ and $h\in\mathbb{N}_+$, define
\begin{equation}\label{eqn:define_Gamma}
\Gamma^{h,\phi}\coloneqq  \frac{1}{h}\sum_{\begin{subarray}~t\in[T],\\i\in[I],\\j\in[J],\\k\in[h]\end{subarray}}\mathbb{E}\Bigl[c_{i,j}\Bigl(S^{\phi,k}_{i,j}(t),\mathcal{e}^{\phi,k}_{i,j}\bigl(\bm{S}^{\phi,h}(t),t\bigr),t\Bigr)\Bigr],
\end{equation}
as the expected cumulative cost of process $\{\bm{S}^{\phi,h}(t):~t\in[T]\}$ under policy $\phi$ and given initial distribution of $\bm{S}^{\phi,h}(1)$. 
Let $\hat{\Phi}^h\subset \tilde{\Phi}^h$ represent the set of all policies in $\tilde{\Phi}^h$ satisfying Constraint~ \eqref{eqn:constraint:neighbourhood:h:relax} (applicable to the relaxed problem), and define
$\Gamma^{h,*}\coloneqq \min_{\phi\in\hat{\Phi}^h} \Gamma^{h,\phi}$ as the minimum of the relaxed MAB-ML problem described in \eqref{eqn:obj:h:relax} and \eqref{eqn:constraint:neighbourhood:h:relax}.

\begin{theorem}[Asymptotic Optimality]\label{theorem:asym_opt}
For a policy $\phi\in \Phi^h$ satisfying \eqref{eqn:constraint:exclusive:h} and \eqref{eqn:constraint:neighbourhood:h}, if
\begin{equation}\label{eqn:theorem:asym_opt:1}
   \lim_{h \rightarrow +\infty}\bigl\lvert \Gamma^{h,\phi} - \Gamma^{h,\psi^*}\bigr\rvert = 0,
\end{equation}
then
\begin{equation}\label{eqn:theorem:asym_opt:2}
\lim_{h \rightarrow +\infty}\bigl\lvert \Gamma^{h,\phi} - \Gamma^{h,*}\bigr\rvert = 0.
\end{equation}
\end{theorem}
The proof of Theorem~\ref{theorem:asym_opt} is provided in Appendix~\ref{app:theorem:asym_opt}.
A policy $\phi\in\Phi^h$ satisfying  \eqref{eqn:constraint:exclusive:h} and \eqref{eqn:constraint:neighbourhood:h} is in fact applicable to the original MAB-ML described in \eqref{eqn:obj:h}, \eqref{eqn:constraint:exclusive:h} and \eqref{eqn:constraint:neighbourhood:h}. 
For such a policy $\phi$ applicable to the original MAB-ML, since $\phi\in\hat{\Phi}^h$ and $ \Gamma^{h,*}$ is a lower bound of the minimum of the original MAB-ML, Theorem~\ref{theorem:asym_opt} indicates that $\phi$ approaches optimality of the original MAB-ML as $h \rightarrow +\infty$; that is, such a policy $\phi$ is \emph{asymptotically optimal} to MAB-ML.

\begin{theorem}\label{theorem:IND_asym_opt}
When the MAB-ML process is~\partialref{definition:indexability}{strictly indexable}, if 
\begin{enumerate}[label=(\roman*)]
    \item a policy $\phi\in\Phi^h$ satisfies \eqref{eqn:constraint:exclusive:h} and \eqref{eqn:constraint:neighbourhood:h}, and
    \item for any $t\in[T]$ and $\pmb{s}\in\mathscr{S}^h$, 
    \begin{equation}\label{eqn:theorem:IND_asym_opt:1}
    \Bigl\lVert\pmb{a}^{\phi,h}(\pmb{s},t) - \pmb{a}^{\text{IND},h}(\pmb{s},t)\Bigr\rVert = o(h),  
    \end{equation}
    where $\lVert\cdot \rVert$ is the Euclidean norm,
\end{enumerate}
then,
\begin{equation}\label{eqn:theorem:IND_asym_opt:2}
    \lim_{h \rightarrow +\infty}\bigl\lvert \Gamma^{h,\phi} - \Gamma^{h,\psi^*}\bigr\rvert = 0.
\end{equation}
\end{theorem}
The proof of Theorem~\ref{theorem:IND_asym_opt} is provided in Appendix~\ref{app:theorem:IND_asym_opt}.
Such a $\phi\in\Phi^h$ satisfying \eqref{eqn:constraint:exclusive:h}, \eqref{eqn:constraint:neighbourhood:h} and \eqref{eqn:theorem:IND_asym_opt:1} is applicable to MAB-ML and, from Theorems~\ref{theorem:IND_asym_opt} and~\ref{theorem:asym_opt}, achieves \eqref{eqn:theorem:asym_opt:2}.
Since $\Gamma^{h,*}$ is a lower bound to the minimum of MAB-ML, such a policy $\phi$ is asymptotically optimal.
\begin{lemma}\label{lemma:MAI_asym_opt}
When MAB-ML reduces to the patrol problem (a special case explained in Section~\ref{subsec:patrolling}), the MAI policy satisfies \eqref{eqn:constraint:exclusive:h}, \eqref{eqn:constraint:neighbourhood:h} and \eqref{eqn:theorem:IND_asym_opt:1}.
\end{lemma}
The proof of Lemma~\ref{lemma:MAI_asym_opt} is provided in Appendix~\ref{app:lemma:MAI_asym_opt}.
\begin{corollary}\label{coro:MAI_asym_opt}
When MAB-ML reduces to the patrol problem (a special case explained in Section~\ref{subsec:patrolling}), the MAI policy satisfies \eqref{eqn:theorem:IND_asym_opt:2}.
\end{corollary}
\plainproof Corollary~\ref{coro:MAI_asym_opt} is a direct result of Theorem~\ref{theorem:IND_asym_opt} and Lemma~\ref{lemma:MAI_asym_opt}.\endproof

Since MAI is applicable to the patrol problem (satisfying~\eqref{eqn:constraint:exclusive:h} and \eqref{eqn:constraint:neighbourhood:h}, under the strict indexability, MAI is asymptotically optimal to the patrol problem.

More importantly, we obtain the following bound on the performance deviation between a policy and an optimal solution of the original MAB-ML for a finite $h < +\infty$.
\begin{theorem}[Exponential Convergence]\label{theorem:exp_convergence}
If the MAB-ML process is strictly indexable and a policy $\phi\in\Phi^h$ satisfies \eqref{eqn:constraint:exclusive:h}, \eqref{eqn:constraint:neighbourhood:h} and \eqref{eqn:theorem:IND_asym_opt:1}, then, for any $\delta>0$, there exist $H>0$ and $s>0$ such that, for all $h > H$,
\begin{equation}\label{eqn:exp_convergence:1}
    \bigl\lvert \Gamma^{h,\phi} - \Gamma^{h,*}\bigr\rvert \leq e^{-sh} +\delta.
\end{equation}
\end{theorem}
The theorem is proved in Appendix~\ref{app:theorem:exp_convergence}.

{
\begin{figure*}[t]
\centering
\subfigure[]{\includegraphics[width=0.32\linewidth]{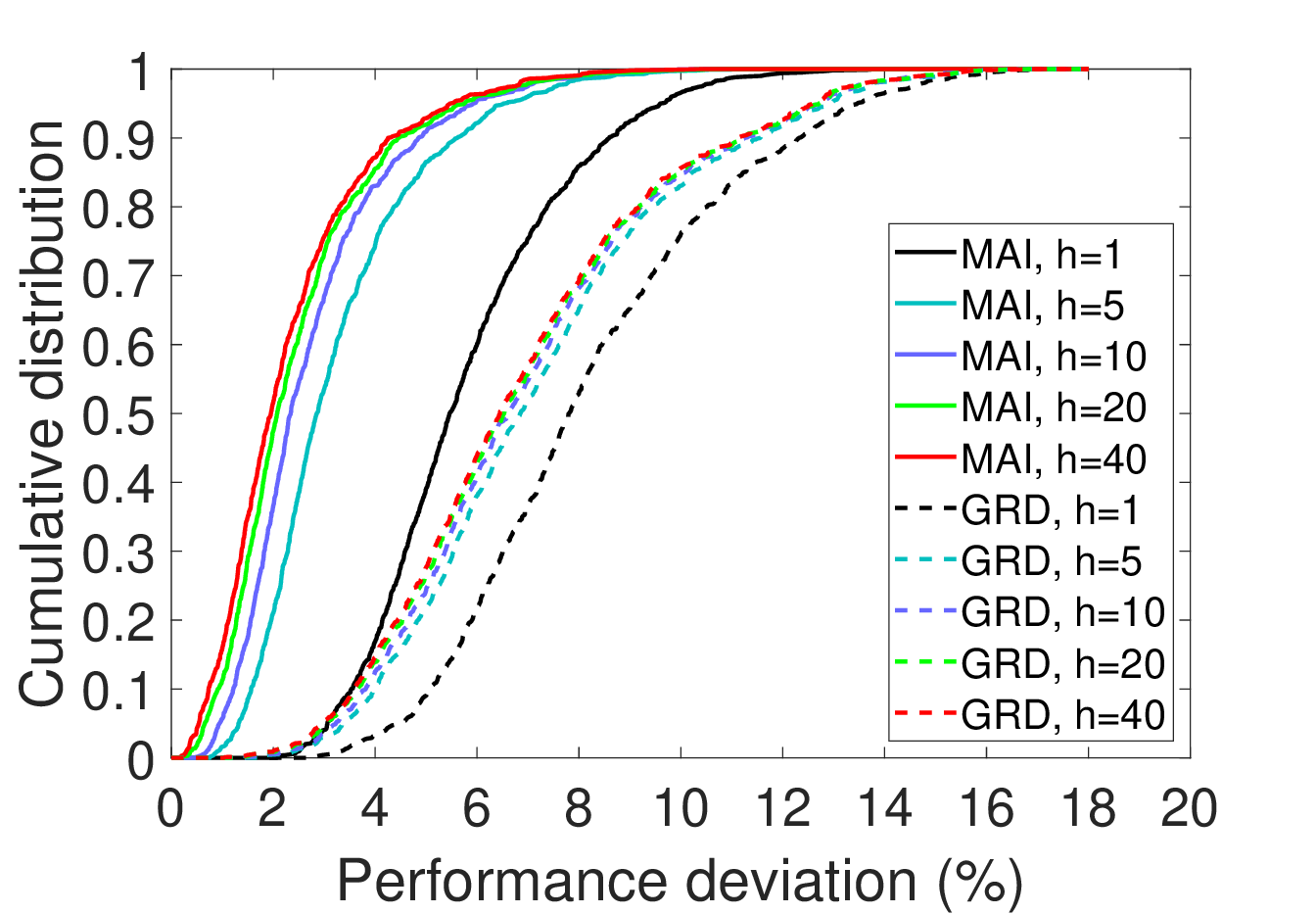}}
\subfigure[]{\includegraphics[width=0.32\linewidth]{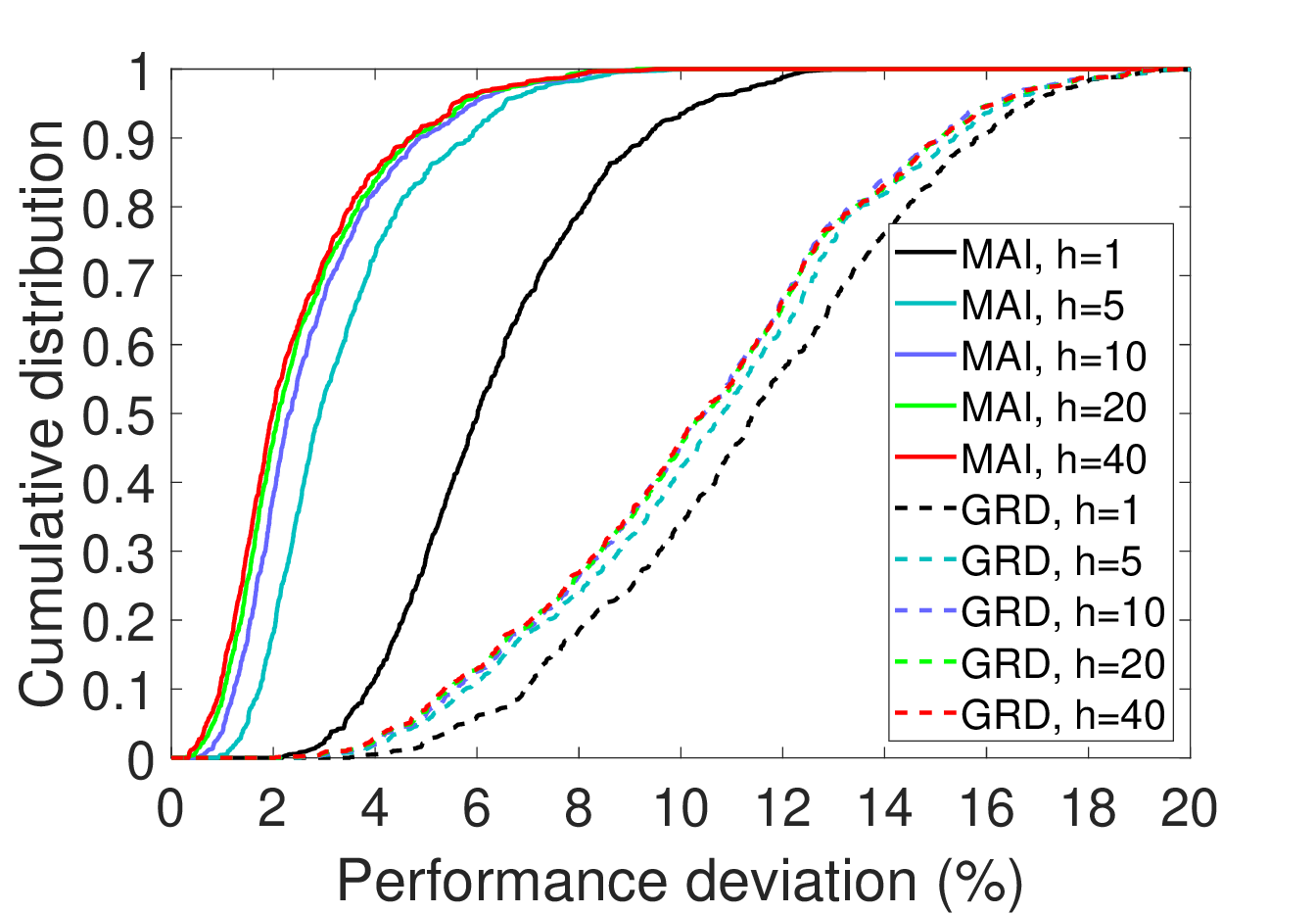}}
\subfigure[]{\includegraphics[width=0.32\linewidth]{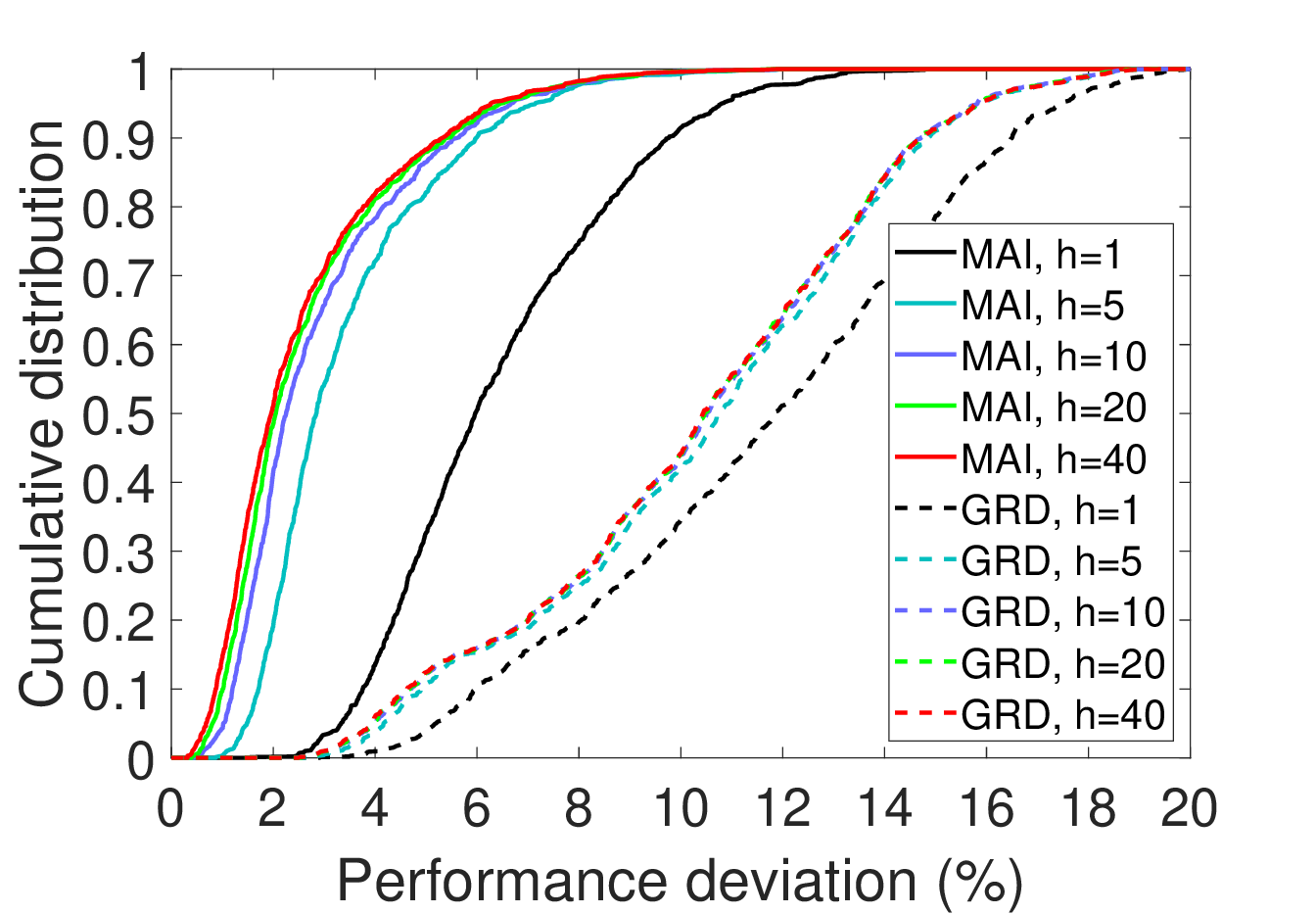}}
\caption{Cumulative distribution of performance deviations of MAI and GRD (Greedy) for (a) Case~I, (b) Case~II, and (c) Case~III. \label{fig:cdf}}
\end{figure*}
}

When MAB-ML reduces to the patrol problem (a special case explained in Section~\ref{subsec:patrolling}), 
since the MAI policy satisfies \eqref{eqn:constraint:exclusive:h}, \eqref{eqn:constraint:neighbourhood:h} and \eqref{eqn:theorem:IND_asym_opt:1} (Lemma~\ref{lemma:MAI_asym_opt}), 
for any $\delta>0$, there exists $H>0$ and $s>0$ such that, for all $h > H$, the MAI policy achieves \eqref{eqn:exp_convergence:1} with $\phi=\text{MAI}$.
Since the standard discrete-time finite-time horizon RMAB is also a special case of MAB-ML - Theorem~\ref{theorem:exp_convergence} is applicable to RMAB.
In Theorem~\ref{theorem:exp_convergence}, the deviation bound is tighter than the $O(\sqrt{N})$ bound for the standard discrete-time finite-time-horizon RMAB problem obtained in \cite{brown2020index}.

In \eqref{eqn:exp_convergence:1},  $\delta$ is a precision parameter that can be arbitrarily small.
Theorem~\ref{theorem:exp_convergence} ensures that the performance deviation between $\phi$ and an optimal solution of MAB-ML diminishes exponentially in the scaling parameter $h$. It implies that $\phi$ becomes sufficiently close to optimality even in a relatively small system and will be better when the problem size increases.

The strict indexability (Definition~\partialref{definition:indexability}{Indexability}) can be achieved by appropriately attaching a sufficiently small $\epsilon_{i,j}(s,a,t) > 0$ to the cost rate $c_{i,j}(s,a,t)$ ($(i,j)\in[I]\times [J]$, $s\in\mathscr{S}_{i,j}$, $a\in\{0,1\}$ and $t\in[T]$).
Compared to the absolute value of $\Gamma^{h,*}$, such $\epsilon_{i,j}(s,a,t) > 0$ ($(i,j)\in[I]\times [J]$, $s\in\mathscr{S}_{i,j}$, $a\in\{0,1\}$ and $t\in[T]$) are designed to be negligible, while they are crucial to determining the ranking of all the movements with achieved strict indexability.

Note that the asymptotically optimal policy is not necessarily unique. We have proved in Theorem~\ref{theorem:asym_opt} that any policies satisfying \eqref{eqn:theorem:asym_opt:1}, or those satisfies \eqref{eqn:theorem:IND_asym_opt:1} in Theorem~\ref{theorem:IND_asym_opt} for a strictly indexable MAB-ML process, are asymptotically optimal.

In the special case for multi-agent patrol, with strict indexability, the MAI policy is asymptotically optimal with exponentially fast convergence to optimality. 
Hence, while MAI can be implemented through various ways of adapting the movements, and may potentially lead to different action variables and different performance in the non-asymptotic regime, 
Theorem~\ref{theorem:exp_convergence} indicates that MAI will quickly converge to an optimal solution.



\section{Numerical Results}\label{sec:simulation}
For the patrol case, we numerically demonstrate the effectiveness of  MAI by comparing with a greedy approach and a lower bound of the minimum cost.
The $95\%$ confidence intervals of all the simulated results in this section, based on the Student t-distribution, are within $\pm3\%$ of the observed mean. 
The simulation programs are coded with g++ 11.4.0 and CPLEX 22.1.1.

As in Example~\ref{fig:model_example}, we consider the situation of crime detection with two ($J=2$) types of agents patrolling for two type of crimes: possession of weapon (type $j=1$) and vehicle crimes (type $j=2$).
We consider three patrol regions with $I=6$, $10$, and $14$ areas, respectively, which, as defined in Section~\ref{subsec:asym_regime}, are scaled by the scaling parameter $h \in \mathbb{N}_+$, and consider the minimization of the expected total crime rate of the patrol region with time horizon $T=10$. The topology of the three regions is provided in Appendix~\ref{app:simulation:topology}.
We refer to the three system settings with $I=6$, $10$, and $14$ as Cases I, II, and III, respectively.
At time $t$, the crime rate of type-$j$ in each sub-area $(i,k)\in[I]\times[h]$ follows a beta distribution, as in~\cite{ensign2018runaway}, with parameters $\mathcal{k}(S^{\phi,k}_{i,j}(t))=K^{\phi,k}_{i,j}(t)=(\alpha^{\phi,k}_{i,j}(t),\beta^{\phi,k}_{i,j}(t))$, where $\phi$ is a patrol policy, $\alpha^{\phi,k}_{i,j}(t),\beta^{\phi,k}_{i,j}(t)\in\mathbb{N}_+$, and we fix the sum $\alpha^{\phi,k}_{i,j}(t)+\beta^{\phi,k}_{i,j}(t) = 50$. 
In this context, with fine-tuned values of the initial parameters $(\alpha^{\phi,k}_{i,j}(1),\beta^{\phi,k}_{i,j}(1))$, such $\alpha^{\phi,k}_{i,j}(t)$ intuitively represents the expected number of crimes per five thousand population. 
Here, we consider the 2022 crime rates of $I$ counties or main cities of United Kingdom~\cite{UKCrime} as the initial crime rates of sub-areas associated with the $I$ areas. 
For each type $j\in[J]$, the number of agents is randomly generated with mean $M^0_j h$ where $M^0_j$ is also a random variable uniformly randomly generated from $\{1,2,\ldots, I/2\}$.
The details of the initial values of $(\alpha^{\phi,k}_{i,j}(1),\beta^{\phi,k}_{i,j}(1))$ for all $(i,j,k)\in[I]\times[J]\times[h]$, together with randomly generated agent numbers and initial positions of all the agents, are provided in Appendix~\ref{app:simulation:initialization:knowledge_position}.
In particular, for sub-area $(i,k)\in[I]\times[h]$ and agent type $j\in[J]$, recall the state variable $S^{\phi,k}_{i,j}(t)=\bigl(K^{\phi,k}_{i,j}(t),I^{\phi,k}_{i,j}\bigr)$, where $K^{\phi,k}_{i,j}(t)$ is the controller's knowledge of sub-area $(i,k)$ and $I^{\phi,k}_{i,j}(t)$ indicates whether sub-area $(i,k)$ is patrolled by a type-$j$ agent at time $t$.
In this section, let $K^{\phi,k}_{i,j}(t) = \bigl(\alpha^{\phi,k}_{i,j}(t), \beta^{\phi,k}_{i,j}(t)\bigr)$, and we rewrite the state variable $S^{\phi,k}_{i,j}(t) = \Bigl(\bigl(\alpha^{\phi,k}_{i,j}(t),\beta^{\phi,k}_{i,j}(t)\bigr),I^{\phi,k}_{i,j}(t)\Bigr)$ for the sub-process associated with sub-area $(i,k)$ and agent type $j$.

Define a random signal $D^{\phi,k}_{i,j}(t)$ that takes binary numbers drawn from a Bernoulli distribution with success probability $\alpha^{\phi,k}_{i,j}(t)/\bigl(\alpha^{\phi,k}_{i,j}(t)+\beta^{\phi,k}_{i,j}(t)\bigr)$, and define a function for normalization $\mathcal{N}(\alpha,\beta)\coloneqq \lfloor 50\alpha/(\alpha+\beta)\rfloor$.
For sub-area $(i,k)$ patrolled at time $t$ with $\bigl(\alpha^{\phi,k}_{i,j}(t),\beta^{\phi,k}_{i,j}(t)\bigr) =(\alpha,\beta)$, if $D^{\phi,k}_{i,j}(t) = 1$ (the patrol agent reports a crime and arrests the criminals), then $\alpha^{\phi,k}_{i,j}(t+1) = \mathcal{N}(\alpha,\beta+\Delta_{\beta})$; otherwise (no reported crime), $\alpha^{\phi,k}_{i,j}(t+1) = \mathcal{N}(\alpha+\Delta^1_{\alpha},\beta)$, where $\Delta_{\alpha}^1,\Delta_{\beta}\in\mathbb{N}_+$.
Similarly, for sub-area $(i,k)$ without any agent at time $t$, if $D^{\phi,k}_{i,j}(t) = 1$, then $\alpha^{\phi,k}_{i,j}(t+1) = \mathcal{N}(\alpha+\Delta^2_{\alpha},\beta)$ with $\Delta_{\alpha}^2\in\mathbb{N}_+$; otherwise, $\alpha^{\phi,k}_{i,j}(t+1)=\alpha^{\phi,k}_{i,j}(t)$. 
With determined $\alpha^{\phi,k}_{i,j}(t+1)$, $\beta^{\phi,k}_{i,j}(t+1) = 50 - \alpha^{\phi,k}_{i,j}(t+1)$.
The increment parameters $\Delta_{\alpha}^1$, $\Delta_{\alpha}^2$, and $\Delta_{\beta}$ are uniformly randomly generated from $\{2,3,\ldots,6\}$, $\{5,6,\ldots,9\}$, and $\{1,2,\ldots,5\}$, respectively.

In our simulations, we consider a lower bound of the minimum cost of the problem described in \eqref{eqn:obj:h}, \eqref{eqn:constraint:exclusive:h}, and \eqref{eqn:constraint:neighbourhood:h}, where the parameters are specified for the patrol and crime detection case.
More precisely, for any $(i,j)\in[I]\times[J]$, $\mathscr{B}_{i,j}=\bar{\mathscr{B}}_{i,j} = \bar{\mathscr{B}}^{-1}_{i,j}$ is specified as the set of area $i$ and its adjacent areas based on the topology provided in Appendix~\ref{app:simulation:topology}, $\mathcal{w}_{i,i',j}=1$, the cost function $c_{i,j}(S^{\phi,k}_{i,j}(t),e,t)$ is set to be the normalized crime rate $100\times \alpha^{\phi,k}_{i,j}(t)/\bigl(\alpha^{\phi,k}_{i,j}(t)+\beta^{\phi,k}_{i,j}(t)\bigr)$, and taking action $a^{\phi,k}_{i,i',j}(\cdot,\cdot)=1$ ($0$) means (not) moving a type-$j$ agent to sub-area $(i,k)$ from a sub-area of area $i'$. 
The lower bound is achieved by the maximum of the linear problem described in \eqref{eqn:dual_problem:2:1}-\eqref{eqn:dual_problem:2:3}.
Based on Corollary~\ref{coro:MAI_asym_opt}, the MAI policy approaches this lower bound and the minimum of the original patrol problem as $h\rightarrow +\infty$.
Based on Theorem~\ref{theorem:exp_convergence}, there exists $H>0$ such that, for all $h>H$, the performance deviation between MAI and this lower bound, as well as the minimum of the original patrol problem, diminishes exponentially as $h\rightarrow +\infty$.

Consider a baseline policy, which is a greedy approach that always prioritizes sub-areas with the highest crime rates and is denoted by $\phi=$GRD. That is, the greedy approach acts in the same way as MAI except that, instead of the indices, it ranks all the movements in the descending order of $\alpha^{h,\text{GRD}}_{i,j}(t)/\bigl(\alpha^{h,\text{GRD}}_{i,j}(t)+\beta^{h,\text{GRD}}_{i,j}(t)\bigr)$.

In Figure~\ref{fig:cdf}, for the above-mentioned three cases, we present the cumulative distributions of the performance deviations between specified policies, $\phi=$MAI and GRD, and the lower bound of optimality for $h=1,5,10,20,40$ in a thousand runs, each of which is with randomly generated simulation settings. 
Given a set of instances of the randomly generated system parameters, the performance deviation of policy $\phi$ is defined as $(\Gamma^{h,\phi}-\Gamma^*)/\Gamma^*$, where $\Gamma^{h,\phi}$ is approximated through the observed mean cost of a thousand simulation runs under $\phi$, and $\Gamma^*$ is the lower bound of optimality.
In all the tested cases, performance deviations of MAI drop quickly as $h$ increases, and, for Case I with $h=40$, they become less then $3\%$ in over $75\%$ of the tested simulations.
Recall that the confidence interval of the simulations is maintained within $\pm3\%$ of the observed mean. 
It implies that MAI has already been close to optimality with relatively small $h$, which is consistent with the theoretical conclusions in Corollary~\ref{coro:MAI_asym_opt} and Theorem~\ref{theorem:exp_convergence}.
Unlike MAI, in all the three cases in Figure~\ref{fig:cdf}, whatever $h$ is, the greedy approach continues to cause a significant drop in performance and demonstrates no trend of converging to optimality.

\begin{figure}[t]
\centering
\includegraphics[width=0.5\linewidth]{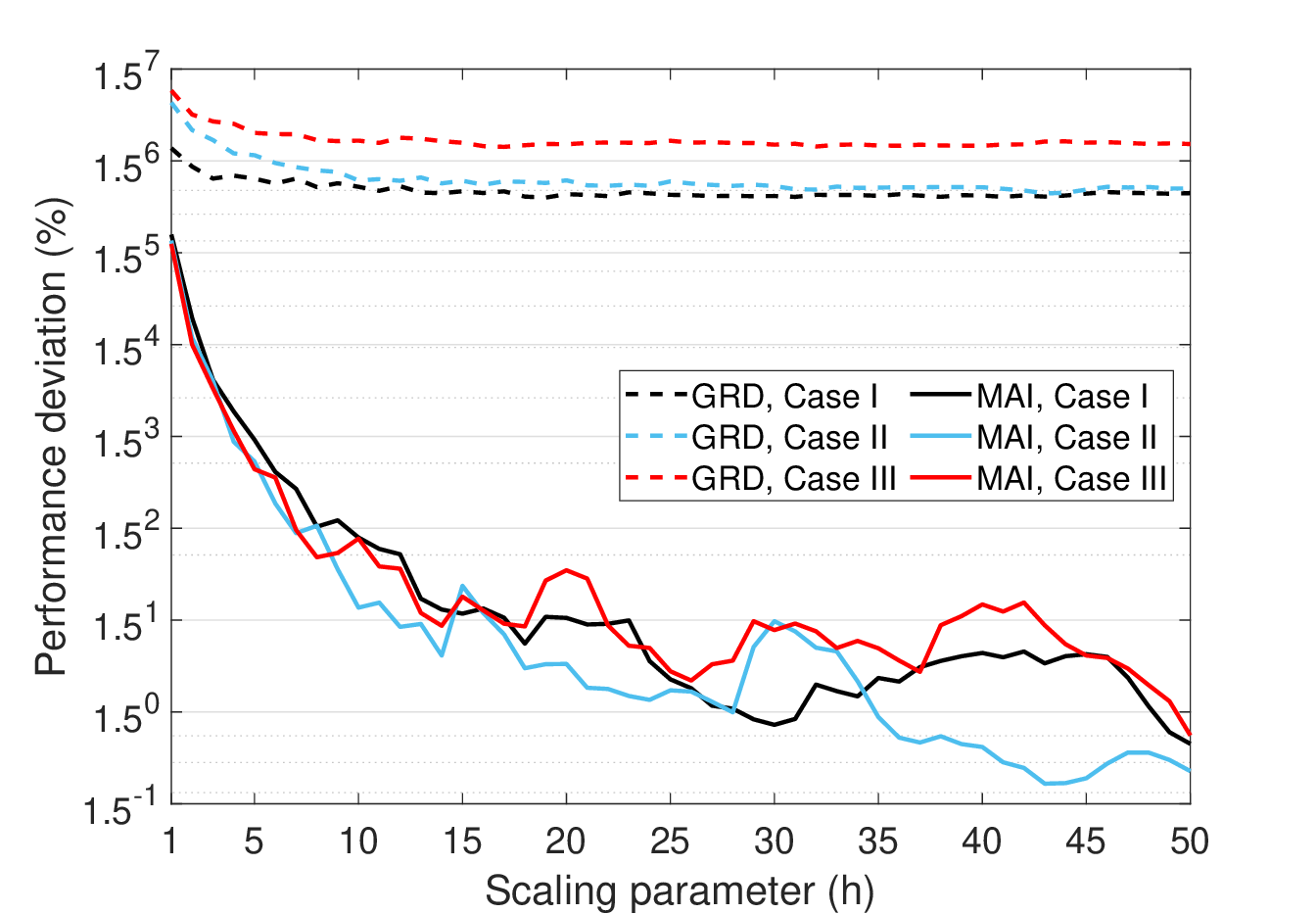}
\caption{Performance deviation of policy $\phi=$ MAI and GRD (Greedy) to the lower bound of optimality. \label{fig:case_study}}
\end{figure}

In Figure~\ref{fig:case_study}, for each of Cases I, II, and III, we consider an instance of the simulation system with random parameters and plot the performance deviation of specified policies to the lower bound of optimality $\Gamma^*$, given by $(\Gamma^{h,\phi}-\Gamma^*)/\Gamma^*$, against the scaling parameter $h$.
The detailed settings of the simulations in Figure~\ref{fig:case_study} are provided in Appendix~\ref{app:simulation:case_study}.
In Figure~\ref{fig:case_study}, we observe that the performance deviation of MAI diminishes quickly against $h$ and reduces to be less than $3\%$ for all $h\geq 30$. 
On the other hand, the performance deviation of the greedy approach is not sensitive to the scaling parameter, and, compared to the greedy approach, the advantages of MAI increase significantly as $h$ increases.  
It is consistent with Corollary~\ref{coro:MAI_asym_opt} and Theorem~\ref{theorem:exp_convergence} - MAI is asymptotically optimal and, for sufficiently large $h$ (the right tails of the curves in Figure~\ref{fig:case_study} for $h\geq 35$), the performance deviation of MAI diminishes quickly (exponentially) in $h$.

\section{Conclusions}\label{sec:conclusion}

We have proposed the MAB-ML problem that consists of multiple multi-action bandit processes, which are MDPs with finite state and action spaces.
These processes are coupled through multiple constraints that are linear to their action variables. 
MAB-ML is a generalized version of finite-time-horizon RMAB, which assumes a simple structure where each bandit process has only binary actions and is related to at most one constraint.
Recall that RMAB is already considered to be hard with intractable solutions~\cite{whittle1988restless,weber1990index,brown2020index} and is proved to be PSPACE-hard for the infinite time horizon case~\cite{papadimitriou1999complexity}.
Unlike RMAB, MAB-ML is capable of modeling the multi-agent patrol problem. 
We have modeled the patrol problem as a special case of MAB-ML that extends the assumptions on the past patrolling models, such as identical agents, submodular reward functions and capability of sensing any location at any time.
We have captured the relationship between the eligible movements, positions and profiles of the agents through  \eqref{eqn:constraint:neighbourhood}, leading to a stronger dependency among the action and state variables of the multi-action bandit processes than that of the conventional RMAB problem.
The patrolling problem is complicated by the high dimension of its state space that increases at least exponentially in the number of agents and the size of the patrol region.
The past RMAB methods cannot ensure the existence of scalable patrol policies with theoretically bounded performance degradation.
Nor can the past methods be used for analyzing the MAB-ML problem.

For MAB-ML, we have decomposed its state space, of which the size increases exponentially in the scaling parameter, into those small ones attached for the sub-problems by adapting the Whittle relaxation and Lagrangian DP techniques.
We have proposed the index policy for the general MAB-ML and proved Theorems~\ref{theorem:asym_opt}, \ref{theorem:IND_asym_opt}, and \ref{theorem:exp_convergence}  that theoretically bound the performance deviation between a range of  index-based policies (the policies $\phi\in\Phi^h$ satisfying \eqref{eqn:constraint:exclusive:h}, \eqref{eqn:constraint:neighbourhood:h} and \eqref{eqn:theorem:IND_asym_opt:1}) and optimality. 
These index-based policies are heuristics with reduced computational complexity - ensuring that they are scalable no matter how large MAB-M problem size, measured by $h$, is.
Theorems~\ref{theorem:asym_opt} and \ref{theorem:IND_asym_opt} guarantees that, when the system is strictly indexable, such index-based policies converges to optimality when $h\rightarrow +\infty$, and Theorem~\ref{theorem:exp_convergence} further shapes the convergence rate to be exponential in $h$.
In particular, for the patrol problem, a special case of MAB-ML, we have proposed the MAI policy with computational \color{black}complexity $O(MB\log(MB)+I^2JBh\bar{M})$, where $\bar{M}=o(h)$, for a system with scaling parameter $h$.
Recall that increasing $h$ indicates larger geographical patrol region or the same region with higher precision of collected information and more careful exploration.
Based on Corollary~\ref{coro:MAI_asym_opt}, we have proved that MAI is an index-based policy.
Together with Theorems~\ref{theorem:asym_opt} and \ref{theorem:exp_convergence}~we have proved that, when the areas are strictly indexable, the MAI policy approaches optimality as $h\rightarrow +\infty$ and the performance deviation between MAI and optimality diminishes exponentially as $h$ increases.
We have also provided numerical results for visualization of the theoretical discussions.

Based on the new knowledge developed in this paper, 
exploring specific application scenarios with practical configurations in more detail would be valuable for future research. 
For example, the knowledge variable $K^{\phi,h}_{i,j}(t)$ can be specified as parameters for Bayesian probabilities, Kalman filters, and other estimators in various scenarios. 
Also, the constraints can be adjusted to exclude collision for agents of different types, allow more than one agents of the same type to simultaneously patrol the same area, or include the time horizon $T$ as part of the control variables.
This will lead to more flexible and practical models and enhance the depth of practice implication of all the relevant studies.

\appendices

\section{Proof of Proposition~\ref{prop:decomposition}}
\label{app:prop:decomposition}

\proofRemark{Proposition~\ref{prop:decomposition}}
Consider the sum of minimum
\begin{equation}\label{eqn:dual_func:2_1}
\sum_{i\in[I]}\sum_{j\in[J]}\min_{\phi_{i,j}\in\tilde{\Phi}}   L^{\phi_{i,j}}_{i,j}(\pmb{\gamma}_{i,j}) \leq \min_{\phi\in\tilde{\Phi}} \sum_{i\in[I]}  \sum_{j\in[J]}  L^{\phi}_{i,j}(\pmb{\gamma}_{i,j}) = L(\pmb{\gamma}).
\end{equation}
The inequality comes from the concavity of minimization operation. For each of the minimum associated with $(i,j)$ at the left hand side of \eqref{eqn:dual_func:2_1}, 
\begin{equation}\label{eqn:dual_func:2_2}
\min_{\phi_{i,j}\in\tilde{\Phi}}   L^{\phi_{i,j}}_{i,j}(\pmb{\gamma}_{i,j}) = \sum_{\pmb{s}\in\mathscr{S}}\pi^0(\pmb{s})V^{\pmb{\gamma}_{i,j}}_{i,j}(\pmb{s},1),
\end{equation}
where $\pi^0(\pmb{s})$ is the given probability of the initial state $\bm{S}^{\phi}(1) = \pmb{s}$, and $V^{\pmb{\gamma}_{i,j}}_{i,j}(\pmb{s},t)$ ($i\in[I]$, $j\in[J]$, $t\in[T]$, $\pmb{s}\in\mathscr{S}$) is the minimal expected cumulative cost from time $t$ to $T$ given the state $\bm{S}^{\phi}(t)=\pmb{s}$ at time $t$ of process $\bigl\{S^{\phi}_{i,j}(t),t\in[T]\bigr\}$.
Recall that a policy $\phi_{i,j}\in\tilde{\Phi}$ is determined by its action variables $\bar{a}^{\phi_{i,j}}_{i,i',j}(\pmb{s},t)$ ($i,i'\in[I]$,$j\in[J]$, $t\in[T]$, $\pmb{s}\in\mathscr{S}$).
Consider the Bellman equation for process $\bigl\{S^{\phi}_{i,j}(t),t\in[T]\bigr\}$, $(i,j)\in[I]\times[J]$.
\begin{equation}\label{eqn:dual_func:2_3}
V^{\pmb{\gamma}_{i,j}}_{i,j}(\pmb{s},t) = \min_{\phi\in\tilde{\Phi}}\Bigl\{C^{\pmb{\gamma}_{i,j}}_{i,j}(s_{i,j},\bar{\pmb{a}}^{\phi}_{i,j}(\pmb{s},t),t)\\ + \sum_{\pmb{s}'\in\mathscr{S}}P_t(\pmb{s},\bar{\pmb{a}}^{\phi}(\pmb{s},t),\pmb{s}')V^{\pmb{\gamma}_{i,j}}_{i,j}(\pmb{s}',t+1)\Bigr\},
\end{equation}
with $V^{\pmb{\gamma}_{i,j}}_{i,j}(\pmb{s},T+1)\equiv 0$ for all $\pmb{s}\in\mathscr{S}$, where $\bar{\pmb{a}}^{\phi}(\pmb{s},t)\coloneqq (\bar{a}^{\phi}_{i,i',j}(\pmb{s},t):i,i'\in[I],j\in[J])$, and $P_t(\pmb{s},\pmb{a},\pmb{s}')$ is the transition probability from $\pmb{s}$ to $\pmb{s}'$ after action $\pmb{a}$ is taken.

For $i\in[I]$, $j\in[J]$, and given $s\in\mathscr{S}_{i,j}$, if $\bar{\pmb{a}}^{\phi}_{i,j}(\pmb{s},T)=\pmb{a}$ ($\pmb{a}\in\mathscr{A}_{i,j}$) minimizes the right hand side of \eqref{eqn:dual_func:2_3}, then there exists an optimal policy $\phi_{i,j}\in\tilde{\Phi}$ for $\min_{\phi\in\tilde{\Phi}}   L^{\phi}_{i,j}(\pmb{\gamma}_{i,j})$ such that, for any $\pmb{s}'\in\mathscr{S}$ with $s'_{i,j}=s_{i,j}$, $\bar{\pmb{a}}^{\phi_{i,j}}_{i,j}(\pmb{s}',T) = \pmb{a}$. That is, $\bar{\pmb{a}}^{\phi_{i,j}}_{i,j}(\pmb{s},T)$ and $V^{\pmb{\gamma}_{i,j}}_{i,j}(\pmb{s},T)$ dependent on $\pmb{s}$ through only $s_{i,j}$.
For any $t\in[T-1]$, assume that there exists an optimal policy $\phi_{i,j}\in\tilde{\Phi}$ for $\min\nolimits_{\phi\in\tilde{\Phi}}   L^{\phi}_{i,j}(\pmb{\gamma}_{i,j})$ where, for $\tau= t+1,t+2\ldots,T$ and any $\pmb{s}'\in\mathscr{S}$,  $\bar{\pmb{a}}^{\phi_{i,j}}_{i,j}(\pmb{s},\tau)$ and $V^{\pmb{\gamma}_{i,j}}_{i,j}(\pmb{s},\tau)$ are dependent on $\pmb{s}$ through only $s_{i,j}$. In this context, given $s_{i,j}\in\mathscr{S}_{i,j}$, if $\bar{\pmb{a}}^{\phi}_{i,j}(\pmb{s},t)=\pmb{a}$ ($\pmb{a}\in\mathscr{A}_{i,j}$) minimizes the right hand side of \eqref{eqn:dual_func:2_3} for some $\pmb{s}\in\mathscr{S}$, then, for any $\pmb{s}'\in\mathscr{S}$ with $s'_{i,j}=s_{i,j}$, there exists an optimal policy $\phi_{i,j}$ such that $\bar{\pmb{a}}^{\phi_{i,j}}_{i,j}(\pmb{s},t) = \pmb{a}$. 
Hence, for any $i\in[I]$, $j\in[J]$, $t\in[T]$ and $\pmb{s}\in\mathscr{S}$, there exists an optimal policy $\phi_{i,j}\in\tilde{\Phi}$ for  $\min\nolimits_{\phi\in\tilde{\Phi}}   L^{\phi}_{i,j}(\pmb{\gamma}_{i,j})$ such that $\bar{\pmb{a}}^{\phi_{i,j}}_{i,j}(\pmb{s},t)$ and $V^{\pmb{\gamma}_{i,j}}_{i,j}(\pmb{s},t)$ are dependent on $\pmb{s}$ through only $s_{i,j}$.
We rewrite $V^{\pmb{\gamma}_{i,j}}_{i,j}(\pmb{s},t)$ as $V^{\pmb{\gamma}_{i,j}}_{i,j}(s_{i,j},t)$ and obtain
\begin{equation}\label{eqn:dual_func:2_4}
V^{\pmb{\gamma}_{i,j}}_{i,j}(\pmb{s},t) = V^{\pmb{\gamma}_{i,j}}_{i,j}(s_{i,j},t) = \min_{\pmb{a}\in\mathscr{A}_{i,j}}\Bigl\{C^{\pmb{\gamma}_{i,j}}_{i,j}(s_{i,j},\pmb{a},t) \\ + \sum_{s'_{i,j}\in\mathscr{S}_{i,j}}P_{i,j,t}(s_{i,j},\mathcal{e}(\pmb{a}),s'_{i,j})V^{\pmb{\gamma}_{i,j}}_{i,j}(s'_{i,j},t+1)\Bigr\},
\end{equation}
where recall that, for $\pmb{a}=(a_{i'}:i'\in[I])\in\mathscr{A}_{i,j}$, $\mathcal{e}_{i,j}(\pmb{a})= 1-\prod_{i'\in[I]}(1-a_{i'})$ representing the probability of taking a non-zero action $(i,i',j)$ with $i'\in\mathscr{B}_{i,j}$,
and $P_{i,j,t}(s_{i,j},\mathcal{e}(\pmb{a}),s'_{i,j})$ is the transition probability from state $s_{i,j}$ to $s'_{i,j}$.
We can solve \eqref{eqn:dual_func:2_4} for all $(i,j)\in[I]\times[J]$ independently.
Construct a policy $\phi\in\tilde{\Phi}$ consisting of the optimal actions that minimize the right hand side of \eqref{eqn:dual_func:2_4} for all $(i,j)\in[I]\times[J]$, $t\in[T]$ and $s\in\mathscr{S}_{i,j}$.
Such a $\phi\in\tilde{\Phi}$ achieves the minimum of $\min_{\phi_{i,j}\in\tilde{\Phi}} L^{\phi_{i,j}}_{i,j}(\pmb{\gamma}_{i,j})$ for all $(i,j)\in[I]\times[J]$.
Together with \eqref{eqn:dual_func:2_1}, we prove the proposition.
\endproof

\section{Proof of Lemma~\ref{lemma:indexability}}
\label{app:lemma:indexability}
\proofRemark{Lemma~\ref{lemma:indexability}}
Recall that, based on the definition of $\mathscr{A}_{i,j}$ (in \eqref{eqn:define_A}), for any action $\pmb{a}\in\mathscr{A}_{i,j}$, $a_{i'} =0$ for all $i'\notin \mathcal{B}_{i,j}$.

From \eqref{eqn:dual_func:2_4}, there exists an optimal solution $\phi_{i,j}(\pmb{\gamma}_{i,j})\in\tilde{\Phi}$ such that  $\alpha^{\phi_{i,j}(\pmb{\gamma}_{i,j})}_{i,i',j}(s,t) = 0$ for all $i'\in \mathscr{B}_{i,j}$ if
\begin{multline}\label{eqn:lemma:indexability:2}
c_{i,j}(s,1,t)+\gamma_{i,j,t}g_{i,j}(s)- \max\nolimits_{i'\in\mathscr{B}_{i,j}}\theta_{i,i',j,t}(\gamma_{i',j,t}) 
+ \sum\nolimits_{s'\in\mathscr{S}_{i,j}}P_{i,j,t}(s,1,s') V^{\pmb{\gamma}_{i,j}}_{i,j}(s',t+1)\\
\geq c_{i,j}(s,0,t)+\gamma_{i,j,t}g_{i,j}(s)
+ \sum_{s'\in\mathscr{S}_{i,j}}P_{i,j,t}(s,0,s') V^{\pmb{\gamma}_{i,j}}_{i,j}(s',t+1).
\end{multline}
where $V^{\pmb{\gamma}_{i,j}}_{i,j}(s',t+1)$, for all $s'\in\mathscr{S}_{i,j}$, are given values and are independent from $\gamma_{i,j,t}$. 
Inequality \eqref{eqn:lemma:indexability:2} is equivalent to
\begin{equation}\label{eqn:lemma:indexability:3}
\max_{i'\in\mathscr{B}_{i,j}}\theta_{i,i',j,t}(\gamma_{i',j,t})\leq \vartheta^{\pmb{\gamma}_{i,j}}_{i,j}(s,t).
\end{equation}
Similarly,  if  
\begin{equation}\label{eqn:lemma:indexability:4}
\max_{i'\in\mathscr{B}_{i,j}}\theta_{i,i',j,t}(\gamma_{i',j,t})\geq \vartheta^{\pmb{\gamma}_{i,j}}_{i,j}(s,t),
\end{equation}
then there exists an optimal solution $\phi_{i,j}(\pmb{\gamma}_{i,j})\in\tilde{\Phi}$ such that 
$\alpha^{\phi_{i,j}(\pmb{\gamma}_{i,j})}_{i,i^*,j}(s,t) =1$ for an $i^*\in\arg\max_{i'\in\mathscr{B}_{i,j}}\theta_{i,i',j,t}(\gamma_{i',j,t})$. It proves the lemma.
\endproof

\section{Proof of Proposition~\ref{prop:strong_duality}}
\label{app:prop:strong_duality}

\begin{lemma}\label{lemma:continuity}
The dual function $L(\pmb{\gamma})$ is continuous and piece-wise linear in $\pmb{\gamma}\in\mathbb{R}^{IJT}$.
\end{lemma}
\proofRemark{Lemma~\ref{lemma:continuity}}
Define 
$L_{i,j}(\pmb{\gamma}_{i,j}) \coloneqq \min_{\phi\in\tilde{\Phi}}L^{\phi}_{i,j}(\pmb{\gamma}_{i,j})$.
From \eqref{eqn:dual_func:3}, 
\begin{equation}
    L(\pmb{\gamma}) 
= \sum_{i\in[I]}\sum_{j\in[J]} L_{i,j}(\pmb{\gamma}_{i,j}) 
= \sum_{i\in[I]}\sum_{j\in[J]} L_{i,j}^{\phi_{i,j}(\pmb{\gamma}_{i,j})}(\pmb{\gamma}_{i,j}),
\end{equation}
where $\phi_{i,j}(\pmb{\gamma}_{i,j})\in\tilde{\Phi}$ is a policy for sub-problem $(i,j)$ with given $\pmb{\gamma}_{i,j}$ that satisfy \eqref{eqn:lemma:indexability}.
From Lemma~\ref{lemma:indexability}, the policy $\phi_{i,j}(\pmb{\gamma}_{i,j})$ is optimal for sub-problem $(i,j)$.
For any $\pmb{\gamma}\in\mathbb{R}^{IJT}$, 
\begin{multline}\label{eqn:lemma:continuity:1}
L(\pmb{\gamma}) 
= \sum_{\begin{subarray}~t\in[T]\\i\in[I]\\j\in[J]\end{subarray}}\mathbb{E}\Bigl[c_{i,j}\bigl(S^{\phi_{i,j}(\pmb{\gamma}_{i,j})}_{i,j}(t),\pmb{\alpha}^{\phi_{i,j}(\pmb{\gamma}_{i,j})}_{i,j}(S^{\phi_{i,j}(\pmb{\gamma}_{i,j})}_{i,j}(t),t),t\bigr)\Bigr] 
- \sum_{\begin{subarray}~i\in[I]\\j\in[J]\\t\in[T]\end{subarray}}\mathbb{E}\Bigl[\sum_{i'\in\bar{\mathscr{B}}_{i,j}}\mathcal{w}_{i',i,j}\alpha^{\phi_{i',j}(\pmb{\gamma}_{i',j})}_{i',i,j}\bigl(S^{\phi_{i',j}(\pmb{\gamma}_{i',j})}_{i',j}(t),t\bigr)\\-g_{i,j}\bigl(S^{\phi_{i,j}(\pmb{\gamma}_{i,j})}_{i,j}(t)\bigr)\Bigr]\gamma_{i,j,t}
\end{multline}
From Lemma~\ref{lemma:indexability}, since the action variables $\alpha^{\phi_{i,j}(\pmb{\gamma}_{i,j})}_{i,i',j}(s,t)$ for $i,i'\in[I]$, $j\in[J]$, $t\in[T]$ and $s\in\mathscr{S}_{i,j}$ are piece-wise constant in $\pmb{\gamma}_{i,j}\in\mathbb{R}^{\lvert\{i\}\cup\bar{\mathscr{B}}^{-1}_{i,j}\rvert T}$, $L(\pmb{\gamma})$ is piece-wise linear in $\pmb{\gamma}\in\mathbb{R}^{IJT}$.

It remains to show the continuity of $L(\pmb{\gamma})$ at some special values of $\pmb{\gamma}\in\mathbb{R}^{IJT}$, for which the optimal action variables $\alpha^{\phi_{i,j}(\pmb{\gamma}_{i,j})}_{i,j}(s,t)$ for some $s\in\mathscr{S}_{i,j}$ and $t\in[T]$ flip from $0$ to $1$ or from $1$ to $0$.
We refer to these special values of $\pmb{\gamma}$ as the \emph{turning points}.
From Lemma~\ref{lemma:indexability}, at each turning point, there are some $i\in[I]$, $j\in[J]$, $t\in[T]$, $i'\in\mathscr{B}_{i,j}$ and $s\in\mathscr{S}_{i,j}$ such that 
\begin{equation}
\theta_{i,i',j,t}(\gamma_{i',j,t})=\max_{i''\in\mathscr{B}_{i,j}}\theta_{i,i'',j,t}(\gamma_{i'',j,t}) = \vartheta^{\pmb{\gamma}_{i,j}}_{i,j}(s,t).
\end{equation}
In this case, the value of $L(\pmb{\gamma})$ does not change for $\alpha^{\phi_{i,j}(\pmb{\gamma}_{i,j})}_{i,i',j}(s,t)$ taking any value in $[0,1]$. 
That is, at the turning point, $L(\pmb{\gamma})$ is continuous if any $\alpha^{\phi_{i,j}(\pmb{\gamma}_{i,j})}_{i,j}(s,t)$ flips its value.
This proves the lemma.
\endproof

\proofRemark{Proposition~\ref{prop:strong_duality}}
Recall that a policy $\phi\in\tilde{\Phi}$ is determined by its action vector $\pmb{\alpha}^{\phi}\coloneqq \Bigl(\alpha^{\phi}_{i,i',j}(s,t):i,i'\in[I],j\in[J],s\in\mathscr{S}_{i,j},t\in[T]\Bigr)$, which takes values in $\mathscr{A}\coloneqq \prod_{(i,j)\in[I]\times[J]} \Bigl(\mathscr{A}_{i,j}\Bigr)^T$.
Following similar ideas for proving \cite[Proposition 4]{brown2020index}, 
for $\pmb{\alpha}^{\phi}\in \mathscr{A}$ and $\pmb{\gamma}\in\mathbb{R}^{IJT}$, define
\begin{equation}\label{eqn:prop:strong_duality:1}
\bar{L}(\pmb{\gamma},\pmb{\alpha}^{\phi}) \coloneqq -\sum_{i\in[I]}\sum_{j\in[J]}L^{\phi}_{i,j}(\pmb{\gamma}_{i,j}),
\end{equation}
where the policy $\phi$ is determined by its action variables $\pmb{\alpha}^{\phi}$.
It is jointly continuous in $\pmb{\alpha}^{\phi}\in\mathscr{A}$ and $\pmb{\gamma}\in\mathbb{R}^{IJT}$ and, from \eqref{eqn:lemma:continuity:1}, is linear in $\pmb{\gamma}\in\mathbb{R}^{IJT}$. In this context, $\bar{L}(\pmb{\gamma}) = \max_{\pmb{\alpha}^{\phi}\in\mathscr{A}}\bar{L}(\pmb{\gamma},\pmb{\alpha}^{\phi}) = -L(\pmb{\gamma})$.
From Lemma~\ref{lemma:continuity}, the Lagrange dual function $L(\pmb{\gamma})$ is continuous and piece-wise linear in $\pmb{\gamma}\in\mathbb{R}^{IJT}$. Since $L(\pmb{\gamma})$ is also concave in $\pmb{\gamma}\in\mathbb{R}^{IJT}$, $\bar{L}(\pmb{\gamma}) = -L(\pmb{\gamma})$ is subdifferentiable.

For given  $\pmb{\gamma}\in\mathbb{R}^{IJT}$, define a set $\mathscr{A}(\pmb{\gamma})\subset \mathscr{A}$ of action vectors $\pmb{\alpha}^{\phi}$ of which the elements $\alpha^{\phi}_{i,i',j}(s,t)$ ($i,i'\in[I]$, $j\in[J]$, $t\in[T]$, $s\in\mathscr{S}_{i,j}$) satisfy \eqref{eqn:lemma:indexability} by substituting $\phi$ with $\phi_{i,j}(\pmb{\gamma}_{i,j})$. 
From Lemma~\ref{lemma:indexability}, for any $\pmb{\alpha}^{\phi}\in\mathscr{A}(\pmb{\gamma})$, $\bar{L}(\pmb{\gamma}) = \bar{L}(\pmb{\gamma},\pmb{\alpha}^{\phi})  = \max_{\pmb{\alpha}^{\phi'}\in\mathscr{A}} \bar{L}(\pmb{\gamma},\pmb{\alpha}^{\phi'})$. 
From the definition of $\bar{L}(\pmb{\gamma})$, for any $\pmb{\alpha}^{\phi}\in\mathscr{A}(
\pmb{\gamma})$,  $L(\pmb{\gamma})=-\bar{L}(\pmb{\gamma})= -\bar{L}(\pmb{\gamma},\pmb{\alpha}^{\phi})$.

Recall the definition of $L^{\phi}_{i,j}(\pmb{\gamma}_{i,j})$ in \eqref{eqn:define_func_i} and expand $\bar{L}(\pmb{\gamma},\pmb{\alpha}^{\phi})$ as follows.
\begin{multline}\label{eqn:prop:strong_duality:2}
\bar{L}(\pmb{\gamma},\pmb{\alpha}^{\phi}) =-\sum_{i\in[I]}\sum_{j\in[J]}L^{\phi}_{i,j}(\pmb{\gamma}_{i,j})\\
= -\sum_{t\in[T]}\sum_{i\in[I]}\sum_{j\in[J]}\mathbb{E}\Bigl[c_{i,j}\Bigl(S^{\phi}_{i,j}(t),\mathcal{e}\bigl(\pmb{\alpha}^{\phi}_{i,j}(S^{\phi}_{i,j}(t),t)\bigr),t\Bigr)\Bigr] 
+  \sum_{t\in[T]}\sum_{i\in[I]}\sum_{j\in[J]}\biggl(\sum_{i'\in\bar{\mathscr{B}}_{i,j}}\mathcal{w}_{i',i,j}\mathbb{E}\Bigl[\alpha^{\phi}_{i',i,j}\bigl(S^{\phi}_{i',j}(t),t\bigr)\Bigr]\\-\mathbb{E}\Bigl[g_{i,j}\bigl(S^{\phi}_{i,j}(t)\bigr)\Bigr]\biggr)\gamma_{i,j,t}.
\end{multline}
From \eqref{eqn:prop:strong_duality:2}, for any given $\pmb{\alpha}^{\phi}\in\mathscr{A}$, $\bar{L}(\pmb{\gamma},\pmb{\alpha}^{\phi})$ is differentiable in $\pmb{\gamma}\in\mathbb{R}^{IJT}$. 
The gradient $\nabla_{\pmb{\gamma}} \bar{L}(\pmb{\gamma},\pmb{\alpha}^{\phi})$ are continuous in $\pmb{\alpha}^{\phi}\in\mathscr{A}$.
Based on Danskin's Theorem \cite[Theorem 4.5.1]{bertsekas2003convex}, 
\begin{equation}\label{eqn:prop:strong_duality:3}
\partial \bar{L}(\pmb{\gamma}) = \bm{conv} \Bigl\{\nabla_{\pmb{\gamma}} \bar{L}(\pmb{\gamma},\pmb{\alpha}^{\phi}) | \pmb{\alpha}^{\phi}\in\mathscr{A}(\pmb{\gamma})\Bigr\},
\end{equation}
where $\partial \bar{L}(\pmb{\gamma})$ represents the subdifferential (the set of all subgradients) of $\bar{L}(\pmb{\gamma})$, and $\bm{conv} \mathscr{X}$ is the convex hull of the set $\mathscr{X}$.

From \cite[Theorem 4.7.2]{bertsekas2003convex}, $\pmb{\gamma}^*\in\mathbb{R}^{IJT}$ minimizes $\bar{L}(\pmb{\gamma})$, or, equivalently, maximizes $L(\pmb{\gamma})$, if and only if there exists a subgradient $\bm{g} \in \partial\bar{L}(\pmb{\gamma}^*)$ such that, for all $\pmb{\gamma}\in\mathbb{R}^{IJT}$, 
\begin{equation}\label{eqn:prop:strong_duality:4}
\bm{g}^T \Bigl(\pmb{\gamma}-\pmb{\gamma}^*\Bigr)\geq 0.
\end{equation}
That is, together with \eqref{eqn:prop:strong_duality:2} and \eqref{eqn:prop:strong_duality:3}, there exist $M\in\mathbb{N}_+$, $\pmb{\alpha}^{\phi^*_1},\pmb{\alpha}^{\phi^*_2},\ldots,\pmb{\alpha}^{\phi^*_M}\in\mathscr{A}(\pmb{\gamma}^*)$ and a probability vector $\bm{\pi}^*\in[0,1]^M$ such that $\bm{g}=\bm{g}^*\coloneqq \sum_{m\in[M]} \pi^*_m \nabla_{\pmb{\gamma}} \bar{L}(\pmb{\gamma}^*,\pmb{\alpha}^{\phi^*_m})$ satisfies \eqref{eqn:prop:strong_duality:4}.
It proves the proposition.
\endproof

\section{Proof of Theorem~\ref{theorem:asym_opt}}
\label{app:theorem:asym_opt}
\begin{lemma}\label{lemma:equivalence_2}
\begin{equation}\label{eqn:equivalence_2:1}
\lim_{h\rightarrow +\infty}\Bigl\lvert \Gamma^{h,\psi^*} - \Gamma^{h,*}\Bigr\rvert = 0,
\end{equation}
if and only if there exists $\phi\in\hat{\Phi}^h$ such that 
\begin{equation}\label{eqn:equivalence_2:2}
\lim_{h\rightarrow +\infty}\Bigl\lvert \Gamma^{h,\psi^*} - \Gamma^{h,\phi}\Bigr\rvert = 0
\end{equation}
\end{lemma}
\proofRemark{Lemma~\ref{lemma:equivalence_2}}
If \eqref{eqn:equivalence_2:1} holds, then there exists an policy $\phi^*(h)\in\hat{\Phi}^h$ with $\Gamma^{h,\phi^*(h)} = \Gamma^{h,*}$, which achieves \eqref{eqn:equivalence_2:2}. 
It remains to show the necessity of \eqref{eqn:equivalence_2:1} to the existence of $\phi\in\hat{\Phi}^h$ satisfying \eqref{eqn:equivalence_2:2}.

When there is a $\phi\in\hat{\Phi}^h$ with \eqref{eqn:equivalence_2:2},  we obtain
\begin{equation}\label{eqn:equivalence_2:3}
0 = \lim_{h\rightarrow +\infty} \Bigl(\Gamma^{h,\psi^*} - \Gamma^{h,\phi}\Bigr) \leq \lim_{h\rightarrow +\infty}\Bigl(\Gamma^{h,\psi^*} - \Gamma^{h,*}\Bigr).
\end{equation}
On the other hand, define the Lagrange dual function of the relaxed problem with scaling parameter $h\in\mathbb{N}_+$, 
\begin{equation}\label{eqn:equivalence_2:4}
L^{h,\phi}(\pmb{\gamma})\coloneqq \Gamma^{h,\phi} -\mathcal{g}^{h,\phi}(\pmb{\gamma}),
\end{equation}
where $\phi\in\tilde{\Phi}^h$, $\pmb{\gamma}\in\mathbb{R}_0^{IJT}$, 
and
\begin{multline}
\mathcal{g}^{h,\phi}(\pmb{\gamma})\coloneqq \sum_{i\in[I]}\sum_{j\in[J]}\sum_{t\in[T]}\frac{\gamma_{i,j,t}}{h}\biggl(\sum_{i'\in\bar{\mathscr{B}}_{i,j}}\sum_{k\in[h]}\mathcal{w}_{i',i,j}\mathbb{E}\Bigl[\alpha^{\phi,k}_{i',i,j}\bigl(S^{\phi,k}_{i',j}(t)\bigr)\Bigr]-\sum_{k\in[h]}\mathbb{E}\Bigl[g_{i,j}\bigl(S^{\phi,k}_{i,j}(t)\bigr)\Bigr]\biggr).
\end{multline}
Recall that $\psi^*$ consists of the policies $\phi^*_m$ ($m\in[M]$) that are optimal to $\min_{\phi\in\tilde{\Phi}^h}L^{h,\phi}(\pmb{\gamma}^*)$. It follows that
\begin{multline}\label{eqn:equivalence_2:5}
\lim_{h\rightarrow +\infty}\Bigl(\Gamma^{h,\psi^*} - \Gamma^{h,*}\Bigr)
= \lim_{h\rightarrow +\infty}\Bigl(L^{h,\psi^*}(\pmb{\gamma}^*) 
+ \mathcal{g}^{h,\psi^*}(\pmb{\gamma}^*) - \Gamma^{h,*}\Bigr)\\
\leq \lim_{h\rightarrow +\infty} \Bigl(L^{h,\phi^*(h)}(\pmb{\gamma}^*) + \mathcal{g}^{h,\psi^*}(\pmb{\gamma}^*) - \Gamma^{h,*}\Bigr) 
= \lim_{h\rightarrow +\infty} \Bigl(g^{h,\psi^*}(\pmb{\gamma}^*) -\mathcal{g}^{h,\phi^*(h)}(\pmb{\gamma}^*)\Bigr)= 0,
\end{multline}
where $\phi^*(h)\in\hat{\Phi}^h$ is an optimal policy for the relaxed problem, which satisfies $\Gamma^{h,\phi^*(h)} = \Gamma^{h,*}$.
The last equality  in \eqref{eqn:equivalence_2:5} holds because $\mathcal{g}^{h,\psi^*}(\pmb{\gamma}^*)=o(h)/h$ and $\mathcal{g}^{h,\phi^*(h)}(\pmb{\gamma}^*)=0$, where $o(h)$ represents a real number such that $\lim_{h\rightarrow \infty} o(h)/h = 0$. 
From \eqref{eqn:equivalence_2:3} and \eqref{eqn:equivalence_2:5}, \eqref{eqn:equivalence_2:1} holds. It proves the lemma.
\endproof

\proofRemark{Theorem~\ref{theorem:asym_opt}}
If \eqref{eqn:theorem:asym_opt:1} holds, then there exists $\phi\in\Phi^h\subset\hat{\Phi}^h$ such that \eqref{eqn:equivalence_2:2} is satisfied.
From Lemma~\ref{lemma:equivalence_2}, \eqref{eqn:equivalence_2:1} holds, leading to \eqref{eqn:theorem:asym_opt:2}.
\endproof

\section{Proof of Theorem~\ref{theorem:IND_asym_opt}}
\label{app:theorem:IND_asym_opt}

We order the process-state-time (PST) tuples $(i,j,s,t)$ ($(i,j)\in[I]\times[J]$, $s\in\mathscr{S}_{i,j}$, $t\in[T]$) according to the ascending order of their indices $\min_{i'\in\mathscr{B}_{i,j}}\eta_{i,i',j}(s,t)$ and the descending order of the time label $t$.
Here, we use the index of the \emph{best} action, $\min_{i'\in\mathscr{B}_{i,j}}\eta_{i,i',j}(s,t)$, as the \emph{index of PST} $(i,j,s,t)$.
That is, for any $t\in[T]$, all the PST tuples $(i_1,j_1,s_1,t+1)$ ($(i_1,j_1)\in[I]\times[J],s_1\in\mathscr{S}_{i_1,j_1}$) are listed in bulk preceding to those labeled by $(i_2,j_2,s_2,t)$ ($(i_2,j_2)\in[I]\times[J],s_2\in\mathscr{S}_{i_2,j_2}$). 
For a given $t\in[T]$, $(i_1,j_1,s_1,t)$ is ordered before $(i_2,j_2,s_2,t)$ if $\min_{i'_1\in\mathscr{B}_{i_1,j_1}}\eta_{i_1,i'_1,j_1}(s_1,t) < \min_{i'_2\in\mathscr{B}_{i_2,j_2}}\eta_{i_2,i'_2,j_2}(s_2,t)$.
The tie is broken in an arbitrary manner.
We alternatively refer to the $\iota$th PST tuple as the PST tuple $\iota$ or $(i_{\iota},j_{\iota},s_{\iota},t_{\iota})$. Let $N\coloneqq T\sum_{i\in[I]}\sum_{j\in[J]}|\mathscr{S}_{i,j}|$ represent the total number of PSTs in the system.

For $\iota\in[N]$ and $t\in[T]$, if $t_{\iota} = t$, then let $Z^{h,\phi}_{\iota}(t)$ represent the proportion of sub-processes associated with $(i_{\iota},j_{\iota})$ that are in state $s_{\iota}$ at time $t$ under policy $\phi$. 
More precisely, for $t\in[T]$ and $\pmb{s}=(s^k_{i,j}: (i,j)\in[I]\times[J],k\in[h])\in\mathscr{S}^h$, define $\pmb{z}(\pmb{s},t)\coloneqq (z_{\iota}(\pmb{s},t): \iota\in[N])$ satisfying that, if $\iota\in[N]$ with $t_{\iota}=t$, 
\begin{equation}\label{define:proportion_subprocesses}
    z_{\iota}(\pmb{s},t)\coloneqq \frac{1}{IJh}\biggl\lvert\Bigl\{k\in[h]~|~s^k_{i_{\iota},j_{\iota}}= s_{\iota}\Bigr\}\biggr\rvert;
\end{equation}
otherwise, $z_{\iota}(\pmb{s},t)=0$.
It follows that $\bm{Z}^{h,\phi}(t) = \pmb{z}(\bm{S}^{\phi,h}(t),t)$.


Let $\bar{\phi}\in\Phi^h$ be a policy satisfying \eqref{eqn:constraint:exclusive:h}, \eqref{eqn:constraint:neighbourhood:h}, and \eqref{eqn:theorem:IND_asym_opt:1}.

\begin{proposition}\label{prop:equivalence_1}
If the MAB-ML process is strictly indexable, then, for any $T>0$ and $t\in[T]$, 
\begin{equation}\label{eqn:equivalence_1}
\lim_{h\rightarrow+\infty} \mathbb{E}\Bigl[\bm{Z}^{h,\rm IND}(t)\Bigr] = \lim_{h\rightarrow+\infty} \mathbb{E}\Bigl[\bm{Z}^{h,\psi^*}(t)\Bigr] \\= \lim_{h\rightarrow +\infty}\mathbb{E}\Bigl[\bm{Z}^{h,\bar{\phi}}(t)\Bigr],
\end{equation}
and, for any $\delta >0$,
\begin{equation}\label{eqn:equivalence_1:0}
    \lim_{h\rightarrow +\infty}\mathbb{P}\biggl\{\Bigl\lVert \bm{Z}^{h,\bar{\phi}}(t) - \mathbb{E}\bigl[\bm{Z}^{h,\bar{\phi}}(t)\bigr]\Bigr\rVert  > \delta\biggr\} = 0.
\end{equation}
where $\bm{Z}^{h,\rm IND}(t) $ and $\bm{Z}^{h,\bar{\phi}}(t)$ are the random variables $\bm{Z}^{h,\phi}(t)$ under the index policy and $\bar{\phi}$ ($\phi=\rm IND$ and $\bar{\phi}$), respectively.
\end{proposition}
Theorem~\ref{theorem:IND_asym_opt} is a direct result of Theorem~\ref{theorem:asym_opt} and Proposition~\ref{prop:equivalence_1}. 
We postpone the proof of Proposition~\ref{prop:equivalence_1} to the end of this appendix, and, here, start with the following discussions that are inevitable for completing the proof.


Similar to the PST tuples, we rank the action-process-state-time (APST) tuples $(i,i',j,s,t)$ ($i,i'\in[I],j\in[J],s\in\mathscr{S}_{i,j},t\in[T]$) according to the ascending order the action indices $\eta_{i,i',j}(s,t)$ and the descending order of the time label $t$. 
For any $t\in[T]$, all the APST tuples $(i_1,i'_1,j_1,s_1,t+1)$ ($i_1,i'_1\in[I],j_1\in[J],s_1\in\mathscr{S}_{i_1,j_1}$) are listed in bulk preceding to those labeled by $(i_2,i'_2,j_2,s_2,t)$ ($i_2,i'_2\in[I],j_2\in[J],s_2\in\mathscr{S}_{i_2,j_2}$). 
For a given $t\in[T]$, $(i_1,i'_1,j_1,s_1,t)$ is ordered before $(i_2,i'_2,j_2,s_2,t)$ if $\eta_{i_1,i'_1,j_1}(s_1,t) < \eta_{i_2,i'_2,j_2}(s_2,t)$.
For the tie case with the same index value, if $(i_1,j_1,s_1,t)=(i_2,j_2,s_2,t)$ and $i'_1 < i'_2$, then APST $(i_1,i'_1,j_1,s_1,t)$ proceeds $(i_2,i'_2,j_2,s_2,t)$; 
otherwise, the tie is broken in an arbitrary manner.
We alternatively refer to the $\kappa$th APST tuple as the APST tuple $\kappa$ or $(i_{\kappa},i'_{\kappa},j_{\kappa},s_{\kappa},t_{\kappa})$. Let $K\coloneqq T\sum_{i,i'\in[I]}\sum_{j\in[J]}|\mathscr{S}_{i,j}|$ represent the total number of APSTs in the system.
In this context, for any $\kappa\in[K]$, there exists a unique $\iota\in[N]$ for which $(i_{\iota},j_{\iota},s_{\iota},t_{\iota}) = (i_{\kappa},j_{\kappa},s_{\kappa},t_{\kappa})$.
We use $\iota(\kappa)$ to represent such a $\iota$ for a given $\kappa\in[K]$.

For a policy $\phi\in\tilde{\Phi}^h$, $\kappa\in[K]$ and $\pmb{z} \in [0,1]^N$ with $z_{\iota}>0$, let 
$\zeta^{h,\phi}_{\kappa}(\pmb{z})$ represent the proportion of sub-processes associated with $(i_{\kappa},j_{\kappa})$ that are in state $s_{\kappa}$ and select action $(i_{\kappa},i'_{\kappa},j_{\kappa})$ under policy $\phi$ when $\bm{Z}^{h,\phi}(t_{\kappa})=\pmb{z}$.
That is, for $\kappa\in[K]$, $h\in\mathbb{N}_+$, and $\phi\in\tilde{\Phi}^h$, given $\bm{Z}^{h,\phi}(t_{\kappa})=\pmb{z}$,
\begin{equation}
    \zeta^{h,\phi}_{\kappa}(\pmb{z})\coloneqq \frac{1}{IJh}\biggl\lvert
    \Bigl\{k\in[h]~\Bigl|~S^{\phi,k}_{i_{\kappa},j_{\kappa}}(t_{\kappa})=s_{\kappa}, \\a^{\phi,k}_{i_{\kappa},i'_{\kappa},j_{\kappa}}(\bm{S}^{\phi,h}(t_{\kappa}),t_{\kappa})=1\Bigr\}\biggr\rvert.
\end{equation}

In particular, for the index policy described in \eqref{eqn:index_policy}, we obtain, for $\pmb{z}\in [0,1]^N$ and $\kappa\in[K]$  with $i'_{\kappa}\in\mathscr{B}_{i_{\kappa},j_{\kappa}}$,
\begin{multline}\label{eqn:zeta}
    \zeta^{h, \rm IND}_{\kappa}(\pmb{z})\coloneqq \max\Biggl\{0,
    \min\biggl\{
    \sum_{\begin{subarray}
    ~\iota\in[N]:\\
    i_{\iota} = i'_{\kappa},\\
    j_{\iota} = j_{\kappa},\\
    t_{\iota} = t_{\kappa}
    \end{subarray}}z_{\iota}g_{i'_{\kappa},j_{\kappa}}(s_{\kappa})
    - \sum_{\begin{subarray}
    ~\kappa'\in [K]:\\
    \kappa' < \kappa,\\
    i'_{\kappa'} = i'_{\kappa},\\
    j_{\kappa'} = j_{\kappa},
    \end{subarray}}\mathcal{w}_{i_{\kappa'},i'_{\kappa'},j_{\kappa'}}\zeta^{h,\rm IND}_{\kappa'}(\pmb{z}), 
    ~z_{\iota(\kappa)}
    - \sum_{\begin{subarray}
    ~\kappa'\in[K]:\\
    \kappa' < \kappa,\\
    i_{\kappa'} = i_{\kappa},\\
    j_{\kappa'} = j_{\kappa},\\
    t_{\kappa'} = t_{\kappa}
    \end{subarray}}\zeta^{h,\rm IND}_{\kappa'}(\pmb{z})
    \biggr\}\Biggr\},
\end{multline}
where, for $(i,j)\in[I]\times[J]$ and any $i'\notin\bar{\mathscr{B}}_{i,j}$, we define $\mathcal{w}_{i',i,j} = 0$.
For $\kappa\in[K]$ with $i'_{\kappa}\notin \mathscr{B}_{i_{\kappa},j_{\kappa}}$, $\zeta^{h,\rm IND}_{\kappa}(\pmb{z}) \equiv 0$ for all $\pmb{z}\in[0,1]^N$.
For any given $\pmb{z}\in [0,1]^N$, the value of $\zeta^{h,\rm IND}_{\kappa}(\pmb{z})$ can be computed iteratively from $\kappa=1$ to $K$.
For $\iota \in[N]$, the probability of taking an action $(i_{\iota}, i', j_{\iota})$ with $i'\in\mathscr{B}_{i_{\iota},j_{\iota}}$ under the index policy when $\bm{Z}^{h,\rm IND}(t_{\iota}) = \pmb{z}$ is 
\begin{equation}\label{eqn:u}
    u^{\rm IND}_{\iota}(\pmb{z}) = \frac{1}{z_{\iota}}\sum\nolimits_{\begin{subarray}
        ~\kappa \in [K]:\\
        i_{\kappa}=i_{\iota},\\
        j_{\kappa} = j_{\iota},\\
        s_{\kappa} = s_{\iota},\\
        t_{\kappa} = t_{\iota}
    \end{subarray}}\zeta^{h,\rm IND}_{\kappa}(\pmb{z}),
\end{equation}
if $z_{\iota} > 0$.
For $\iota\in[N]$ and $\pmb{z}\in[0,1]^N$ with $z_{\iota} = 0$, 
define $\pmb{z}_{\iota}(a)\coloneqq (z_1,z_2,\ldots,z_{\iota-1},a,z_{\iota+1},\ldots,z_N)$, and $u^{\rm IND}_{\iota}(\pmb{z}) = \lim_{a\downarrow 0} u^{\rm IND}_{\iota}(\pmb{z}_{\iota}(a))$.

For $\iota,\iota' \in [N]$, $t\in \mathbb{N}_+$, $e\in\{0,1\}$ and $k\in[h]$, define a random variable 
$\xi^{k,e}_{\iota,\iota',t}$ that takes values in $\{0,1\}$:
if $i_{\iota}=i_{\iota'}$, $j_{\iota}=j_{\iota'}$, $t_{\iota}=t$ and $t_{\iota'} =  t+1 $, then the probability 
\begin{equation}
    \mathbb{P}\{\xi^{k,e}_{\iota,\iota',t}=1\} = \mathbb{P}\biggl\{S^{\phi,k}_{i_{\iota},j_{\iota}}(t+1) = s_{\iota'}\\~\Bigl|~S^{\phi,k}_{i_{\iota},j_{\iota}}(t)=s_{\iota},\mathcal{e}\Bigl(\pmb{a}^{\phi,k}_{i_{\iota},j_{\iota}}\bigl(\bm{S}^{\phi,h}(t),t\bigr)\Bigr) = e\biggr\},
\end{equation}
and $\mathbb{P}\{\xi^{k,e}_{\iota,\iota',t}=0\} = 1- \mathbb{P}\{\xi^{k,e}_{\iota,\iota',t}=1\}$; otherwise, $\xi^{k,e}_{\iota,\iota',t}=0$ with probability $1$. 
For $t < 1 $, define $\xi^{k,e}_{\iota,\iota',t}=0$ ($\iota,\iota'\in[N]$, $k\in[h]$, $e\in\{0,1\}$) with probability $1$, and, for $t\in[1,+\infty)\backslash \mathbb{N}_+$, let $\xi^{k,e}_{\iota,\iota',t} = \xi^{k,e}_{\iota,\iota',\lfloor t\rfloor}$. 
In this context, the trajectory of $\xi^{k,e}_{\iota,\iota',t}$ is almost continuous in $t\geq 0$ with finitely many discontinuities of the first kind within every finite time interval.
Let $\bm{\xi}^h_t\coloneqq (\xi^{k,e}_{\iota,\iota',t} : k\in[h],e\in\{0,1\}, \iota,\iota'\in [N])$.

For $t\in[T+1]$, the value of $\int_t^{t+1} \xi^{k,e}_{\iota,\iota',\tau} d\tau$ is an integer, representing the potential number of transitions for sub-process $\{S^{\phi,k}_{i,j}(t),t\in[T]\}$ in state $s_{\iota}$ at time $t$ under action $\pmb{a}$ with $\mathcal{e}(\pmb{a})=e$ transitioning to $s_{\iota'}$ at time $t+1$.  
We use the word ``potential'' because sub-process $\{S^{\phi,k}_{i,j}(t),t\in[T]\}$ may not in state $s_{\iota}$ at time $t$, $t_{\iota}$ is not $t$, nor the action taken has $\mathcal{e}(\pmb{a})=e$, so the real number of transitions may be less than the integral.

For $\bm{x}\in\mathbb{R}_0^N$, $\bm{\xi}\in\mathbb{R}^{2hN^2}$ and $\iota,\iota'\in[N]$, if $(i_{\iota},j_{\iota}) \neq (i_{\iota'},j_{\iota'})$ or $t_{\iota}+1\neq t_{\iota'}$, then define $Q^h(\iota,\iota',\bm{x},\bm{\xi})\coloneqq 0$.
For $\bm{x}\in\mathbb{R}_0^N$, $\bm{\xi}\in\mathbb{R}^{2hN^2}$ and $\iota,\iota'\in[N]$ with $(i_{\iota},j_{\iota}) = (i_{\iota'},j_{\iota'})$ and $t_{\iota}+1=t_{\iota'}$, define
\begin{equation}\label{eqn:Q}
Q^h(\iota,\iota',\bm{x},\bm{\xi}) \coloneqq 
\sum_{\highlight{k=h\lceil x^-_{\iota-1}/h\rceil+1}}^{\highlight{\mathcal{L}^h_{\iota}(\bm{x})}}\xi^{k,1}_{\iota,\iota'} \\+ 
\sum_{\highlight{k=\mathcal{L}^h_{\iota}(\bm{x})+1}}^{\highlight{h\lceil x^-_{\iota}/h\rceil}}\xi^{k,0}_{\iota,\iota'} +f^{h,\mathcal{a}}_{\iota,\iota'}(\bm{x},\bm{\xi}),
\end{equation}
where $x^-_{\iota} = \sum_{\iota'=1,i_{\iota'}=i_{\iota},j_{\iota'}=j_{\iota}}^{\iota} x_{\iota'}$ with $x^-_0 = 0$, 
\[\highlight{\mathcal{L}^h_{\iota}(\bm{x})\coloneqq h\lceil x^-_{\iota-1}/h\rceil+h\Bigl\lceil(\lceil x^-_{\iota}/h\rceil-\lceil x^-_{\iota-1}/h\rceil)u^{\text{IND}}_{\iota}(\frac{\bm{x}}{hIJ})\Bigr\rceil,}\] 
and $f^{h,\mathcal{a}}_{\iota,\iota'}(\bm{x},\bm{\xi})$ are appropriate functions to make $Q^h(\iota,\iota',\bm{x},\bm{\xi})$ Lipschitz continuous in $\bm{x}$ for any given $\bm{\xi}$ and $a\in(0,1)$. Similar to \cite{fu2020resource,fu2020energy}, such $f^{h,\mathcal{a}}_{\iota,\iota'}$ can be obtained by utilizing the Dirac delta function and we choose those satisfying $\lim_{a\downarrow 0} d~f^{h,\mathcal{a}}_{\iota,\iota'}/d \mathcal{a}  = 0$.


For $h=1$, $\mathcal{a}\in(0,1)$ and $\sigma >0$, define a variable $\bm{X}^{\sigma}_t$ on $t>0$, satisfying
\begin{equation}\label{eqn:x_derivative:1}
\dot{\bm{X}}^{\sigma}_t \coloneqq b(\bm{X}^{\sigma}_t,\bm{\xi}^1_{t/\sigma})
\end{equation}
where $\bm{\xi}^1_{t/\sigma}$ represents $\bm{\xi}^h_{t/\sigma}$ with $h=1$, and $b(\bm{X}^{\sigma}_t,\bm{\xi}^1_{t/\sigma})$ is a vector taking values in $\mathbb{R}^N$, for which the $\iota$th element is such that, for $\bm{x}\in\mathbb{R}_0^N$ and $\bm{\xi}\in\mathbb{R}^{^{2N^2}}$,
\begin{equation}\label{eqn:x_derivative:2}
b_{\iota}(\bm{x},\bm{\xi}) \coloneqq \sum_{\iota'\in[N]} \Bigl(Q^1(\iota',\iota,\bm{x},\bm{\xi}) - Q^1(\iota,\iota',\bm{x},\bm{\xi})\Bigr),
\end{equation}
with $Q^1$ the special case of $Q^h$ for $h=1$. Following the Lipschitz continuity of $Q^h$, $b(\bm{x},\bm{\xi})$ is Lipschitz continuous in $\bm{x}\in\mathbb{R}_0^N$ and $\bm{\xi}\in\mathbb{R}^{^{2hN^2}}$. From the definitions \eqref{eqn:Q} and \eqref{eqn:x_derivative:2}, for $\bm{x}\in\mathbb{R}_0^N$ and $\bm{\xi}\in\mathbb{R}^{^{2hN^2}}$, there exists a matrix $\tilde{\mathcal{Q}}(\bm{x})$ such that $b(\bm{x},\bm{\xi}) = \tilde{\mathcal{Q}}(\bm{x})\bm{\xi}$. We obtain that, for any $\bm{x}\in\mathbb{R}_0^N$, $h=1$, $\delta > 0$, $t,T>0$ and any function $\bar{b}(\bm{x})\in\mathbb{R}^N$, 
\begin{multline}\label{eqn:converge_b:1}
\mathbb{P}\biggl\{\Bigl\lVert\frac{1}{T} \int_t^{t+T}b(\bm{x},\bm{\xi}^1_{\tau})d\tau - \bar{b}(\bm{x})\Bigr\rVert > \delta\biggr\}\\
\leq \mathbb{P}\biggl\{\Bigl\lVert\tilde{\mathcal{Q}}(\bm{x})\frac{1}{T} \Bigl\lfloor\int_t^{t+T}\bm{\xi}^1_{\tau}d\tau \Bigr\rfloor - \bar{b}(\bm{x})\Bigr\rVert 
+\Bigl\lVert\tilde{\mathcal{Q}}(\bm{x})\frac{1}{T} \bigl(\int_t^{t+T}\bm{\xi}^1_{\tau}d\tau  - \Bigl\lfloor\int_t^{t+T}\bm{\xi}^1_{\tau} d\tau \Bigr\rfloor\bigr)\Bigr\rVert 
> \delta\biggr\}\\
\leq \mathbb{P}\biggl\{\Bigl\lVert\tilde{\mathcal{Q}}(\bm{x})\frac{1}{T} \Bigl\lfloor\int_t^{t+T}\bm{\xi}^1_{\tau}d\tau \Bigr\rfloor - \bar{b}(\bm{x})\Bigr\rVert 
+\frac{o(T)}{T}
> \delta\biggr\},
\end{multline}
where $\pmb{\xi}^1_t$ is a special case of $\pmb{\xi}^h_t$ with $h=1$.
Recall that, for $t\in\mathbb{N}_+$, $\Bigl\lfloor\int_t^{t+T}\bm{\xi}^1_{\tau}d\tau \Bigr\rfloor = \Bigl(\Bigl\lfloor\int_t^{t+T}\xi^{1,e}_{\iota,\iota',\tau}d\tau\Bigr\rfloor: \iota,\iota'\in[N],e\in\{0,1\}\Bigr)$, where each element $\Bigl\lfloor\int_t^{t+T}\xi^{1,e}_{\iota,\iota',\tau}d\tau\Bigr\rfloor$ is Poisson distributed with a rate, $\lambda_{\iota,\iota'}^e$, representing the expected number of potential transitions for process $\{S^{\phi,k}_{i_{\iota},j_{\iota}}(t),t\in[T]\}$ 
transitioning from state $s_{\iota}$ to $s_{\iota'}$ at time $t_{\iota}$. Let $\bm{\lambda}\coloneqq (\lambda_{\iota,\iota}^e: \iota,\iota'\in[N],e\in\{0,1\})$.
For $\bm{x}\in\mathbb{R}_0^N$, $h=1$, $\delta > 0$, and $t>0$, there exists $\bar{b}(\bm{x})=\mathbb{E}[b(\bm{x},\bm{\xi}^1_t)]$ such that
\begin{equation}\label{eqn:converge_b:2}
\tilde{\mathcal{Q}}(\bm{x})\frac{1}{T} \Bigl\lfloor\int_t^{t+T}\bm{\xi}^1_{\tau}d\tau \Bigr\rfloor - \bar{b}(\bm{x}) = \tilde{\mathcal{Q}}(\bm{x})
\biggl(\frac{1}{T} \Bigl\lfloor\int_t^{t+T}\bm{\xi}^1_{\tau}d\tau \Bigr\rfloor - \mathbb{E}\bigl[\bm{\xi}^1_t\bigr]\biggr)\\
=\tilde{\mathcal{Q}}(\bm{x})\biggl(\frac{1}{T} \Bigl\lfloor\int_t^{t+T}\bm{\xi}^1_{\tau}d\tau \Bigr\rfloor - \bm{\lambda}\biggr).
\end{equation}
Hence, by the Law of Large Numbers, for any $\bm{x}\in\mathbb{R}_0^N$, $h=1$, $\delta > 0$, and $t>0$, there exists $\bar{b}(\bm{x})=\mathbb{E}[b(\bm{x},\bm{\xi}^1_t)]$ such that
\begin{equation}\label{eqn:converge_b:3}
\lim_{T\rightarrow \infty}\sup_{t > 0}\biggl\lvert\mathbb{P}\Bigl\{\Bigl\lVert\tilde{\mathcal{Q}}(\bm{x})\Bigl(\frac{1}{T} \Bigl\lfloor\int_t^{t+T}\bm{\xi}^1_{\tau}d\tau \Bigr\rfloor - \bm{\lambda}\Bigr)\Bigr\rVert+\frac{o(T)}{T} > \delta\Bigr\}\biggr\rvert = 0.
\end{equation}
From \eqref{eqn:converge_b:1} and \eqref{eqn:converge_b:3}, for any $\bm{x}\in\mathbb{R}_0^N$, $h=1$ and $\delta > 0$, there exists $\bar{b}(\bm{x})=\mathbb{E}[b(\bm{x},\bm{\xi}^1_t)]$ such that
\begin{equation}\label{eqn:converge_b:4}
\lim_{T\rightarrow \infty}\mathbb{P}\biggl\{\Bigl\lVert\frac{1}{T} \int_t^{t+T}b(\bm{x},\bm{\xi}^1_{\tau})d\tau - \bar{b}(\bm{x})\Bigr\rVert > \delta\biggr\} = 0,
\end{equation}
uniformly in $t>0$. We will take $\bar{b}(\bm{x})=\mathbb{E}[b(\bm{x},\bm{\xi}^1_t)]$ for discussions in the sequel.
Let $\bar{\bm{x}}_t$ be the solution of $\dot{\bar{\bm{x}}}_t = \bar{b}(\bar{\bm{x}}_t)$ and $\bar{\bm{x}}_0 = \bm{X}^{\sigma}_0$. By Picard’s Existence Theorem \cite{coddington1955theory}, the solution $\bar{\bm{x}}_t$ uniquely exists.

Based on \cite[Chapter 7, Theorem 2.1]{freidlin2012random}, if there exists $\bar{\bm{x}}_t$ satisfying \eqref{eqn:converge_b:4} and $\lVert b(\bm{x},\bm{\xi}^1_t)\rVert^2 <+\infty$ for all $\bm{x}\in\mathbb{R}^N$, then, for any $T>0$ and $\delta>0$,
\begin{equation}\label{eqn:converge_x:1}
\lim_{\sigma \rightarrow 0}\mathbb{P} \Bigl\{\sup\nolimits_{1\leq t\leq T+1}\bigl\lVert \bm{X}^{\sigma}_t - \bar{\bm{x}}_t\bigr\rVert > \delta\Bigr\} = 0.
\end{equation}
We can translate the scaling effects in \eqref{eqn:converge_x:1} in another way. For $\bm{x}\in\mathbb{R}^N$, and $\bm{\xi}^h\in\mathbb{R}^{2hN^2}$, define a vector $b^h(\bm{x},\bm{\xi}^h) \in \mathbb{R}^N$, of which the $\iota$th element is
\begin{equation}
b^h_{\iota}(\bm{x},\bm{\xi}^h) \coloneqq \sum_{\iota'\in[N]} \left(Q^h(\iota',\iota,\bm{x},\bm{\xi}^h) - Q^h(\iota,\iota',\bm{x},\bm{\xi}^h)\right).
\end{equation}
Let $\sigma = 1/h$, where $h$ is the scaling parameter defined in Section~\ref{subsec:asym_regime}, and we can obtain, for any $\bm{x}\in\mathbb{R}^N$, $\sigma > 0 $, and $T>0$,
\begin{equation}\label{eqn:equivalence_b:1}
\int_0^T b(\bm{x},\bm{\xi}^1_{t/\sigma}) dt = \sigma \int_0^{T/\sigma}b(\bm{x},\bm{\xi}^1_t) dt = \frac{1}{h}\int_0^{hT}b(\bm{x},\bm{\xi}^1_t) dt.
\end{equation}
Since $\bm{\xi}^1_t$ is identically distributed for any $t>0$, for $T\in\mathbb{N}_+$, $\frac{1}{h}\int_0^{hT}b(\bm{x},\bm{\xi}^1_t) dt$ is in distribution equivalent to 
$\frac{1}{h} \int_0^T \sum_{n=1}^h \tilde{b}_n (\bm{x},\bm{\xi}^1_t) dt$,
where $\tilde{b}_n (\bm{x},\bm{\xi}^1_t)$ ($n\in[h]$) are independently identically distributed as $b(\bm{x},\bm{\xi}^1_t)$.
On the other hand, we observe that, for any $t>0$, $\sum_{n=1}^h \tilde{b}_n(\bm{x},\bm{\xi}^1_t)$ is in distribution equivalent to $b^h(h\bm{x},\bm{\xi}^h_t)$.
Together with \eqref{eqn:equivalence_b:1}, we obtain that $\int_0^T b(\bm{x},\bm{\xi}^1_{t/\sigma}) dt$ and $\frac{1}{h}\int_0^T b^h(h\bm{x},\bm{\xi}^h_t) dt$ are identically distributed.

Define $\bm{Z}^{\sigma}_t$ ($t>0$) as the solution of $\dot{\bm{Z}}^{\sigma}_t = \frac{1}{IJ}b(IJ\bm{Z}^{\sigma}_t,\bm{\xi}^1_{t/\sigma})$ and $\dot{\bm{Z}}^{h}_t = \frac{1}{hIJ}b^h(hIJ\bm{Z}^{\sigma}_t,\bm{\xi}^h_t)$, respectively, with $\bm{Z}^{\sigma}_1 = \bm{Z}^{h}_1 = \frac{\bm{x}_1}{IJ}$. Based on \eqref{eqn:converge_x:1}, for any \highlight{$T\in\mathbb{N}_+$} and $\delta>0$, 
\begin{equation}\label{eqn:converge_x:2}
\lim_{h \rightarrow +\infty}\mathbb{P} \Bigl\{\sup\nolimits_{1\leq t\leq T+1}\bigl\lVert \bm{Z}^h_t - \frac{\bar{\bm{x}}_t}{IJ}\bigr\rVert > \delta\Bigr\} = 0.
\end{equation}
That is, scaling $\sigma \rightarrow 0$ and $h\rightarrow +\infty$ leads to the same effect.
\highlight{For $T\in\mathbb{N}_+$ and $t\in[T]$,} the elements of the random vector $\bm{Z}^h_t$ in fact is \highlight{in distribution} equivalent to $\bm{Z}^{h,\rm IND}(t)$, representing the proportions of sub-processes in all the $N$ PST tuples at time $t$ given scaling parameter $h$, when the index policy characterized by $\bm{u}^{\text{IND}}(\bm{Z}^{h,\rm IND}(t))$ is employed.


\begin{lemma}\label{lemma:converge_z_opt}
For any $t\in[T]$ and $\delta>0$, 
\begin{equation}\label{eqn:converge_z_opt:1}
\lim_{h\rightarrow +\infty} \mathbb{P}\biggl\{\Bigl\lVert \bm{Z}^{h,\psi^*}(t) - \mathbb{E}\bigl[\bm{Z}^{h,\psi^*}(t)\bigr] \Bigr\rVert > \delta\biggr\} = 0.
\end{equation}
\end{lemma}
\proofRemark{Lemma~\ref{lemma:converge_z_opt}}
From the definitions, for $\iota\in[N]$ and $t\neq t_{\iota}$, $Z^{h,\psi^*}_{\iota}(t) \equiv 0$, and \eqref{eqn:converge_z_opt:1} straightforwardly holds.


For $i\in[I]$, $j\in[J]$, $k\in[h]$ and $t\in[T]$, define $\pmb{\Theta}^{\psi^*,k}_{i,j}(t) \coloneqq \Bigl(\Theta^{\psi^*,k}_{s,i,j}(t): s\in\mathscr{S}_{i,j} \Bigr)$, where, for $s\in\mathscr{S}_{i,j}$, if the state variable  $S^{\psi^*,k}_{i,j}(t)=s$ then $\Theta^{\psi^*,k}_{s,i,j}(t) = 1$; otherwise, $\Theta^{\psi^*,k}_{s,i,j}(t) = 0$.
For given $i\in[I]$, $j\in[J]$ and $t\in [T]$, $S^{\psi^*,k}_{i,j}(t)$ are independently and identically distributed for all $k\in[h]$, and so as $\pmb{\Theta}^{\psi^*,k}_{i,j}(t)$. 
For $k\in[h]$, let $\bar{\pmb{\Theta}}^{\psi^*}_{i,j}(t) \coloneqq \mathbb{E}\pmb{\Theta}^{\psi^*,k}_{i,j}(t)$.
From the law of large numbers, for any $\delta>0$,
\begin{equation}\label{eqn:converge_z_opt:2}
\lim_{h\rightarrow +\infty}\mathbb{P}\biggl\{\Bigl\lVert \frac{1}{h}\sum_{k\in[h]}\bm{\Theta}^{\psi^*,k}_{i,j}(t) -\bar{\bm{\Theta}}^{\psi^*}_{i,j}(t) \Bigr\rVert > \delta\biggr\}=0.
\end{equation}
For $\iota\in[N]$ and $t = t_{\iota}$, the number of sub-processes in PST~$\iota$ at time $t$ is $IJh Z^{h,\psi^*}_{\iota}(t) = \sum_{k\in[h]}\Theta^{\psi^*,k}_{s_{\iota},i_{\iota},j_{\iota}}(t)$, and the expected number of the sub-processes is $IJh \mathbb{E} Z^{h,\psi(M,\bm{\phi}^*)}_{\iota}(t) = h\bar{\Theta}^{\psi^*}_{s_{\iota},i_{\iota},j_{\iota}}(t)$. Substituting these two equations into \eqref{eqn:converge_z_opt:2}, we obtain \eqref{eqn:converge_z_opt:1}.
The lemma is proved.
\endproof


\begin{lemma}\label{lemma:equivalence_3}
If the MAB-ML process is strictly indexable, then, for  $\kappa\in[K]$ and $\pmb{z} = \lim_{h\rightarrow +\infty}\mathbb{E}\bm{Z}^{h,\psi^*}(t)$,
\begin{equation}\label{eqn:equivalence_3}
    \lim_{h\rightarrow +\infty}\Bigl[\zeta^{h,\psi^*}_{\kappa}\bigl(\pmb{z}\bigr) - \zeta^{h,\rm IND}_{\kappa}\bigl(\pmb{z}\bigr)\Bigr] = 0,
\end{equation}
where the existence of $\lim_{h\rightarrow +\infty}\mathbb{E}\bm{Z}^{h,\psi^*}(t)$ is guaranteed by Lemma~\ref{lemma:indexability} and the definition of $\psi^*$.
\end{lemma}
\proofRemark{Lemma~\ref{lemma:equivalence_3}}
Based on Proposition~\ref{prop:strong_duality} and the definition of the $\psi^*$ policy, for any $h\in\mathbb{N}_+$, we obtain that, for $i\in[I]$, $j\in[J]$ and $t\in[T]$, 
\begin{multline}\label{eqn:equivalence_3:1}
\sum_{i'\in \bar{\mathscr{B}}_{i,j}}\sum_{k\in[h]}\mathcal{w}_{i',i',j} \mathbb{E}\Bigl[\alpha^{\psi^*,k}_{i',i,j}\bigl(S^{\psi^*,k}_{i,j}(t),t\bigr)\Bigr]  - \sum_{k\in[h]}\mathbb{E}\Bigl[g_{i,j}\bigl(S^{\psi^*,k}_{i,j}(t)\bigr)\Bigr]\\
=
\sum_{i'\in\bar{\mathscr{B}}_{i,j}}\sum_{m\in[M]}\sum_{k=k^-(m-1)+1}^{k^-(m)} \mathcal{w}_{i',i',j} \mathbb{E}\Bigl[\alpha^{\phi^*_m,k}_{i',i,j}\bigl(S^{\phi^*_m,k}_{i,j}(t),t\bigr)\Bigr] 
-\sum_{m\in[M]}\sum_{k=k^-(m-1)+1}^{k^-(m)}\mathbb{E}\Bigl[g_{i,j}\bigl(S^{\phi^*_m,k}_{i,j}(t)\bigr)\Bigr]\\
= \sum_{i'\in\bar{\mathscr{B}}_{i,j}}\sum_{m\in[M]} \pi^*_m h\mathcal{w}_{i',i,j} \mathbb{E}\Bigl[\alpha^{\phi^*_m,k_m}_{i',i,j}\bigl(S^{\phi^*_m,k_m}_{i,j}(t),t\bigr)\Bigr] -
\sum_{m\in[M]} \pi^*_m h \mathbb{E}\Bigl[g_{i,j}\bigl(S^{\psi^*_m,k_m}_{i,j}(t)\bigr)\Bigr]+ o(h)\\
= o(h),
\end{multline}
where $k_m$ is any integer in $\{k^-(m-1)+1,k^-(m-1)+2,\ldots,k^-(m)\}$, and the $o(h)$ represents a real number (possibly negative) that satisfies $\lim_{h\rightarrow \infty} o(h)/h = 0$.
The second equality holds because processes $\{S^{\phi^*_m,k}_{i,j}(t),t\in[T]\}$ with $k\in\{k^-(m-1)+1,k^-(m-1)+2,\ldots,k^-(m)\}$ are stochastically identical (that is, the processes have the same state and action spaces, transition rules under various actions, and cost functions.).
Recall that $\phi^*_m$ for all $m\in[M]$ satisfy \eqref{eqn:lemma:indexability} when replacing $\phi_{i,j}(\pmb{\gamma}_{i,j})$ and $\pmb{\gamma}_{i,j}$ with $\phi^*_m$ and $\pmb{\gamma}^*_{i,j}$, respectively, and that $\zeta^{h,\phi}_{\kappa}(\pmb{z})$  represents the proportion of the selected APSTs $\kappa\in [K]$ under policy $\phi\in\tilde{\Phi}^h$ at time $t$ given $\bm{Z}^{h,\phi}(t) = \pmb{z}$.

Let $\mathscr{K}_t\in 2^{[K]}$ represent the set of APST tuples $\kappa$ with $t_{\kappa}=t$. Note that, from the definition of APST tuples, the elements in $\mathscr{K}_t$ are successive integers. 
In this context, let $\mathscr{K}_t = \{\kappa_t+1,\kappa_t+2,\ldots,\kappa_t+ K_t\}$ where $K_t = |\mathscr{K}_t|$, and $\kappa_t = \sum_{\tau = T}^{t+1} K_{\tau}$ with $\kappa_T= 0$.
For any $h\in\mathbb{N}_+$, $t\in[T]$ and $\mathbb{E}\bm{Z}^{h,\psi^*}(t) = \pmb{z}\in [0,1]^N$, there exist $K^1_t$ and $ K^2_t \in \mathscr{K}_t$ such that, for any $\kappa\in\mathscr{K}_t$, 
\begin{enumerate}[label=(\roman*)]
\item \label{case:1} if $\kappa \leq \kappa_t+K^1_t$, then 
 $\eta_{i_{\kappa},i'_{\kappa},j_{\kappa}}(s_{\kappa},t) < 0$, and 
 \begin{equation} \label{eqn:equivalence_3:2}
     \zeta^{h,\psi^*}_{\kappa}(\pmb{z}) = \max\Bigl\{z_{\iota(\kappa)} - \sum_{
     \begin{subarray}
         ~\kappa' \in[K]:\\
         \kappa' < \kappa,\\
         i_{\kappa'} = i_{\kappa},\\
         j_{\kappa'} = j_{\kappa},\\
         s_{\kappa'} = s_{\kappa},\\
         t_{\kappa'} = t_{\kappa}
     \end{subarray}
     }\zeta^{h,\psi^*}_{\kappa'}(\pmb{z}),0\Bigr\};
 \end{equation}
\item \label{case:2}if $\kappa_t+K^1_t < \kappa \leq \kappa_t+K^2_t$, then $\eta_{i_{\kappa},i'_{\kappa},j_{\kappa}}(s_{\kappa},t) = 0$, and 
 \begin{multline}\label{eqn:equivalence_3:3}
    \!\!\!\!\!\!\!\!\!\!\zeta^{h,\psi^*}_{\kappa}(\pmb{z}) = \max\biggl\{
     \min\Bigl\{\sum_{
     \begin{subarray}
         ~\iota\in[N]:\\
         i_{\iota} = i'_{\kappa},\\
         j_{\iota} = j_{\kappa},\\
         t_{\iota} = t_{\kappa},
     \end{subarray}
     }z_{\iota}g_{i_{\iota},j_{\iota}}\bigl(s_{\iota}\bigr) +\frac{o(h)}{h}
     - \sum_{
     \begin{subarray}
         ~\kappa'\in[K]:\\
         \kappa' < \kappa,\\
         i'_{\kappa'} = i'_{\kappa},\\
         j_{\kappa'} = j_{\kappa}
     \end{subarray}
     }    \mathcal{w}_{i_{\kappa'},i'_{\kappa'},j_{\kappa'}}\zeta^{h,\psi^*}_{\kappa'}(\pmb{z}),
     ~z_{\iota(\kappa)} - \sum_{
     \begin{subarray}
         ~\kappa' \in[K]:\\
         \kappa' < \kappa,\\
         i_{\kappa'} = i_{\kappa},\\
         j_{\kappa'} = j_{\kappa},\\
         s_{\kappa'} = s_{\kappa},\\
         t_{\kappa'} = t_{\kappa}
     \end{subarray}
     }\zeta^{h,\psi^*}_{\kappa'}(\pmb{z})\Bigr\},0\biggr\},
\end{multline}
where recall that, for any $(i,j)\in[I]\times[J]$ and $i'\notin \bar{\mathscr{B}}_{i,j}$, $\mathcal{w}_{i',i,j}=0$;
\item \label{case:3}otherwise, $\eta_{i_{\kappa},i'_{\kappa},j_{\kappa}}(s_{\kappa},t) > 0 $, and $\zeta^{h,\psi^*}_{\kappa}(\pmb{z}) = 0$.
\end{enumerate}
For any $\kappa $ in Case~\ref{case:1}, because of \eqref{eqn:constraint:exclusive:h},  if \eqref{eqn:equivalence_3:2} is satisfied, then $\zeta^{h,\psi^*}_{\kappa}(\pmb{z}) = z_{\iota(\kappa)}$ (or equivalently, $a^{\psi^*,k}_{i_{\kappa},i'_{\kappa},j_{\kappa}}(s_{\kappa},t)=1$ for all $k\in[h]$) if and only if $i'_{\kappa} = \min\arg\max_{i'\in\mathscr{B}_{i_{\kappa},j_{\kappa}}}\theta_{i_{\kappa},i',j_{\kappa}}(\gamma^*_{i',j,t})$, which is consistent with \eqref{eqn:lemma:indexability}. 
Also, based on Proposition~\ref{prop:strong_duality} and the definition of $\psi^*$, we obtain that, for all $\kappa$ in Case~\ref{case:1},
\begin{equation}
    \sum_{
     \begin{subarray}
         ~\iota\in[N]:\\
         i_{\iota} = i'_{\kappa},\\
         j_{\iota} = j_{\kappa},\\
         t_{\iota} = t_{\kappa},\\
     \end{subarray}
     }z_{\iota}g_{i_{\iota},j_{\iota}}(s_{\iota}) +\frac{o(h)}{h}
     - \sum_{
     \begin{subarray}
         ~\kappa'\in[K]:\\
         \kappa' < \kappa,\\
         i'_{\kappa'} = i'_{\kappa},\\
         j_{\kappa'} = j_{\kappa}
     \end{subarray}
     }   \mathcal{w}_{i_{\kappa'},i'_{\kappa'},j_{\kappa'}} \zeta^{h,\psi^*}_{\kappa'}(\pmb{z}) \geq 0.
\end{equation}
For $\kappa$ in Case~\ref{case:2}, based on the strict indexability of each area, for any $i'\in[I]$, $j\in[J]$ and $t\in[T]$, there exists at most a $\kappa$ such that $i'_{\kappa} = i'$, $t_{\kappa} = t$, $j_{\kappa} = j$ and $\eta_{i_{\kappa},i'_{\kappa},j_{\kappa}}(s_{\kappa},t) = 0$. 
Together with \eqref{eqn:lemma:indexability} and \eqref{eqn:constraint:exclusive:h}, 
we obtain \eqref{eqn:equivalence_3:3}.
Case~\ref{case:3} is led by \eqref{eqn:lemma:indexability}.
Compare to $\zeta^{h,\rm IND}_{\kappa}(\pmb{z})$ described in \eqref{eqn:zeta}, it proves \eqref{eqn:equivalence_3}
for all $\kappa\in[K]$.
\endproof

\proofRemark{Proposition~\ref{prop:equivalence_1}}
For $t=1$, \eqref{eqn:equivalence_1} holds, and for any $\delta>0$, \eqref{eqn:equivalence_1:0} holds.
For all $\tau\in[t]$, assume \eqref{eqn:equivalence_1} and \eqref{eqn:equivalence_1:0} hold with substituted $\tau$ for $t$. Based on \eqref{eqn:converge_x:2}, 
\begin{equation}\label{eqn:equivalence_1:2}
\lim_{h\rightarrow +\infty}\mathbb{E}[\bm{Z}^{h,\rm IND}(t)] = \frac{\bar{\bm{x}}_t}{IJ} = \lim_{h\rightarrow +\infty}\mathbb{E}[\bm{Z}^{h,\psi^*}(t)]\\
= \lim_{h\rightarrow +\infty}\mathbb{E}[\bm{Z}^{h,\bar{\phi}}(t)].
\end{equation}
Let $\pmb{z} = \lim_{h\rightarrow +\infty}\mathbb{E}[\bm{Z}^{h,\psi^*}(t)] = \lim_{h\rightarrow +\infty}\mathbb{E}[\bm{Z}^{h,\rm IND}(t)]=\lim_{h\rightarrow +\infty}\mathbb{E}[\bm{Z}^{h,\bar{\phi}}(t)]$.
Since $\bar{\phi}$ satisfies \eqref{eqn:theorem:IND_asym_opt:1}, 
\begin{multline}
\sum_{\kappa\in[K]}\Bigl\lvert \zeta^{h,\text{IND}}_{\kappa}(\bm{Z}^{h,\text{IND}}(t)) -\zeta^{h,\bar{\phi}}_{\kappa}(\bm{Z}^{h,\bar{\phi}}(t))\Bigr\rvert  \\\leq 
\sum_{\kappa\in[K]}\Bigl\lvert \zeta^{h,\text{IND}}_{\kappa}(\bm{Z}^{h,\text{IND}}(t)) -\zeta^{h,\text{IND}}_{\kappa}(\bm{Z}^{h,\bar{\phi}}(t))\Bigr\rvert+
\sum_{\kappa\in[K]}\Bigl\lvert \zeta^{h,\text{IND}}_{\kappa}(\bm{Z}^{h,\bar{\phi}}(t))
-\zeta^{h,\bar{\phi}}_{\kappa}(\bm{Z}^{h,\bar{\phi}}(t))\Bigr\rvert. \\
=\sum_{\kappa\in[K]}\Bigl\lvert \zeta^{h,\text{IND}}_{\kappa}(\bm{Z}^{h,\text{IND}}(t)) -\zeta^{h,\text{IND}}_{\kappa}(\bm{Z}^{h,\bar{\phi}}(t))\Bigr\rvert
+\frac{o(h)}{h},
\end{multline}
where the equality is achieved because $\bar{\phi}$ satisfies \eqref{eqn:theorem:IND_asym_opt:1}.
Together with the definition of $\zeta^{h,\text{IND}}_{\kappa}(\pmb{z})$ in \eqref{eqn:zeta}, \eqref{eqn:converge_x:2}, the assumed \eqref{eqn:equivalence_1:0}, and \eqref{eqn:equivalence_1:2}, we obtain
\begin{equation}\label{eqn:equivalence_1:3_5}
 \sum_{\kappa\in[K]}\Bigl\lvert \zeta^{h,\text{IND}}_{\kappa}(\bm{Z}^{h,\text{IND}}(t)) -\zeta^{h,\bar{\phi}}_{\kappa}(\bm{Z}^{h,\bar{\phi}}(t))\Bigr\rvert  =
\frac{o(h)}{h}.   
\end{equation}

Since \eqref{eqn:equivalence_1:0} is assumed to be true for all $\tau \in[t]$,  \eqref{eqn:converge_x:2}, Lemma~\ref{lemma:converge_z_opt}, \eqref{eqn:equivalence_1:3_5}, and  Lemma~\ref{lemma:equivalence_3},
for $\kappa\in[K]$, $\phi=\rm IND$ and $\bar{\phi}$, and $\delta >0$,
\begin{equation}\label{eqn:equivalence_1:4}
    \lim_{h\rightarrow +\infty}\mathbb{P}\biggl\{\Bigl\lvert \zeta^{h,\rm \phi}_{\kappa}\bigl(\bm{Z}^{h,\phi}(t)\bigr) - \zeta^{h,\psi^*}_{\kappa}\bigl(\bm{Z}^{h,\psi^*}(t)\bigr)\Bigr\rvert > \delta\biggr\} = 0.
\end{equation}
It follows the same probability of activating a sub-process in each PST under the index policy, $\bar{\phi}$ and $\psi^*$ in the asymptotic regime. 
Hence, 
\begin{equation}\label{eqn:equivalence_1:5}
    \frac{\bar{\bm{x}}_{t+1}}{IJ} = \lim_{h\rightarrow +\infty}\mathbb{E}\bigl[\bm{Z}^{h,\rm IND}(t+1)\bigr] = \lim_{h\rightarrow +\infty}\mathbb{E}\bigl[\bm{Z}^{h,\bar{\phi}}(t+1)\bigr] \\= \lim_{h\rightarrow +\infty}\mathbb{E}[\bm{Z}^{h,\psi^*}(t+1)],
\end{equation}
where the first equality comes from \eqref{eqn:converge_x:2}. It remains to show that \eqref{eqn:equivalence_1:0} also holds when substituting $t+1$ with $t$.

For $t\in[T]$, given assumed \eqref{eqn:equivalence_1} and \eqref{eqn:equivalence_1:0} and proved \eqref{eqn:equivalence_1:4}, for each PST $\iota\in[N]$, there are at most $o(h)/h$ proportion of sub-processes that evolve differently under $\bar{\phi}$ and $\psi^*$.
In other words, for all $\iota\in[N]$, there are $hIJ Z^{h,\bar{\phi}}_{\iota}(t) - o(h)$ sub-processes that have the same transition probabilities as those under the policy $\psi^*$ in PST $\iota$. Since all the sub-processes evolve independently under the policy $\psi^*$, for any $\iota\in[N]$, we obtain that $Z^{h,\bar{\phi}}(t+1) = Z^{h,\psi^*}(t+1) + o(h)/h$, where $o(h)$ is a real number, possibly negative, with $\lim_{h\rightarrow +\infty}o(h)/h = 0$.
In this context,
for $\iota\in[N]$ and any $\delta>0$,
\begin{multline}\label{eqn:equivalence_1:6}
    \lim_{h\rightarrow +\infty}\mathbb{P}\biggl\{\Bigl\lvert Z^{h,\bar{\phi}}_{\iota}(t+1) - \mathbb{E}\bigl[Z^{h,\bar{\phi}}_{\iota}(t+1)\bigr] \Bigr\rvert > \delta\biggr\}
    = \lim_{h\rightarrow +\infty}\mathbb{P}\biggl\{\Bigl\lvert Z^{h,\psi^*}_{\iota}(t+1)   - \mathbb{E}\bigl[Z^{h,\psi^*}_{\iota}(t+1)\bigr]  + \frac{o(h)}{h}\Bigr\rvert > \delta\biggr\} \\
    \leq \lim_{h\rightarrow +\infty}\mathbb{P}\biggl\{\Bigl\lvert Z^{h,\psi^*}_{\iota}(t+1) -\mathbb{E}\bigl[Z^{h,\psi^*}_{\iota}(t+1)\bigr] \Bigr\rvert > \delta/2\biggr\} +  \lim_{h\rightarrow +\infty}\mathbb{P}\biggl\{\Bigl\lvert\frac{o(h)}{h}\Bigr\rvert > \delta/2\biggr\},
\end{multline}
where the first equality is derived by \eqref{eqn:equivalence_1:5}.Together with Lemma~\ref{lemma:converge_z_opt}, we obtain \eqref{eqn:equivalence_1:0} for the $t+1$ case. The theorem is proved. 
\endproof

\section{Proof of Lemma~\ref{lemma:MAI_asym_opt}}
\label{app:lemma:MAI_asym_opt}
\proofRemark{Lemma~\ref{lemma:MAI_asym_opt}}
Based on Lemma~\eqref{lemma:equivalence_3}, Proposition~\ref{prop:strong_duality} and the definition of $\psi^*$, for any $\kappa\in[K]$ and $\pmb{z}\in[0,1]^N$,
\begin{equation}
    \sum_{
    \begin{subarray}
    ~\kappa'\in[K]:\\
    \kappa' < \kappa,\\
    i'_{\kappa'} = i'_{\kappa},\\
    j_{\kappa'} = j_{\kappa}
    \end{subarray}
    }\mathcal{w}_{i_{\kappa'},i'_{\kappa'},j_{\kappa'}}\zeta^{h,\rm IND}_{\kappa} (\pmb{z})+\frac{o(h)}{h}
    =\sum_{
    \begin{subarray}
    ~\kappa'\in[K]:\\
    \kappa' < \kappa,\\
    i'_{\kappa'} = i'_{\kappa},\\
    j_{\kappa'} = j_{\kappa}
    \end{subarray}
    }\mathcal{w}_{i_{\kappa'},i'_{\kappa'},j_{\kappa'}}\zeta^{h,\psi^*}_{\kappa} (\pmb{z})\\
    =\sum_{\begin{subarray}
    ~\iota\in[N]:\\
    i_{\iota} = i'_{\kappa},\\
    j_{\iota} = j_{\kappa},\\
    t_{\iota} = t_{\kappa}
    \end{subarray}}z_{\iota}g_{i_{\iota},j_{\iota}}(s_{\iota}) + \frac{o(h)}{h},
\end{equation}
where recall that, for any $(i,j)\in[I]\times[J]$ and $i'\notin \bar{\mathscr{B}}_{i,j}$, $\mathcal{w}_{i',i,j}=0$.
In this context, for any $\bm{Z}^{h,\text{MAI}}(t) = \bm{Z}^{h,\text{IND}}(t) = \pmb{z}$, based on the definition of MAI in  Section~\ref{subsubsec:MAI}, the total proportion of the sub-processes with adapted action variables under the MAI policy is $o(h)/h$. That is, for $t\in[T]$, any $\pmb{s}\in\mathscr{S}^h$ and its associated proportions of sub-processes $\pmb{z}(\pmb{s},t)\in[0,1]^N$ (see the relationship in \eqref{define:proportion_subprocesses}),
\begin{equation}\label{eqn:lemma:MAI_asym_opt:0}
    \frac{1}{h}\lVert \pmb{a}^{\text{MAI},h}(\pmb{s},t) - \pmb{a}^{\text{IND},h}(\pmb{s},t)\rVert^2\\
    =\sum_{\kappa\in[K]}\Bigl\lvert\zeta^{h,\text{IND}}_{\kappa}(\pmb{z}(\pmb{s},t)) - \zeta^{h,\text{MAI}}_{\kappa}(\pmb{z}(\pmb{s},t))\Bigr\rvert = \frac{o(h)}{h}.
\end{equation} 
which proves the lemma.\endproof



\section{Proof of Theorem~\ref{theorem:exp_convergence}}
\label{app:theorem:exp_convergence}

For $\bm{x},\bm{\omega}\in\mathbb{R}^N$, define
\begin{equation}\label{eqn:definition_H}
H(\bm{x},\bm{\omega}) \coloneqq \lim_{T\rightarrow \infty} \frac{1}{T}\ln \mathbb{E} exp \biggl\{\int_0^T \Bigl<\bm{\omega},b(\bm{x},\bm{\xi}^1_t)\Bigr> dt\biggr\},
\end{equation}
where $<\cdot,\cdot>$ is the dot product of two vectors. 
Since $b(\cdot,\cdot)$ is Lipschitz continuous in both arguments, we obtain that $b(\bm{x},\bm{\xi}^1_t) = \bar{\mathcal{Q}}(\bm{x}) \bm{\xi}^1_t$ is bounded and Lipschitz continuous in both arguments.
From \cite[Lemma 4.1, Chapter 7]{freidlin2012random}, $H$ is jointly continuous in both arguments and convex in the second argument.
\begin{lemma}\label{lemma:existance_H}
For any $T\in\mathbb{R}_+$, and vectors $\bm{x}_t,\bm{\omega}_t\in\mathbb{R}^N$ for $t\in[0,T]$, the $H$ defined in \eqref{eqn:definition_H} satisfies
\begin{equation}\label{eqn:existance_H}
\int_0^T H(\bm{x}_t,\bm{\omega}_t) dt =\lim_{\epsilon \rightarrow 0} \ln \mathbb{E} exp\biggl\{\frac{1}{\epsilon}\int_0^T \Bigl<\bm{\omega}_t, b(\bm{x}_t,\bm{\xi}^1_{t/\epsilon})\Bigr>d t\biggr\}.
\end{equation}
\end{lemma}
\proofRemark{Lemma~\ref{lemma:existance_H}}
For any $\bm{w}\in \mathbb{R}^{2N^2}$, define
\begin{equation}\label{eqn:existance_H:1}
H_{\xi}(\bm{w}) \coloneqq \lim_{T\rightarrow \infty} \frac{1}{T}\ln \mathbb{E} exp\biggl\{\int_0^T \Bigl<\bm{w},\bm{\xi}^1_t\Bigl> dt\biggr\}.
\end{equation}
Recall that the integral \highlight{$\lfloor \int_0^T\bm{\xi}^1_t dt \rfloor$} is a random vector, Poisson distributed with some parameter $\bm{\lambda}$ that are determined by the probabilities of $S^{\phi,k}_{i,j}(t)$ ($(i,j)\in[I]\times[J]$ with any $k\in[h]$) transitioning between different states in $\mathscr{S}_{i,j}$ under given actions.
It follows that 
\begin{equation}\label{eqn:existence_H:2}
H_{\xi}(\bm{w}) = \sum_{\ell=1}^{2N^2} \lambda_{\ell} (e^{w_{\ell}}-1),
\end{equation}
where $e$ is the Euler number. For any bounded $\bm{\omega}$, $H_{\xi}(\bm{\omega})$ is bounded.
From \eqref{eqn:definition_H} and, for any $\bm{x}\in\mathbb{R}^N$, $b(\bm{x},\bm{\xi}^1_t) = \bar{\mathcal{Q}}(\bm{x}) \bm{\xi}^1_t$, we obtain that 
\begin{equation}\label{eqn:continuity_in_derivative}
    H(\bm{x},\bm{\omega})= H_{\xi}\Bigl(\bm{\omega}\bar{\mathcal{Q}}(\bm{x})\Bigr)<\infty.
\end{equation}
For any compact set $\mathscr{A}^c\in\mathbb{R}^N$ and $\bm{\omega}\in\mathscr{A}^c$, $H(\bm{x},\bm{\omega})$ obtained from \eqref{eqn:definition_H} is bounded and, because of its joint continuity, is Riemann integrable. Thus, for  $T\in\mathbb{R}_+$ and vectors $\bm{x}_t,\bm{\omega}_t\in\mathbb{R}^N$ for $t\in[0,T]$, it satisfies
\eqref{eqn:existance_H}.
\endproof

For $\bm{x},\bm{\beta}\in\mathbb{R}^N$, define the Legendre transform of $H(\bm{x},\bm{\omega})$
\begin{equation}\label{eqn:legendre_L}
L(\bm{x},\bm{\beta}) \coloneqq \sup\nolimits_{\bm{\omega}\in\mathbb{R}^N} \Bigl[<\bm{\omega},\bm{\beta}> - H(\bm{x},\bm{\omega})\Bigr],
\end{equation}
which is non-negative, because $<\bm{\omega},\bm{\beta}> - H(\bm{x},\bm{\omega}) = 0$ when $\bm{\omega} = \bm{0}$.

\begin{lemma}\label{lemma:unique_zero_L}
For $\bm{x}\in\mathbb{R}^N$ and any $\delta > 0$,  
if $\lVert \bm{\beta} - \mathbb{E} b(\bm{x},\bm{\xi}^1_t)\rVert\geq \delta$ then $L(\bm{x},\bm{\beta}) > 0$.
\end{lemma}
\proofRemark{Lemma~\ref{lemma:unique_zero_L}}
From \cite[Chapter 7, Section 4]{freidlin2012random},  if $\bm{\beta} = b(\bm{x},\bm{\xi}^1_t)$, then $L(\bm{x},\bm{\beta}) = 0$. It remains to show that there exists a unique $\bm{\beta}$ such that $L(\bm{x},\bm{\beta}) = 0$.

Rewrite \eqref{eqn:legendre_L} as, for $\bm{x},\bm{\beta}\in\mathbb{R}^N$,
\begin{equation}\label{eqn:L}
L(\bm{x},\bm{\beta}) = <\omega(\bm{x},\bm{\beta}), \bm{\beta}> - H(\bm{x},\omega(\bm{x},\bm{\beta})),
\end{equation}
where $\omega(\bm{x},\bm{\beta}) \in \mathbb{R}^N$ represents an extreme point satisfying $\bm{\beta} = \frac{\partial H}{\partial \bm{\omega}}\Bigl|_{\bm{\omega}=\omega(\bm{x},\bm{\beta})}$. Note that such extreme point is not necessarily unique.
For any $\bm{x},\bm{\beta}\in\mathbb{R}^N$, 
\begin{equation}\label{eqn:unique_zero_L:1}
L(\bm{x},\bm{\beta})  
=\Bigl[<\omega(\bm{x},\bm{\beta}),\bm{\beta}> - H_{\xi}(\omega(\bm{x},\bm{\beta})\bar{\mathcal{Q}}(\bm{x}))\Bigr],
\end{equation}
where the equality is a direct result of \eqref{eqn:continuity_in_derivative}.

Substitute \eqref{eqn:existence_H:2} in \eqref{eqn:unique_zero_L:1}, we obtain, for $\bm{x},\bm{\beta}\in\mathbb{R}^N$ satisfying $L(\bm{x},\bm{\beta})=0$,
\begin{equation}
\left\{
\begin{array}{l}
\sum_{\iota=1}^N \omega_{\iota}(\bm{x},\bm{\beta})\beta_{\iota} = \sum_{\ell=1}^{2N^2} \lambda_{\ell}\Bigl(e^{\bigl<\omega(\bm{x},\bm{\beta}),\bm{q}_{\ell}\bigr>}-1\Bigr),\\
\beta_{\iota} = \sum_{\ell=1}^{2N^2} \lambda_{\ell}q_{\iota,\ell}e^{\bigl<\omega(\bm{x},\bm{\beta}),\bm{q}_{\ell}\bigr>}\text{, for all $\iota\in[N]$,}
\end{array}\right.
\end{equation}
where $\bm{q}_{\ell}$ is the $\ell$th column vector of the matrix $\bar{\mathcal{Q}}(\bm{x})$.
It follows that
\begin{equation}
\sum_{\ell=1}^{2N^2} \lambda_{\ell}\biggl(e^{\bigl<\omega(\bm{x},\bm{\beta}),\bm{q}_{\ell}\bigr>} \Bigl(\bigl<\omega(\bm{x},\bm{\beta}),\bm{q}_{\ell}\bigr>-1\Bigr) +1 \biggr) = 0.
\end{equation}
Define a function $f(z) = \sum_{\ell=1}^{2N^2} \lambda_{\ell}(e^{z} (z-1) +1 )$. By taking the first derivative of $f(z)$, we observe that $f(z) \geq 0 $ and $f(z)=0$ if and only if $z=0$. That is, for given $\bm{x}\in\mathbb{R}^N$, $L(\bm{x},\bm{\beta}) = 0$ only if $\bigl<\omega(\bm{x},\bm{\beta}),\bm{q}_{\ell}\bigr> = 0$ and $\beta_{\iota} = \sum_{\ell=1}^{2N^2} \lambda_{\ell}q_{\iota,\ell}$ for all $\iota\in[N]$. In other words, for given $\bm{x}\in\mathbb{R}^N$, there is a unique $\bm{\beta}\in\mathbb{R}^N$ satisfying $L(\bm{x},\bm{\beta}) = 0$. This proves the lemma.
\endproof

\proofRemark{Theorem~\ref{theorem:exp_convergence}}
Define $S_{0,T}(\varphi) = \int_0^T L(\bm{\varphi}_t,\dot{\bm{\varphi}})dt $, where $\varphi$ represents a trajectory $\{\bm{\varphi}_t, 0\leq t \leq T\}$ with $\bm{\varphi}_t\in\mathbb{R}^N$, and let $C_{0,T}$ represent the compact set of all such trajectories with given initial state $\bm{\varphi}_0 = \bm{x}_0 \in \mathbb{R}^N$. Define a closed set $\mathscr{C}(\delta) \coloneqq \{\varphi \in C_{0,T}~|~\sup_{0\leq t\leq T}\lVert \bm{\varphi}_t-\bar{\bm{x}}_t\rVert \geq \delta\}$ where $\bar{\bm{x}}_t$ is the solution of $\dot{\bar{\bm{x}}}_t = \mathbb{E} b(\bar{\bm{x}}_t, \bm{\xi}^1_t)$ and $\bar{\bm{x}}_0 = \bm{x}_0$. 

From \cite[Theorem 4.1 in Chapter 7 and Theorem 3.3 in Chapter 3]{freidlin2012random}, for any $\delta>0$, 
\begin{equation}\label{eqn:exp_convergence:2}
\overline{\lim_{\sigma\downarrow 0}}~\sigma \ln \mathbb{P}\Bigl\{\sup_{0\leq t\leq T}\lVert \bm{X}^{\sigma}_t-\bar{\bm{x}}_t\rVert > \delta\Bigr\}
\leq -\inf_{\varphi \in \mathscr{C}(\delta)} S_{0,T}(\varphi),
\end{equation}
where process $\bm{X}^{\sigma}_t$ is the solution of $\dot{\bm{X}}^{\sigma}_t = b(\bm{X}^{\sigma,}_t,\bm{\xi}^1_{t/\sigma})$ and $\dot{\bm{X}}^{\sigma}_0 = \bm{x}_0\in\mathbb{R}^N$. 
Based on Lemma~\ref{lemma:unique_zero_L}, for any $T>0$, $\delta>0$, $\inf\nolimits_{\varphi \in \mathscr{C}(\delta)}S_{0,T} > 0$.
It follows that, for any $\delta>0$, there exists $\sigma_0>0$ and $s>0$ such that, for all $0<\sigma<\sigma_0$, 
\begin{equation}
\mathbb{P}\Bigl\{\sup_{0\leq t\leq T}\lVert\bm{X}^{\sigma}_t - \bar{\bm{x}}_t\rVert > \delta \Bigr\} \leq e^{-\frac{s}{\sigma}}.
\end{equation}
In other words, for any $\delta > 0$,
there exists $H>0$ and $s>0$ such that, for all $h>H$,
\begin{multline}\label{eqn:exp_convergence:3}
    \mathbb{P}\Bigl\{\sup_{0\leq t\leq T}\lVert\bm{Z}^{h,\rm IND}(t) - \frac{\bar{\bm{x}}_t}{IJ}\rVert > \delta\Bigr\}
=   \mathbb{P}\Bigl\{\sup_{0\leq t\leq T}\lVert\bm{Z}^h_t - \frac{\bar{\bm{x}}_t}{IJ}\rVert > \delta\Bigr\}\\
= \mathbb{P}\Bigl\{\sup_{0\leq t\leq T}\lVert\bm{X}^{\sigma}_t - \bar{\bm{x}}_t\rVert > IJ\delta \Bigr\} 
\leq e^{-sh},
\end{multline}
where  the inequality is achieved by setting $\sigma = \frac{1}{h}$.

Let $\bar{\phi}\in\Phi^h$ be an policy satisfying \eqref{eqn:constraint:exclusive:h}, \eqref{eqn:constraint:neighbourhood:h}, and \eqref{eqn:theorem:IND_asym_opt:1}.
Based on Proposition~\ref{prop:equivalence_1} and \eqref{eqn:converge_x:2},
for $t\in[T]$ and $\iota\in[N]$, 
\begin{equation}
    \lim_{h\rightarrow +\infty}\Bigl\lvert Z^{h,\bar{\phi}}_{\iota}(t) - Z^{h,\rm IND}_{\iota}(t)\Bigr\rvert = 0.
\end{equation}
Hence, we can write $Z^{h,\bar{\phi}}_{\iota}(t) = Z^{h,\rm IND}_{\iota}(t) + o(h)/h$. 
There exists $\pmb{o}_t(h) \coloneqq (o_{t,\iota}(h): \iota\in[N])$, for which $o_{t,\iota}(h)\in\mathbb{R}$ and $\lim_{h\rightarrow +\infty}o_{t,\iota}(h)/h = 0$, such that $\bm{Z}^{h,\bar{\phi}}(t) = \bm{Z}^{h, \rm IND}(t) +\bm{o}_t(h)/h$.
For $t\in[T]$ and any $\delta>0$,
\begin{multline}\label{eqn:exp_convergence:4}
    \mathbb{P}\biggl\{\Bigl\lVert \bm{Z}^{h,\bar{\phi}}(t) - \lim_{h\rightarrow +\infty}\mathbb{E}\bigl[\bm{Z}^{h,\bar{\phi}}(t)\bigr] \Bigr\rVert > \delta\biggr\}
    = \mathbb{P}\biggl\{\Bigl\lVert \bm{Z}^{h,\rm IND}(t) -\frac{\bar{\bm{x}}_t}{IJ} + \frac{1}{h}\bm{o}_t(h)\Bigr\rVert > \delta\biggr\} \\
    \leq \mathbb{P}\biggl\{\Bigl\lVert \bm{Z}^{h,\rm IND}(t) -\frac{\bar{\bm{x}}_t}{IJ} \Bigr\rVert + \Bigl\lVert \frac{1}{h}\bm{o}_t(h) \Bigr\rVert> \delta\biggr\}
    \leq \mathbb{P}\biggl\{\Bigl\lVert \bm{Z}^{h,\rm IND}(t) -\frac{\bar{\bm{x}}_t}{IJ}\Bigr\rVert > \frac{\delta}{2}\biggr\}+ \mathbb{P}\biggl\{\Bigl\lVert\frac{1}{h}\bm{o}_t(h)\Bigr\rVert > \frac{\delta}{2}\biggr\},
\end{multline}
where the first equality is based on \eqref{eqn:converge_x:2} and Proposition~\ref{prop:equivalence_1}.
Then, together with \eqref{eqn:exp_convergence:3}, for any $\delta>0$, there exists $H>0$ such that, for all $h>H$,
\begin{equation}
    \mathbb{P}\Bigl\{\sup_{0\leq t\leq T}\bigl\lVert\bm{Z}^{h,\bar{\phi}}(t) - \lim_{h\rightarrow+\infty}\mathbb{E}[\bm{Z}^{h,\bar{\phi}}(t)]\bigr\rVert > \delta \Bigr\} \\
    \leq \mathbb{P}\Bigl\{\sup_{0\leq t\leq T}\lVert\bm{Z}^{h,\rm IND}(t) - \frac{\bar{\bm{x}}_t}{IJ}\rVert > \frac{\delta}{2}\Bigr\} \leq e^{-sh}.
\end{equation}
Together with Proposition~\ref{prop:equivalence_1}, it leads to \eqref{eqn:exp_convergence:1}.
We prove the theorem.
\endproof

\begin{table}[t]
\centering
	\caption{Initialization of $\alpha^{\phi,k}_{i,j}(1)=\alpha^0_{i,j}$ and $\beta^{\phi,k}_{i,j}(1)=\beta^0_{i,j}$.}\label{table:alpha_beta}
	\begin{tabular}{c|cc|cc}
		\hline
		 & \multicolumn{2}{c}{$j=1$} &\multicolumn{2}{c}{$j=2$}\\ 
		\hline
         $i$& $\alpha^0_{i,1}$ & $\beta^0_{i,1}$ & $\alpha^0_{i,2}$ & $\beta^0_{i,2}$\\
         \hline
         1 & 2 & 48 & 4 & 46\\
         2 & 2 & 48 & 4 & 46 \\
         3 & 4 & 46 & 32 & 18 \\
         4 & 3 & 47 & 33 & 17 \\
         5 & 4 & 46 & 50 & 0 \\
         6 & 3 & 47 & 34 & 16 \\
         7 & 4 & 46 & 26 & 24\\
         8 & 3 & 47 & 11 & 39\\
         9 & 3 & 47 & 9 & 41 \\
         10 & 5 & 45 & 10 & 40\\
         11 & 4 & 46 & 22 & 28\\
         12 & 2 & 48 & 22 & 28\\
         13 & 2 & 48 & 22 & 46 \\
         14 & 4 & 46 & 25 & 25 \\
		\hline
	\end{tabular}
\end{table}
\begin{table*}[t]
\centering
	\caption{Initial probabilities $\pi^0_{i,j}(s_{i,j})$ for all $(i,j)\in[I]\times[J]$ and $s_{i,j}\in\mathscr{S}_{i,j}$ with $\mathcal{k}(s_{i,j})=(\alpha^0_{i,j},\beta^0_{i,j})$.}\label{table:init_probs}
	\begin{tabular}{c|cc|cc|cc}
        \hline
        \multicolumn{7}{c}{Type $j=1$}\\
		\hline
		& \multicolumn{2}{c}{Case I} &\multicolumn{2}{c}{Case II}&\multicolumn{2}{c}{Case III}\\ 
		\hline
         $i$& $\mathcal{i}(s_{i,j})=0$ & $\mathcal{i}(s_{i,j})=1$ & $\mathcal{i}(s_{i,j})=0$ & $\mathcal{i}(s_{i,j})=1$&$\mathcal{i}(s_{i,j})=0$ & $\mathcal{i}(s_{i,j})=1$\\
         \hline
         1 & 0.9 & 0.1 & 0.95 & 0.05 & 0.75 & 0.25\\
         2 & 0.75 & 0.25 & 0.9 & 0.1 & 1 & 0 \\
         3 & 0.9 & 0.1 & 0.55 & 0.45 & 0.85 & 0.15 \\
         4 & 0.85 & 0.15 & 0.85 & 0.15 & 0.95 & 0.05 \\
         5 & 0.8 & 0.2 & 0.9 & 0.1 & 0.85 & 0.15 \\
         6 & 0.8 & 0.2 & 0.8 & 0.2 & 0.8 & 0.2 \\
         7 &&& 0.85 & 0.15 & 0.85 & 0.15\\
         8 &&& 0.6 & 0.4 & 0.9 & 0.1\\
         9 &&& 0.7 & 0.3 & 0.9 & 0.1 \\
         10 &&& 0.9 & 0.1 & 0.85 & 0.15\\
         11 &&&  & & 0.75 & 0.25\\
         12 &&& & & 0.8 & 0.2\\
         13 &&& & & 0.9 & 0.1 \\
         14 &&& & & 0.85 & 0.15 \\
		\hline
   \hline
        \multicolumn{7}{c}{Type $j=2$}\\
		\hline
		 & \multicolumn{2}{c}{Case I} &\multicolumn{2}{c}{Case II}&\multicolumn{2}{c}{Case III}\\ 
		\hline
         $i$& $\mathcal{i}(s_{i,j})=0$ & $\mathcal{i}(s_{i,j})=1$ & $\mathcal{i}(s_{i,j})=0$ & $\mathcal{i}(s_{i,j})=1$&$\mathcal{i}(s_{i,j})=0$ & $\mathcal{i}(s_{i,j})=1$\\
         \hline
         1 & 0.5 & 0.5 & 0.7 & 0.3 & 0.65 & 0.35\\
         2 & 0.5 & 0.5 & 0.4 & 0.6 & 0.65 & 0.35 \\
         3 & 0.35 & 0.65 & 0.7 & 0.3 & 0.75 & 0.25 \\
         4 & 0.8 & 0.2 & 0.45 & 0.55 & 0.75 & 0.25 \\
         5 & 0.55 & 0.45 & 0.55 & 0.45 & 0.65 & 0.35 \\
         6 & 0.3 & 0.7 & 0.25 & 0.75 & 0.7 & 0.3 \\
         7 &&& 0.35 & 0.65 & 0.55 & 0.45\\
         8 &&& 0.4 & 0.6 & 0.8 & 0.2\\
         9 &&& 0.65 & 0.35 & 0.65 & 0.35 \\
         10 &&& 0.55 & 0.45 & 0.7 & 0.3\\
         11 &&&  &  & 0.65 & 0.35\\
         12 &&&  &  & 0.55 & 0.45\\
         13 &&&  &  & 0.55 & 0.45 \\
         14 &&&  &  & 0.4 & 0.6 \\
         \hline
	\end{tabular}
\end{table*}
\section{Policies in the Scaled System}\label{app:algorithms:h}
For $i,i'\in[I]$, $j\in[J]$ and $k\in[h]$, define $(i,k,i',j)$ as an action (movement) for process $\Bigl\{S^{\phi,k}_{i,j}(t),t\in[T]\bigr\}$ in the scaled system. 
Similar to the $\mathscr{M}^{\phi}_j(t)$ defined in Section~\ref{subsubsec:index_policy}, at each time slot $t\in[T]$, for a policy $\phi\in\Phi$ and $j\in[J]$, we maintain a set $\mathscr{M}^{\phi,h}_j(t)$ of all the actions $(i,k,i',j)$ such that $\sum\nolimits_{k'\in[h]}g_{i',j}(S^{\phi,k'}_{i',j}(t))>0$ and $i'\in\bar{\mathscr{B}}_{i,j}$. 
We rank all actions in $\mathscr{M}^{\phi,h}_j(t)$ in the same rule as the \partialref{define:movement_ranking}{action ranking} in Section~\ref{subsubsec:index_policy} by replacing $(i,i',j)$, $\mathscr{M}^{\phi}_j(t)$ and $S^{\phi}_{i,j}(t)$ with $(i,k,i',j)$, $\mathscr{M}^{\phi,h}_j(t)$ and $S^{\phi,k}_{i,j}(t)$, respectively.
We write the $\mathcal{r}$th action in this ranking as $(i_{\mathcal{r},j},k_{\mathcal{r},j},i'_{\mathcal{r},j},j)$.
Based on the ranked actions in $\mathscr{M}^{\text{IND},h}_j(t)$, the pseudo-code for implementing the index policy in the scaled system is provided in Algorithm~\ref{algo:index_policy:h}.
\IncMargin{1em}
\begin{algorithm}\small
\linespread{0.5}\selectfont

\SetKwProg{Fn}{Function}{}{End}
\SetKwInOut{Input}{Input}\SetKwInOut{Output}{Output}
\SetAlgoLined
\DontPrintSemicolon

\Input{Ranked actions for each $j\in[J]$ and $\bm{S}^{\text{IND},h}(t)$.}
\Output{$a^{\text{IND},k}_{i,i',j}\bigl(\bm{S}^{\rm IND}(t),t\bigr)$ for all $ (i,k,i',j)\in \bigcup_{j\in[J]}\mathscr{M}^{\text{IND},h}_j(t)$.}

\Fn{ScaledIndexPolicy}{
	$a^{\text{IND},k}_{i,i',j}\bigl(\bm{S}^{\text{IND},h}(t),t\bigr)\gets 0$ for all $(i,k,i',j)\in\bigcup_{j\in[J]}\mathscr{M}^{\text{IND},h}_j(t)$\;
	\For{$j=1,2,\ldots, J$}{
        $q(i,k)\gets 0$ for all $(i,k)\in\bigl\{(i,k)| (i,k,i',j)\in\mathscr{M}^{\text{IND},h}_j(t)\bigr\}$\;
        $p(i')\gets \sum\nolimits_{k'\in[h]}g_{i',j}(s^{k'}_{i',j})$ for all $i'\in\bigl\{i'| (i,k,i',j)\in\mathscr{M}^{\text{IND},h}_j(t)\bigr\}$\;
        \For{$\mathcal{r}=1,2,\ldots, |\mathscr{M}^{\text{IND},h}_j(t)|$}{
            $(i,k,i')\gets (i_{\mathcal{r},j},k_{\mathcal{r},j},i'_{\mathcal{r},j})$\;
	 	    \If{($i\notin\bar{\mathscr{B}}_{i',j}$ or $p(i') \geq \mathcal{w}_{i,i',j}$)  and $q(i,k)=0$\label{algo:availability:h}}{                      
                        $a^{\text{ IND},k}_{i,i',j}\bigl(\bm{S}^{\text{IND},h}(t),t\bigr) \gets 1$ \;
                        $q(i,k) \gets 1$\;
                        $p(i') \gets p(i') -\mathcal{w}_{i,i',j}\mathds{1}\{i\in\bar{\mathscr{B}}_{i',j}\}$\;
                }
            }
    }
    \Return $a^{\text{IND},k}_{i,i',j}\bigl(\bm{S}^{\text{IND},h}(t),t\bigr)$ for all $(i,k,i',j)\in\bigcup_{j\in[J]}\mathscr{M}^{\text{IND},h}_j(t)$\;
}
\caption{Pseudo-code for implementing the index policy in the scaled system.}\label{algo:index_policy:h}
\end{algorithm}
\DecMargin{1em}
\begin{figure*}[t]
\centering
\subfigure[]{\includegraphics[width=0.22\linewidth]{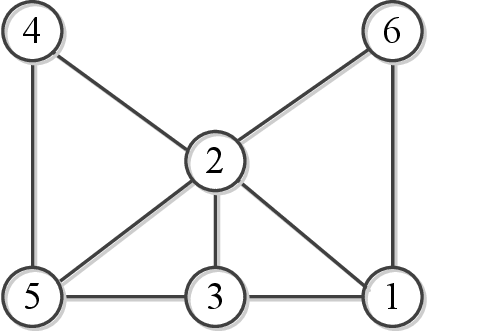}}
\subfigure[]{\includegraphics[width=0.33\linewidth]{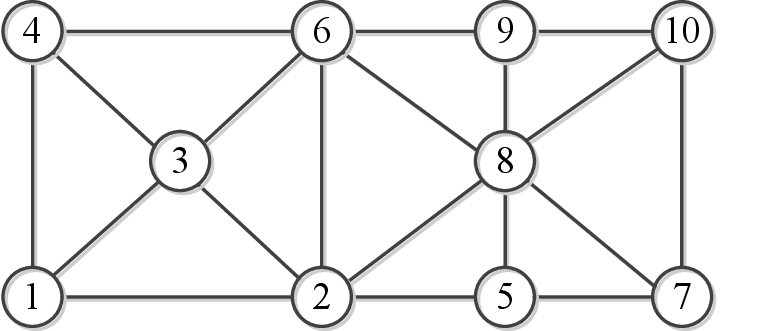}}
\subfigure[]{\includegraphics[width=0.35\linewidth]{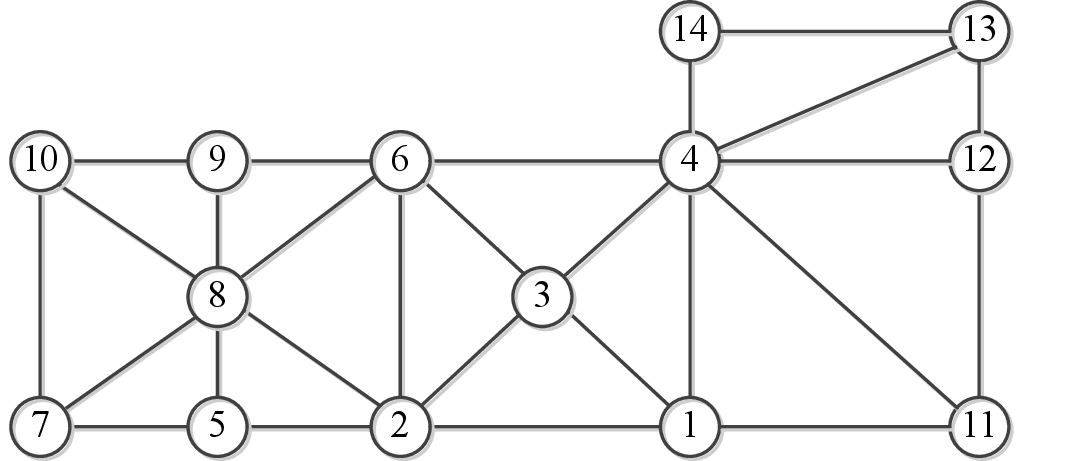}}
\caption{Topology of patrol region for (a) Case~I ($I=6$), (b) Case~II ($I=10$), and (c) Case~III ($I=14$).\label{fig:topology}}
\end{figure*}

Similarly, in the scaled system for patrol, given the action vector $\pmb{a}^{\phi,h}(\bm{S}^{\phi,h}(t),t) = \pmb{a}^h$, define the origin of an agent that moves to sub-area $(i,k)$ ($i\in[I],k\in[h]$) as $\varpi_{i,j}^k(\pmb{a}^h)$. Let $\pmb{\varpi}^h\coloneqq (\varpi_{i,j}^k: (i,j)\in[i]\times[J],k\in[h])$. 
Define $\mathscr{I}^{\phi,h}_j(t)$ the set of all $i\in[I]$ such that  $\sum_{k\in[h]}g_{i,j}(S^{\phi,k}_{i,j}(t))>0$ but \eqref{eqn:constraint:neighbourhood:h} is not satisfied.
Similar to \eqref{eqn:distance} for the special case, the distance of an area $i$ for type-$j$ agent with respect to $\bm{S}^{\phi,h}(t)=\pmb{s}^h$ and $\pmb{\varpi}^h$ is defined as 
\begin{equation}\label{eqn:distance:h}
    d^h_{i,j}(\pmb{s}^h,\pmb{\varpi}^h) = \left\{\begin{cases}
        0,& \text{if }\sum_{k\in[h]}\mathcal{i}(s^k_{i,j}) > 0,\text{ and } \exists (i',k')\in\mathscr{B}_{i,j}\times[h],\varpi^{k'}_{i',j} = 0,\\       \min_{\begin{subarray}~i'\in\mathscr{B}_{i,j}:\\\exists k'\in[h],\\~\varpi^{k'}_{i',j}\neq i\end{subarray}}d_{\varpi_{i',j}}+ 1, & \text{if }\sum_{k\in[h]}\mathcal{i}(s^k_{i,j}) >0,\text{ and } \forall (i',k')\in\mathscr{B}_{i,j}\times[h],\varpi^{k'}_{i',j} > 0, 
        \\
        \infty, &\text{otherwise}.
    \end{cases}\right.
\end{equation}
The implementation of the movement adaption process is in the same vein as Algorithm~\ref{algo:movement_adaption}, and
we propose in Algorithm~\ref{algo:movement_adaption:h} the pseudo-code for the movement adaption in the scaled system.

\IncMargin{1em}
\begin{algorithm*}\small
\linespread{1.5}\selectfont

\SetKwProg{Fn}{Function}{}{End}
\SetKwInOut{Input}{Input}\SetKwInOut{Output}{Output}
\SetAlgoLined
\DontPrintSemicolon

\Input{The state vector $\bm{S}^{\text{MAI},h}(t) = \pmb{s}^h$ at time $t$.}
\Output{The action variables $a^{\text{MAI},k}_{i,i',j}\bigl(\bm{S}^{\text{MAI},h}(t),t\bigr)$ for all sub-movements $(i,k,i',j)\in\bigcup_{j\in[J]}\mathscr{M}^{\text{MAI},h}_j(t)$ for the MAI policy.}

\Fn{ScaledMAIPolicy}{
	$a^{\text{MAI},k}_{i,i',j}(\pmb{s}^h,t)\gets a^{\text{IND},k}_{i,i',j}(\bm{s^h},t)$ for all $(i,k,i',j)\in\bigcup_{j\in[J]}\mathscr{M}^{\text{MAI},h}_j(t)$\; \tcc*{Initialize the action variables by calling  Algorithm~\ref{algo:index_policy:h}.}
    Based on the action variables $\pmb{a}^{\rm MAI,h}(\pmb{s}^h,t)$, initialize $\mathscr{I}^{\rm MAI,h}_j(t)$ for all $j\in[J]$, $\varpi^k_{i,j}\bigl(\pmb{a}^{\rm MAI,h}(\pmb{s}^h,t)\bigr)$ for all $(i,j,k)\in[I]\times[J]\times[h]$, and $d^h_{i,j}\Bigl(\pmb{s}^h,\pmb{\varpi}^h\bigl(\pmb{a}^{\text{MAI},h}(\pmb{s}^h,t)\bigr)\Bigr)$ described in \eqref{eqn:distance:h} for all $(i,j)\in[I]\times [J]$.\;
    $\varpi^k_{i,j}\gets \varpi^k_{i,j}\bigl(\pmb{a}^{\text{MAI},h}(\pmb{s}^h,t)\bigr)$ for all $(i,j,k)\in[I] \times [J]\times[h]$\;
    $d^h_{i,j}\gets d^h_{i,j}\Bigl(\pmb{s}^h,\pmb{\varpi}^h\bigl(\pmb{a}^{\text{MAI},h}(\pmb{s}^h,t)\bigr)\Bigr)$ for all $(i,j)\in[I] \times [J]$\;
	\For{$j=1,2,\ldots,J$}{
        \For{$\forall i\in\mathscr{I}^{\text{MAI},h}_j(t)$}{
            $K\gets \sum\nolimits_{k\in[h]}\mathcal{i}(s^k_{i,j}) - \sum\nolimits_{s'\in\mathscr{B}_{i,j}}\sum\nolimits_{k'\in[h]}a^{\text{MAI},k'}_{i',i,j}(\pmb{s}^h,t)$\;
            \tcc*{The number of agents in area $i$ that awaiting to be moved.}
            \For{$k=1,2,\ldots,K$}{
                $\bar{i}\gets i$\;            
                \While{$\bar{i} >0$}{\label{line:while:h}
                    $(i^*,k^*)\gets \arg\min_{i'\in\mathscr{B}_{\bar{i},j},k'\in[h]}d^h_{\varpi^{k'}_{i',j},j}$\label{line:min:h}\;
                    \tcc*{If tie case, select the smallest $i'$ and arbitrary $k'$. }
                    $\varpi \gets \varpi^{k^*}_{i^*,j}$ \label{line:select_neighbour}\; 
                    $a^{\text{MAI},k^*}_{i^*,\bar{i},j}(\pmb{s}^h,t)\gets 1$\;
                    $a^{\text{MAI},k^*}_{i^*,\varpi,j}(\pmb{s}^h,t)\gets 0$
                    \tcc*{Replace movement $(i^*,k^*,\varpi,j)$ with $(i^*,k^*,\bar{i},j)$}
                    \uIf {$d^h_{\varpi,j}=0$}{ \label{line:reach_edge:h}
                        Find an $(i',k')\in\mathscr{B}_{\varpi,j}\times [h]$ such that $\varpi^{k'}_{i',j}= \infty$\;       
                        $a^{\text{MAI},k'}_{i',\varpi,j}(\pmb{s}^h,t)\gets 1$
                        \tcc*{Take movement $(i',k',\varpi,j)$}
                        $\varpi^{k'}_{i',j}\gets \varpi$\;
                        $\varpi^{k^*}_{i^*,j}\gets \bar{i}$\;
                        $\bar{i}\gets 0$\;
                    }\Else{
                        Exchange the values of variables $\varpi^{k^*}_{i^*,j}$ and $\bar{i}$\;
                    }
                \label{line:endWhile:h}}
                Based on the updated $\pmb{\varpi}^h$, update the values of $d^h_{i',j}(\pmb{s}^h,\pmb{\varpi}^h)$ for all $i'\in[I]$ that satisfy \eqref{eqn:distance:h}.\;
                $d^h_{i',j}\gets d^h_{i'j}(\pmb{s}^h,\pmb{\varpi}^h)$ for all $i'\in[I]$\;
            }
        }        
    }
    \Return $a^{\text{MAI},k}_{i,i',j}(\pmb{s}^h,t)$ for all $(i,k,i',j)\in\bigcup_{j\in[J]}\mathscr{M}^{\text{MAI},h}_j(t)$\;
}
\caption{Pseudo-code for the movement adaption in the scaled system.}\label{algo:movement_adaption:h}
\end{algorithm*}
\DecMargin{1em}

\section{Topology of the Three Cases}\label{app:simulation:topology}

The topology for Cases~I ($I=6$), Case~II ($I=10$), and Case~III ($I=14$) are plotted in Figure~\ref{fig:topology}.

\begin{table*}[t]\color{black}
\caption{Important Symbols: Random Variables for MAB-ML}\label{table:symbol:random}
\begin{tabular}{p{2cm}p{5cm}p{2cm}p{6cm}}
\hline
$S^{\phi,k}_{i,j}(t)$ & State variable of sub-process $(i,j,k)\in[I]\times[J]\times[h]$. It is a tuple with two parts, $ (K^{\phi,k}_{i,j}(t), I^{\phi,k}_{i,j}(t))$, 
where the former is the knowledge variable and the latter is the indication variable. When the superscript $k$ is omitted, it indicates the special case with $k=h=1$.&
$\bm{S}^{\phi,h}(t)$ & State vector of MAB-ML scaled by scaling parameter $h$, $\bm{S}^{\phi,h}(t) \coloneqq  (S^{\phi,k}_{i,j}(t): (i,j,k)\in[I]\times[J]\times[h])$. When the superscript $k$ or $h$ is omitted, it indicates the special case with $k=h=1$.\\
$K^{\phi,k}_{i,j}(t) = \mathcal{k}\bigl(S^{\phi,k}_{i,j}(t)\bigr)$& The knowledge variable of sub-process $(i,j,k)\in [I]\times[J]\times[h]$. When the superscript $k$ is omitted, it indicates the special case with $k=h=1$.&
$I^{\phi,k}_{i,j}(t) = \mathcal{i}\bigl(S^{\phi,k}_{i,j}(t)\bigr)$
& The indication variable of sub-process $(i,j,k)\in[I]\times[J]\times[h]$. When the superscript $k$ is omitted, it indicates the special case with $k=h=1$.\\
$\mathcal{e}^{\phi,k}_{i,j}(\bm{S}^{\phi,h}(t),t)$ & Indicator on the sub-process $(i,j,k)\in[I]\times[J]\times[h]$ under policy $\phi$, which is equal to $\mathcal{e}\Bigl(\pmb{a}^{\phi,k}_{i,j}\bigl(\bm{S}^{\phi,h}(t),t\bigr)\Bigr)$. When the superscript $k$ or $h$ is omitted, it indicates the special case with $k=h=1$.\\
\hline
\end{tabular}
\end{table*}

\begin{table*}[t]\color{black}
\caption{Important Symbols: Real Numbers and Vectors for MAB-ML}\label{table:symbol:general}
\begin{tabular}{p{1.6cm}p{6.5cm}p{1.6cm}p{6.5cm}}
\hline
$I\in\mathbb{N}_+$ & Number of areas. &
$J\in\mathbb{N}_+$ & Number of agent types.\\
$T\in\mathbb{N}_+$ & Time horizon. &
$M_j\in\mathbb{N}_+$ & Number of type-$j$ agents and increases in the scaling parameter $h\in\mathbb{N}_+$.\\
$\mathcal{w}_{i,i',j}\in\mathbb{R}_0$ & The weight coefficient in Constraint~\eqref{eqn:constraint:canonical}.&
$h\in\mathbb{N}_+$ & Scaling parameter of MAB-ML.\\
  $a^{\phi,k}_{i,i',j}(\pmb{s},t)$ & Action variable for the $k$th sub-process $\{S^{\phi,k}_{i,j}(t),t\in[T]\}$ associated with $(i,j)\in[I]\times[J]$, which is a function of $\pmb{s}\in\mathscr{S}_{i,j}$ and $t\in[T]$ that takes binary numbers. When the superscript $k$ is omitted, it indicates the special case with $k=h=1$.&
 $\pmb{a}^{\phi,k}_{i,j}(\pmb{s},t)$ & Action vector, $(a^{\phi,k}_{i,i',j}(\pmb{s},t)$, of process $\bigl\{S^{\phi,k}_{i,j}(t),t\in[T]\bigr\}$ when the state vector of the entire MAB-ML $\bm{S}^{\phi,h}(t) = \pmb{s}$. When the superscript $k$ is omitted, it indicates the special case with $k=h=1$.\\
$\mathcal{e}(\pmb{a})$& The probability of taking a non-zero action, given $\pmb{a}\in [0,1]^I$, which is a function $\mathcal{e}(\pmb{a}) = 1-\prod_{i'\in[I]}(1-a_{i'})$ of $\pmb{a}\in[0,1]^I$.&
$c_{i,j}(s,\mathcal{e}(\pmb{a}),t)$ & Instantaneous cost rate of sub-process $\{S^{\phi,k}_{i,j}(t),t\in[T]\}$ for any $k\in[h]$ (or process $\{S^{\phi}_{i,j}(t),t\in[T]\}$ for the special case with $k=h=1$) at time $t$, given $S^{\phi,k}_{i,j}(t)=s$ and $\pmb{a}^{\phi,k}_{i,j}(\cdot,t) = \pmb{a}$. The superscript $k$ is omitted in the special case with $k=h=1$.\\
$\bar{a}^{\phi}_{i,i',j}(\pmb{s},t)$ & The probability of taking $a^{\phi,k}_{i,i',j}(\pmb{s},t)=1$ under policy $\phi$ for any sub-process associated with $(i,j)\in[I]\times[J]$. The superscript $k$ is omitted in the special case with $k=h=1$.&
$\bar{\pmb{a}}^{\phi}_{i,j}(\pmb{s},t)$ & The vector of $\bar{a}^{\phi}_{i,i',j}(\pmb{s},t)$ for all $i'\in[I]$. \\
$\alpha^{\phi,k}_{i,i',j}(s,t)$ & Action variable for sub-process $\Bigl\{S^{\phi,k}_{i,j}(t), t\in[T]\Bigr\}$ of the relaxed problem.
It is the probability of taking action $(i,k,i',j)\in [I]^2\times[J]\times[h]$ at time $t$ under policy $\phi$, given $S^{\phi,k}_{i,j}(t) = s$. The superscript $k$ is omitted in the special case with $k=h=1$. For this special case, its definition is given in \eqref{eqn:define_alpha}.& 
$\pmb{\gamma}\in\mathbb{R}^{IJT}$ & Lagrange multipliers for Constraint~\eqref{eqn:constraint:neighborhood:relax} (or Constraint~\eqref{eqn:constraint:neighbourhood:h:relax} for the scaled case).\\
$L^h(\pmb{\gamma})$ & Dual function of the relaxed problem, given Lagrange multipliers $\pmb{\gamma}\in\mathbb{R}^{IJT}$. When $h$ is omitted, it indicates the special case with $h=1$.&
$V^{\pmb{\gamma}_{i,j}}_{i,j}(s,t)$& Value function for process $\Bigl\{S^{\phi,k}_{i,j}(t),t\in[T]\Bigr\}$ with cost rate $C^{\pmb{\gamma}_{i,j}}_{i,j}(s,\pmb{a},t)$ defined in \eqref{eqn:define_c_gamma}. The superscript $k$ can take any value in $[h]$ and is omitted for the special case with $k=h=1$.  \\
$\pi^0_{i,j}(s)$ & Initial probability for $S^{\phi,k}_{i,j}(1)=s$, where the superscript $k$ can be any value in $[h]$ and is omitted for the special case with $k=h=1$.&
$\pmb{\pi}^0_{i,j}$ & Vector $(\pi^0_{i,j}(s):s\in\mathscr{S}_{i,j})$.\\
$V_{i,j}(s,t)$ & Control variables for linear optimization \eqref{eqn:dual_problem:2:1}-\eqref{eqn:dual_problem:2:3}. &
$V^*_{i,j}(s,t)$ & Optimal control variables for linear optimization \eqref{eqn:dual_problem:2:1}-\eqref{eqn:dual_problem:2:3}. \\
$(\pmb{\gamma}^*,\pmb{V}^*)$ & Optimal solution to linear optimization \eqref{eqn:dual_problem:2:1}-\eqref{eqn:dual_problem:2:3}, where $\pmb{\gamma}^*=(\gamma^*_{i,j,t}: (i,j)\in[I]\times[J], t\in[T])$ and $\pmb{V}^* = (V^*_{i,j}(s,t): (i,j)\in[I]\times[j], s\in\mathscr{S}_{i,j},t\in[T])$. &
$\theta_{i,i',j,t}(\gamma)$& A function defined in \eqref{eqn:define_theta} to normalize $\gamma\in\mathbb{R}$ with $\mathcal{w}_{i,i',j}$. \\
$\vartheta^{\pmb{\gamma}_{i,j}}_{i,j}(s,t)$ & The marginal cost of state $s\in\mathscr{S}_{i,j}$ at time $t$, given the multipliers $\pmb{\gamma}_{i,j}\in \mathbb{R}^{\lvert\{i\}\cup \bar{\mathscr{B}}^{-1}_{i,j}\rvert T}$.&
$\eta_{i,i',j}(s,t)$ & The index assigned to the action $(i,i',j)$ when the state $S^{\phi,k}_{i,j}(t) = s$. The superscript $k$ can be any value in $[h]$ and is omitted for the special case with $k=h=1$.\\ 
$\pmb{\eta}$ & The vector of indices $(\eta_{i,i',j}(s,t): (i,i',j)\in[I]^2\times[J], s\in\mathscr{S}_{i,j}, t\in[T])$.&
$\mathcal{r}^{\phi}_j(i,i',t)$ & The rank of the action $(i,i',j)$ according to the \partialref{define:movement_ranking}{action (movement) ranking}.\\
$\Gamma^{h,\phi}$ & The expected cumulative costs of the entire system normalized by $h$, which is defined in \eqref{eqn:define_Gamma}.&
$\Gamma^{h,*}$ & The minimum of the relaxed problem given scaling parameter $h$.\\
\hline
\end{tabular}
\end{table*}

\begin{table*}[t]\color{black}
\caption{Other Important Symbols for MAB-ML}\label{table:symbol:sets}
\begin{tabular}{p{1.6cm}p{6.5cm}p{1.6cm}p{6.5cm}}
\hline
\multicolumn{4}{l}{Special Policy}\\
\hline
$\psi^*$ & The abbreviation of the mixed policy $\psi(M,\pmb{\phi}^*,\pmb{\pi}^*)$, where $\pmb{\phi}^*=(\phi^*_m: m\in[M])$ and $\pmb{\pi}^*=(\pi^*_m: m\in[M])\in[0,1]^M$, $\phi^*_1,\phi^*_2,\ldots,\phi^*_M\in\tilde{\Phi}^h$ satisfy \eqref{eqn:lemma:indexability} and \eqref{eqn:vartheta} with $\pmb{\gamma}$ and $\phi_{i,j}(\pmb{\gamma}_{i,j})$ replaced by $\pmb{\gamma}^*$ and $\phi^*_m$, respectively, and \eqref{eqn:prop:strong_duality:6} is satisfied with such $\pmb{\phi}^*$ and $\pmb{\pi}^*$.\\
\hline
\multicolumn{4}{l}{Important Labels}\\
\hline
$i \in [I]$ & Label of an area. &
$j \in [J]$ & Label of an agent type. \\
$(i,j,k)\in[I]\times[J]\times[h]$ & Label of process $\Bigl\{S^{\phi,k}_{i,j}(t), t\in[T]\Bigr\}$. When $k$ is omitted, it indicates the special case with $k=h=1$.&
$(i,k,i',j)\in[I]^2\times[J]\times[h]$  & Label of an action (movement) for process $\bigl\{S^{\phi,k}_{i,j}(t),t\in[T]\bigr\}$. When $k$ is omitted, it indicates the special case with $k=h=1$. \\
$t \in [T]$ & Label of a time slot. &
$\phi\in\tilde{\Phi}^h$ & Label of a policy. When $h$ is omitted, it indicates the special case with $h=1$. \\
$s \in \mathscr{S}_{i,j}$  & Label of a state of process $\Bigl\{S^{\phi,k}_{i,j}(t)\Bigr\}$. When $k$ is omitted, it indicates the special case with $k=h=1$.&
$\pmb{\gamma}_{i,j}$ & The vector $(\gamma_{i',j,t}:i'\in\{i\}\cup\bar{\mathscr{B}}^{-1}_{i,j},t\in[T])$ with $\gamma_{i',j,t}$ the corresponding elements of the Lagrange multipliers $\pmb{\gamma}$.\\
\hline
\multicolumn{4}{l}{Sets}\\
\hline
$\mathscr{S}^h$ & The state space of process $\bigl\{\bm{S}^h(t), t \in [T]\bigr\}$. The superscript $h$ is omitted in the special case with $h=1$. & 
$\mathscr{S}_{i,j}$ & The state space of process $\bigl\{S^{\phi,k}_{i,j}(t),t\in[T]\bigr\}$ for any $k\in[h]$. The superscript $k$ is omitted in the special case with $k=h=1$.\\
$\mathscr{B}_{i,j}$ & A subset of $[I]$, indicating $|\mathscr{B}_{i,j}|$ different actions associated with an element in $[I]$. For the patrol case, it represents the neighbourhood of area $i$ for agent type $j$. &
$\bar{\mathscr{B}}_{i,j}$ & A subset of $[I]$ used in \eqref{eqn:constraint:canonical} to establish the connections between different (multi-action) bandit processes. For the patrol case, $\bar{\mathscr{B}}_{i,j} = \mathscr{B}_{i,j}$.\\
$\bar{\mathscr{B}}^{-1}_{i,j}$ & A subset of $[I]$, defined as $\bar{\mathscr{B}}^{-1}_{i,j}=\{i'\in[I] | i \in \bar{\mathscr{B}}_{i',j}\}$. For the patrol case, $\bar{\mathscr{B}}^{-1}_{i,j} = \mathscr{B}_{i,j}$.&
$\mathscr{A}_{i,j}$ &  The set of available actions for the relaxed problem, which is subject to Constraint~\eqref{eqn:constraint:exclusive} (or equivalently, Constraint~\eqref{eqn:constraint:exclusive:h} for scaling parameter $h >1$). The set is defined in \eqref{eqn:define_A}.\\
$\Phi^h$ & The set of all the polices $\phi$ determined by action variables $\pmb{a}^{\phi,h}(\pmb{s},t)$ for all $\pmb{s}\in\mathscr{S}^h$ and $t\in[T]$. The superscript $h$ is omitted in the special case with $h=1$. &
$\tilde{\Phi}^h$ & The set of all the polices $\phi$ determined by action variables $\alpha^{\phi,k}_{i,i',j}(s,t)$ for all $(i,k,i',j)\in[I]^2\times[J]\times[h]$, $t\in[T]$ and $s\in\mathscr{S}_{i,j}$ with \eqref{eqn:constraint:exclusive:h} satisfied. The superscript $k$ or $h$ is omitted in the special case with $k=h=1$. \\ 
$\Phi^h$
& the set of all the polices for the scaled general MAB-ML system \\
$\mathscr{M}^{\phi,h}_j(t)$ & A set of all the actions $(i,k,i',j)$ such that $\sum\nolimits_{k'\in[h]}g_{i',j}(S^{\phi,k'}_{i',j}(t))>0$ and $i'\in\bar{\mathscr{B}}_{i,j}$. The superscript $h$ is the scaling parameter and is omitted for the special case with $h=1$.& 
$\mathscr{I}^{\phi,h}_j(t)$ & The set of all $i\in[I]$ such that  $\sum_{k\in[h]}g_{i,j}(S^{\phi,k}_{i,j}(t))>0$ but \eqref{eqn:constraint:neighbourhood:h} is not satisfied. The superscript $h$ is the scaling parameter and is omitted for the special case with $h=1$.\\

\hline
\end{tabular}
\end{table*}

\begin{table*}[t]\color{black}
\caption{Important Symbols: Special Symbols for the patrol problem}\label{table:symbol:patrol}
\begin{tabular}{p{1.6cm}p{6.5cm}p{1.6cm}p{6.5cm}}
\hline
$\varpi^k_{i,j}(\pmb{a})$ & Original location of a type-$j$ agent that is decided to be moved to sub-area $(i,k)$ according to action vector $\pmb{a}$. The superscript $k\in[h]$ is omitted in the specical case with $k=h=1$.&
$\pmb{\varpi}^h(\pmb{a})$ & The vector of $\varpi^k_{i,j}(\pmb{a})$  for all $(i,j,k)\in[i]\times[J]\times[h]$ according to action vector $\pmb{a}$. The superscript $h$ is the scaling parameter and is omitted in the special case with $h=1$.\\
$d^h_{i,j}(\pmb{s},\pmb{\varpi})$ & The distance of an area $i$ for type-$j$ agent with respect to state vector $\bm{S}^{\phi,h}(t)=\pmb{s}$ and $\pmb{\varpi}^h\bigl(\pmb{a}^{\phi,h}(\pmb{s},t)\bigr) = \pmb{\varpi}$. The superscript $h$ is the scaling parameter and is omitted in the special case with $h=1$.\\
\hline
\end{tabular}
\end{table*}

\section{Settings of Initial States and Agents' Positions}\label{app:simulation:initialization:knowledge_position}

The initial values of $\alpha^{\phi,k}_{i,j}(1)=\alpha^0_{i,j}$ and $\beta^{\phi,k}_{i,j}(t)=\beta^0_{i,j}$ are adapted from the crime rates in United Kingdom, reported in 2022~\cite{UKCrime}.
For  $(i,j)\in[I]\times[J]$, the values of such $\alpha^0_{i,j}$ and $\beta^0_{i,j}$ are provided in Table~\ref{table:alpha_beta}.

Given agent type $j\in[J]$, apart from the determined initial values of $\alpha^{\phi,k}_{i,j}(1)=\alpha^0_{i,j}$ and $\beta^{\phi,k}_{i,j}(t)=\beta^0_{i,j}$, for the initial probabilities $\pi^0_{i,j}(s)$ for all $i\in[I]$ and $s\in\mathscr{S}_{i,j}$, it remains to decide the probabilities of $I^{\phi,k}_{i,j}(1)=0$ or $1$.

As described in Section~\ref{sec:simulation}, for type $j\in[J]$, the number of agents is randomly generated with expectation $M^0_jh$ where $M^0_j\in\mathbb{N}_+$ is also randomly generated from $\{1,2,\ldots, I/2\}$.
Given an instance of $M^0_j$, we consider in total $20M^0_j$ type-$j$ \emph{virtual agents} that are randomly distributed in the $I$ areas. Note that these $20M^0_j$ virtual agents are not any real agent, but used for constructing the initial probabilities $\pi^0_{i,j}(s)$ for all $i\in[I]$ and $s\in\mathscr{S}_{i,j}$. 
More precisely, initialize a set $\mathscr{I}=[I]$, for the first virtual agent, we uniformly randomly pick an area from $\mathscr{I}$ and locate this virtual agent there. Then, repeat the process for the second virtual agent until the last one. During this process, if area $i\in\mathscr{I}$ locates $20$ virtual agents, then remove it from $\mathscr{I}$ for all later virtual agents. 
Let $M_{i,j}$ represent the number of virtual agents in area $i$.
The initial probability for area $i$ to have a type-$j$ agent is set to be $M_{i,j}/20$.
In this context, for $(i,j)\in[I]\times[J]$ and $s\in\mathscr{S}_{i,j}$, if $\mathcal{k}(s) = (\alpha^0_{i,j},\beta^0_{i,j})$ and $\mathcal{i}(s)=1$, then the initial probability $\pi^0_{i,j}(s) = M_{i,j}/20$; if $\mathcal{k}(s) = (\alpha^0_{i,j},\beta^0_{i,j})$ and $\mathcal{i}(s)=0$, then $\pi^0_{i,j}(s) = 1-M_{i,j}/20$; otherwise, $\pi^0_{i,j}(s) = 0$.

In a scaled patrol system with scaling parameter $h\in\mathbb{N}_+$, for any $k\in[h]$, the initial state of a sub-process $\{S^{\phi,k}_{i,j}(t), t\in[T]\}$ follows the probability distribution $\pmb{\pi}^0_{i,j} = (\pi^0_{i,j}(s): s\in\mathscr{S}_{i,j})$.
With given $\pmb{\pi}^0_{i,j}$, we randomly generate the initial value of $S^{\phi,k}_{i,j}(1)$.
If $I^{\phi,k}_{i,j}(1) = \mathcal{i}\bigl(S^{\phi,k}_{i,j}(1)\bigr) = 1$, then there is a type-$j$ agent in sub-area $(i,k)$; otherwise, no type-$j$ agent there.

\section{Detailed Settings of the Sample Systems}\label{app:simulation:case_study}
For the simulations results plotted in Figure~\ref{fig:case_study}, for $(i,j)\in[I]\times[J]$ and $s\in\mathscr{S}_{i,j}$, if $\mathcal{k}(s)\neq(\alpha^0_{i,j},\beta^0_{i,j})$, then the initial probability $\pi^0_{i,j}(s)=0$; otherwise, the initial probability is provided in Table~\ref{table:init_probs}.
Here, $(\alpha^0_{i,j},\beta^0_{i,j})$ are those initial values provided in Table~\ref{table:alpha_beta}.

\section{Notation Tables}\label{app:notation}
We summarize important variables for general MAB-ML in Tables~\ref{table:symbol:random}, \ref{table:symbol:general}, and \ref{table:symbol:sets}, and some special symbols used for only the patrol case in Table~\ref{table:symbol:patrol}.

\clearpage

\bibliographystyle{references/IEEEtran}
\bibliography{references/IEEEabrv,references/bib1}

\end{document}